\documentclass[12pt]{amsart}
\usepackage{amsmath, amssymb, amscd, array}

\usepackage{tabularx,mdwtab}

\usepackage{lscape}

\usepackage{hyperref}

\usepackage{enumitem}

\usepackage{stmaryrd}

\usepackage{DLdef1}

\newcommand{\Size}{\text{Size}}

\newcommand {\fBj} {{\mathfrak{Bj}}}
\newcommand {\fbr} {{\mathfrak{br}}}
\newcommand {\fer} {{\mathfrak{er}}}

\newcommand {\ffr} {{\mathfrak{fr}}}
\newcommand {\fy} {{\mathfrak{y}}}
\newcommand {\fdy} {{\mathfrak{dy}}}
\newcommand {\fby} {{\mathfrak{by}}}

\newcommand {\fme} {{\mathfrak{me}}}
\newcommand {\fMe} {{\mathfrak{Me}}}
\newcommand {\fmy} {{\mathfrak{my}}}
\newcommand {\fsmy} {{\mathfrak{smy}}}

\newcommand{\un}{\underline{N}}
\newcommand{\del}{\partial}
\newcommand {\der} {{\partial}}
\newcommand {\lra} {\leftrightarrow}
\newcommand {\Lra} {\Longrightarrow}

\newcommand{\what}{\widehat}

\newcommand{\Cinf}{C^\infty}

\def\mmat #1,#2,#3,#4,{\text{\small\arraycolsep=3pt $
\begin{pmatrix}#1&#2\\#3&#4\end{pmatrix}$}}

\let\ssec\subsection

\renewcommand {\ssbegin}[2][*]
 {\refstepcounter{subsection}%
\if#1*
\addcontentsline{toc}{subsection}{\thesubsection.\hskip 1pc #2}%
\else
\addcontentsline{toc}{subsection}{\thesubsection.\hskip 1pc #2. #1}%
\fi
 \def \secno {\gdef \secno {}{\ssecfont
\thesubsection.\hskip 2ex}%
 }%
 \begin{#2}}

\renewcommand {\sssbegin}[2][*]
 {\refstepcounter{subsubsection}
\if#1*
\addcontentsline{toc}{subsubsection}{\thesubsubsection.\hskip 1pc #2}%
\else
\addcontentsline{toc}{subsubsection}{\thesubsubsection.\hskip 1pc #2. #1}
\fi
 \def \secno {\gdef \secno {}{\ssecfont \thesubsubsection.\hskip 2ex}%
 }%
 \begin{#2}}

\renewcommand {\parbegin}[2][*]
 {\refstepcounter{paragraph}
\if#1*
\addcontentsline{toc}{paragraph}{\theparagraph.\hskip 1pc #2}%
\else
\addcontentsline{toc}{paragraph}{\theparagraph.\hskip 1pc #2. #1}
\fi
 \def \secno {\gdef \secno {}{\ssecfont \theparagraph.\hskip 2ex}%
 }%
 \begin{#2}}

 \makeatletter
\@namedef{subjclassname@2020}{\textup{2020} Mathematics Subject Classification}
\makeatother

\begin{document}

\title[Non-integrable distributions]{Non-integrable distributions with simple infinite-dimensional Lie (super)algebras of symmetries}

\author[A.~Krutov, D.~Leites, I.~Shchepochkina]{Andrey Krutov${}^{a,*}$, Dimitry Leites${}^{b}$, Irina~Shchepochkina${}^c$}
\address{
 ${}^*$The corresponding author\\
 ${}^a$Mathematical Institute of Charles University, Sokolovsk\'a 83, Prague, Czech Republic \\
 \email{andrey.krutov@matfyz.cuni.cz}\\
${}^b$Department of Mathematics\\
Stockholm University, Albanov\"agen 28, SE-114 19 
Stockholm, Sweden\\
\email{dimleites@gmail.com}\\
${}^c$Independent University of Moscow,
B. Vlasievsky per., d. 11, RU-119 002 Moscow, Russia\\
\email{irina@mccme.ru}\\
}

\begin{abstract}Under usual locality assumptions, we classify all non-integrable distributions with simple infinite-dimensional Lie superalgebra of symmetries
over~$\mathbb{C}$: 
we single out 15 series (containing 2 analogs of contact series and one family of deformations of their divergence-free subalgebras), and 7 exceptional Lie superalgebras. 
Over algebraically closed fields~$\mathbb{K}$ of characteristic $p>0$, we classify the W-gradings (corresponding to a maximal subalgebra of finite codimension) of the known simple vectorial Lie (super)algebras with
unconstrained shearing vector of heights of the indeterminates, distinguish W-gradings of (super)algebras preserving non-integrable distributions. For $p>3$, we get analogs of the result over $\mathbb{C}$. For $p=3$, of all possible W-gradings (12 of Skryabin algebras, 3 of superized Melikyan algebras, and 4 of Bouarroudj superalgebras) most are new, together with the corresponding distributions. For $p=2$, we also get several new examples of distributions and their Lie (super)algebras of symmetries. 
\end{abstract}

\date{}



\maketitle

\markboth{\itshape Andrey Krutov\textup{,} Dimitry Leites\textup{,} Irina Shchepochkina}
{{\itshape Non-integrable distributions}}

\thispagestyle{empty}

\section{Introduction}

Killing and \'E.~Cartan described the finite-dimensional simple Lie algebras, especially exceptional ones, not
in customary now terms of root systems (not yet defined at that time), but in terms of non-integrable distributions each of these algebras preserves.
For infinite-dimensional Lie (super)algebras, the situation is different:
Lie algebras of contract vector fields on odd-dimensional spaces are the only simple infinite-dimensional Lie algebras (of polynomial growth) over~$\Cee$ and $\Ree$ that
preserve non-integrable distributions. 
In this paper, under the usual locality and regularity assumptions, we clasify non-integrable distributions with simple
infinite-dimensional Lie superalgebras of symmetries over $\Cee$ and over algebraically closed fields $\Kee$
of characteristic $p>0$.
The picture over $\Kee$ and in the super case over any field is much richer than non-super over $\Cee$.

The protagonists of this paper are infinite-dimensional filtered Lie (super)algebras~$\cL$ or the associated
with them graded Lie (super)algebras~$\fg$, where $\fg_i=\cL_i/\cL_{i+1}$ is finite-dimensional for every~$i$. We call these algebras \textit{vectorial}
since their elements are vector fields on the linear (super)manifold corresponding (as described, for example, in~\cite{Lo}) to the linear (super)space
$(\fg/\fg_{\geq0})^*$, where $\fg_{\geq0}:=\oplus_{i\geq 0}\ \fg_i$ (resp., $\cL_0$) is a maximal subalgebra
of finite codimension in~$\fg$ (resp.,~$\cL$).
Such filtrations (and gradings) arise from the corresponding distribution~$\cD$ on the
(super)manifold~$\cM$ preserved
by~$\cL$ (and~$\fg$); recall that a~\textit{distribution} is a~subbundle~$\cD$ of the tangent bundle~$T\cM$ over ~$\cM$.

The definitions of this subsection have direct superizations and generalizations for algebraically closed
fields~$\Kee$ of positive characteristics
by replacing the algebra~$\cF$ of smooth functions on a~manifold by the algebra of polynomials, or (tacitly)
of formal power series, or of divided powers;
see~\cite{MaG,SoS,GL}.
To simplify notation, we formulate the definitions in non-super setting 
over $\Cee$.
Let $M$ be an $m$-dimensional manifold. Since we consider only local problems, one can assume that $M=\Kee^m$ or~$\Cee^m$.
The flag associated to~$\cD$ 
is described by its spaces of sections, i.e., the space $\Gamma\cD_{-i}$ of sections of the bundle~$\cD_{-i}$ is spanned by $\leq i$ brackets of vector fields of $\Gamma\cD:=\Gamma\cD_{-1}$:
\be\label{deja}
 \cD = \cD_{-1} \subset \cD_{-2} \subset \ldots, \text{~~where 
 $\Gamma \cD_{-k-1} := [\Gamma \cD, \Gamma \cD_{-k}]+\Gamma \cD_{-k}$}.
\ee
The distribution~$\cD$ is called \textit{regular} if all the dimensions~$m_i:=\dim\cD_{-i}$ are constants on~$M$. We consider only \textit{complete} regular distributions, i.e., such that the sections of~$\cD$ generate by bracketing the whole Lie algebra
$\fvect(m)$ of sections of~$TM$.

The \textit{integral manifold} at point $*\in M$ for the~distribution 
$\cD$ is a~submanifold 
$N\subset M$ such that 
$T_{*}N=\cD_{*}$ for every point $*\in N\subset M$. A distribution is called \textit{integrable} at point $*$ if there is an integral manifold through 
$*\in M$. Hertz coined the term \textit{non-holonomous} for any manifold given with a~non-integrable distribution, see \cite[p.~80]{H}, the modern spelling is \textit{non-holonomic};
the non-integrable distributions are also called \textit{non-holonomic}.

Let the \textit{symbol algebra}~$\gr(TM)$ 
corresponding to $\cD$ be
\[
\gr(TM) := \oplus_{i>0} \gr\cD_{-i}, \text{where $\Gamma (\gr \cD_{-i}) := \Gamma \cD_{-i} / \Gamma \cD_{-i+1}$.}
\]

We assume that 
\begin{equation}\label{assum} 
\begin{minipage}[l]{14cm}
the stalk of the sheaf of sections of $\cD$ at every point $*$ of the manifold $M$ is a~$\Zee$-graded Lie algebra $\fg(*)_-=\oplus_{-d\leq i\leq -1}\ \fg_i(*)$ such that $\dim\fg(*)_-=\dim M$ and $\fg_{-1}(*)=\cD_*$; moreover, all Lie algebras $\fg(*)$ are isomorphic to a~Lie algebra $\fg=\oplus_{-d\leq i\leq -1}\ \fg_i$ and the Lie algebra $\fg_0$ --- a~subalgebra of the Lie algebra of degree-preserving automorphisms of $\fg_{-}$ --- irreducibly acts on $\fg_{-1}$. 
\end{minipage}
\end{equation}
If $\dim\fg_i=m_i$, then the $\Zee$-grading of $\fg$ yields a~nonstandard grading of $\fvect(m)\simeq \gr(TM)|_*$, 
and hence of $\cF$, see~\cite{GL,Shch}:
 \begin{equation}\label{GradingMk}
\begin{array}{ll}
 & \deg x_1 = \ldots = \deg x_{m_1} = 1,\\
 & \deg x_{m_1+1} = \ldots = \deg x_{m_2} = 2,\\
 & \dots \dots \dots \dots\\
 & \deg x_{m-m_{d-1}+1} = \ldots = \deg x_m = d,
\end{array}
\end{equation}
where the $x_i$ are coordinates compatible with the flag associated to~$\cD$.
Denote by $\fv = \oplus_{i \geq -d}\ \fv_i$ or $\fvect(m; v_{\cD})$, where $v_{\cD}:=(m_1, m_2, \dots, m)$ is the \textit{growth vector} of~$\cD$, the Lie algebra $\fvect(m)$ with the grading \eqref{GradingMk}.
The complete Cartan prolong of $\fv_{-}:=\oplus_{i <0} \ \fv_i$, see Subsection~\ref{ssTrans}, is called the \emph{symmetry algebra} of~$\cD$.


In~\cite{VG1} Vershik and Gershkovich stated that 
\begin{equation}\label{VG} 
\begin{minipage}[l]{13.3cm}
an open and dense orbit of germs of a~$k$-dimensional distribution in $\Ree^n$ under the action of the group of germs of smooth diffeomorphisms is either (a) a~contact distribution, where $k=n-1$, or\\ (b) the Engel distribution for which the growth vector is $(2, 3, 4)$.
\end{minipage}
\end{equation}
The direct analog of statement~\eqref{VG} holds over $\Cee$ as well.
Any generic orbit of such distribution is of infinite codimension.
It follows that 
the Lie algebra preserving the non-integrable distribution is of infinite dimension only in the two cases (a) and (b) of statement~\eqref{VG}.
The Lie algebra preserving the contact distribution for $n$ odd is simple;
the Lie algebra preserving the Engel distribution and the Lie algebra preserving the contact distribution for $n$ even are not simple. 
No analog of the statement~\eqref{VG} is formulated yet for supermanifolds or for varieties over fields of characteristic $p>0$; what we did is a step towards formulating such analog.

\'E.~Cartan and Killing considered the non-integrable distributions preserved by finite-dimen\-si\-onal Lie
algebras as well as by infinite-dimen\-si\-onal ones. Both types of distributions can be associated with
interesting differential equations on manifolds as well as on supermanifolds, see \cite{C, Shch, DG, DK, KST, KST1, vC, vCS} and references therein.

\subsection{Weisfeiler filtrations and gradings}\label{WeisF} For vectorial Lie
(super)algebras, the notion invariant with respect to automorphisms is filtration, not grading.
Weisfeiler introduced a~most natural filtrations, see~\cite{W}, bearing his name. For any characteristic~$p$, let~$\cL_{0}$ be a~maximal subalgebra of finite codimension in the vectorial Lie (super)algebra $\cL$, let~$\cL_{-1}$ be a~minimal $\cL_{0}$-invariant subspace strictly
containing~$\cL_{0}$. Then, for $i\geq 1$, we set, see \cite{BGLLS2}:
\begin{equation*}
\label{WeiGr}
\begin{array}{l}
\cL_{-i-1}=\begin{cases}
[\cL_{-1}, \cL_{-i}] +\cL_{-i}+\Span(X^2\mid X\in \cL_{-1}),&\text{if $\cL$ is a~superalgebra, $p=2$ and $i=1$},\\
[\cL_{-1}, \cL_{-i}]+\cL_{-i},&\text{otherwise};\\
\end{cases}\\
\cL_i =\{D\in \cL_{i-1}\mid [D, \cL_{-1}]\subset\cL_{i-1}\}.
\end{array}
\end{equation*}
We thus get a~filtration called, after~\cite{W}, a~\textit{Weisfeiler filtration} with strict inclusions:
\begin{equation}
\label{Wfilt} \cL= \cL_{-d}\supset \cL_{-d+1}\supset \dots \supset
\cL_{0}\supset \cL_{1}\supset \dots
\end{equation}
The number~$d$ in formula~\eqref{Wfilt} is called the \textit{depth} of~$\cL$, and
of the associated \textit{Weisfeiler-graded} Lie superalgebra
$\fg={\oplus}_{-d\leq i}\ \fg_i$, where $\fg_{i}=\cL_{i}/\cL_{i+1}$.
For brevity, we will say \textit{W-graded} and \textit{W-filtered}.

It is clear that $W$-graded vectorial Lie (super)algebras whose depth is $d>1$ preserve non-integrable distributions. For more details on relation between vectorial algebras and distributions, see Subsection~\ref{ssVar}.

The Weisfeiler filtrations were used in~\cite{LS, LS1} in the classification of the simple infinite-dimensional vectorial
Lie (super)algebras~$\cL$ with polynomial or formal coefficients over
$\Cee$; their incarnations over $\Kee$ of characteristic $p>0$ are especially interesting for $p=2$, see~\cite{BGLLS1, BGLLS2}. 

Therefore, to obtain the list of non-integrable distributions (satisfying condition~\eqref{assum})
with infinite-dimensional simple Lie (super)algebra of symmetries
we consider Weisfeiler filtrations of the following Lie (super)algebras:

(a) the simple
 infinite-dimensional vectorial Lie (super)algebras whose coefficients are polynomials over $\Cee$; for their classification, see \cite{LS, LS1}; 
 
 (b) simple
infinite-dimensional vectorial Lie (super)algebras whose coefficients are divided powers (obtained from known simple
infinite-dimensional vectorial Lie (super)algebras by taking unconstrained shearing vector of heights of the indeterminates, see eq.~\eqref{shear}) over an algebraically closed field~$\Kee$ of characteristic $p>0$. The classification of simple finite-dimensional Lie algebras for $p=3$ and $2$ is conjectural (this conjecture is vague and incomplete for $p=2$), together with
the classification of simple finite-dimensional Lie superalgebras for any $p>5$, which is also conjectural (together with an unfathomable superization of~\cite{Sk1}) but most probably correct, see~\cite{BLLS2, BGLLS1, BGLLS2}.

Among these Lie (super)algebras considered with Weisfeiler filtration, we select the ones which preserve non-integrable distributions. Sometimes, we naturally widen the list of thus selected examples to include certain \lq\lq relatives" of simple Lie (super)algebras: the algebras of their derivations. 

If $p>0$, our examples are obtained from the known finite-dimensional Lie (super)algebras by considering the unconstrained shearing
vector of heights of the indeterminates (for definitions, see \S\ref{Sprel}). 
Recall that the classification of simple finite-dimensional Lie algebras is completed for $p>3$,
see~\cite{S, BGP} appended with explicit Poisson brackets corresponding to different normal symplectic forms in \cite{Sk1, Kir}. For basics on Lie superalgebras, see \cite{BGLLS2, Lo}.

\ssec{Remark}\label{Rem1}
Over~$\Cee$, the first example of a~simple Lie superalgebra~$\fg$ preserving a~non-integrable distribution and containing
an infinite-dimensional simple Lie superalgebra~$\fh$ also preserving the same distribution, together with something else
which is taken care of by certain constraints of positive components of~$\fh$, was observed in 
\cite{Sh5,Sh14}, where the exceptional simple superalgebra $\fkas\subset \fk(1|6)$ was discovered. 

The same phenomenon is demonstrated by many new examples found in this paper and the following well-known pairs of simple
modular Lie algebras described, e.g., in~\cite{GL1}:
\[ 
\begin{array}{l}
\fk(5;\un)\supset \fme(\un) \text{~~for $\un=(1, 1,1, N_4, N_5)$},\ \ \
\fk(3; (n, 1,1) )\supset \ffr(n).\\
\end{array}
\]

\ssec{Several Weisfeiler filtrations of the same algebra}\label{SeveralW}
One phenomenon of superization was astounding when first proclaimed in \cite{ALSh}: one simple vectorial Lie superalgebra can have several (finitely many) non-isomorphic W-regradings. For a~proof of completeness of their list over~$\Cee$, see~\cite{LS, LS1} correcting one omission in \cite{Sh5, Sh14}.

In the classification of simple vectorial Lie superalgebras, a~new phenomenon existing for any~$p$ was observed, see \cite{Sh5,Sh14}:
a~simple Lie (super)algebra~$\fg$ preserving a~non-integrable distribution might contain a~smaller simple Lie (super)algebra~$\fh$ also preserving the same distribution together with something else which is taken care of
by positive components of~$\fh$, see Subsection~\ref{ssTrans}. This new construction was called \textit{partial prolongation} in \cite{Shch} to distinguish it from \textit{complete prolongation} considered by \'E.~Cartan and Killing and rediscovered later when their works were successfully forgotten.
In particular, there are many simple Lie (super)algebras (partial prolongs) which are subalgebras of the contact Lie (super)algebras with the same components of negative degrees in their respective $\Zee$-gradings, see Remark~\ref{Rem1}; we list such examples, labeled in tables by ``C", as byproducts of our proofs.

\ssec{On gradings}\label{gradings}
One Lie (super)algebra can have several W-gradings. This phenomenon was encountered in the 1960s in the study of simple modular Lie algebras, but did not draw much attention in the finite-dimensional setting where the importance of W-gradings among other gradings is not manifest.

For reducing a~description of all possible gradings of a~given Lie algebra $\fg$ by abelian groups without
$p$-torsion to the description of the automorphisms of $\fg$, see~\cite{KfiD} with corrections and
generalizations in~\cite{Sk91, Sk95}.
For Lie algebras, Skryabin described the group scheme of automorphisms.
This description implied that any grading of a~given vectorial Lie algebra $\fg$, at least for $p>3$,
corresponds to a~grading of the algebra of truncated polynomials; in particular, the differential forms that
single out a~distribution $\fg$ preserves (contact or symplectic in the simplest cases) are homogeneous with
respect to the induced grading.
It is an~important \textbf{Open problem} to extend this Skryabin's result to $p=3$ and $2$ and superize it for
any $p$, especially $p=0$, where it is a~working hypothesis in~\cite{LS, LS1}.

\ssec{Our results}\label{OurRes}
Here, we classify the known simple infinite-dimensional Lie (super)algebras preserving non-holonomic structures
over algebraically closed fields of characteristic~$p>0$ --- a~part of the super version of classification~\eqref{VG}.
In some cases, the Cartan prolong is a~more natural and easier to describe larger algebra than its simple (with no ideals) derived algebra preserving the same non-integrable distribution. A main tool in our descriptions are Weisfeiler gradings. We do not consider filtered deformations of any W-graded Lie (super)algebra $\fg$ because they preserve the same distribution as $\fg$ does.

\textbf{For $p=0$}, from 
the list of simple Lie superalgebras of vector fields with polynomial coefficients and a~Weisfeiler grading (this is the analog of \'E.~Cartan's list of simple Lie algebras of vector fields with polynomial coefficients over $\Cee$, see~\cite{K1, K2, ALSh, Sh5, Sh14, LS} and numerous corrections listed in \cite{LS1}), we select Lie (super)algebras preserving a~non-integrable distribution. They constitute~2 series similar to the contact Lie algebras, see formulas~\eqref{grVecK} and~\eqref{grVecM}, and~7 exceptions, see tables~\eqref{grVex}.

For $p>5$, the direct versions of the same Lie (super)algebras we considered over $\Cee$, with the algebras of
polynomials and formal series replaced with the algebras of divided powers, provide us with examples of Lie
(super)algebras preserving non-integrable distributions with the same growth vector as over $\Cee$. It is an
important \textbf{Open problem}, with its answer not even conjectured, to extend Skryabin's result on normal
shapes of symplectic forms (see \cite{Sk1}) to non-alternating forms and to extend it to both symplectic and
periplectic cases in the super setting; then, one will find more examples of W-graded Lie superalgebras of
type $\fh$ and $\fle$ than there are in table~\eqref{1114}.

\textbf{For $p=5, 3, 2$}, there appear new examples whereas some analogs disappear or look very differently from their namesakes for $p>5$ and $p=0$.

\textbf{For $p=5$}, we describe all W-gradings of the Melikyan algebra, see Section~\ref{SMel}. 

\textbf{For $p=3$}, we describe all W-gradings of all known simple vectorial Lie algebras (Frank and Ermolaev algebras, see Section~\ref{SFrandEr}, Skryabin algebras, see Section~\ref{SSkr}) and three superizations of the Melikyan algebra, see Section~\ref{SBrj}. (Observe that Skryabin had listed all W-gradings of Frank and Ermolaev algebras, see~\cite{Sk}.) 

\textbf{For $p=2$}, we describe all W-gradings of selected examples of several types, see Section~\ref{Sp=2}.

\textbf{For any $p$}, we describe the linear Lie (super)algebras~$\fg_0$ preserving the non-integrable distributions and the $\fg_0$-modules~$\fg_i$ for~$i$ negative: some of them are exotic.

In Section~\ref{Sequiv}, we discuss several related questions and open problems.

\section{Preliminaries}\label{Sprel} 

The ground field $\Kee$ in this paper is either algebraically closed of characteristic $p>0$ or $\Cee$. The context does not allow to confuse characteristic with parity also denoted by $p$.

\ssec{The Lie
superalgebras}\label{SSlieSuper} We give a~definition of Lie superalgebras working over fields of any characteristic and compare it with the one most often used for $p=0$: superization by means of the Sign Rule. The \textit{Lie superalgebra}
is a~superspace $\fg=\fg_\ev\oplus\fg_\od$ such that the even part
$\fg_\ev$ is a~Lie algebra (in particular, $[x,x]=0$ for any $x\in
\fg_\ev$) the odd part $\fg_\od$ is a
$\fg_\ev$-module (made into the two-sided one by
\textit{anti}-symmetry, i.e., $[y,x]=-[x,y]$ for any $x\in
\fg_\ev$ and $y\in
\fg_\od$), and a~\textit{squaring}
defined on $\fg_\od$ as a~map $S^2(\fg_\od)\tto\fg_\ev;\ \ x\mapsto x^2\in\fg_\ev$ such that 
\begin{equation}\label{squaring}
\begin{array}{l}
\text{$(ax)^2=a^2x^2$ for any $x\in\fg_\od$ and $a\in \Kee$},\\
{}[x,y]:=(x+y)^2-x^2-y^2\text{~is a~bilinear form on $\fg_\od$ with values
in $\fg_\ev$.}\\
\end{array}
\end{equation}

The Jacobi identity involving odd elements splits into two:
\begin{equation}\label{JI}
\begin{array}{l}
~[x^2,y]=[x,[x,y]]\text{~for any~} x\in\fg_\od, y\in\fg_\ev,\\
~[x^2,x]=0\;\text{ for any $x\in\fg_\od$.}
\end{array}
\end{equation}

For $p=2$, we want the space $\fder\ \fg$ of all derivations of $\fg$ be a~Lie superalgebra for any Lie superalgebra $\fg$, so we have to add one more condition 
\[
D(x^{2}) = [D(x),x]\text{~~for all odd elements $x\in\fg$ and any $D\in\fder\ \fg$.}
\]

If $p\neq 2$, then $x^2=\nfrac12 [x,x]$ and conditions \eqref{squaring} are equivalent to bilinearity and symmetry of the bracket of two odd elements; the Jacobi identity for even elements is often written in the symmetric (but meaningless) form 
\[
[x, [y,z]] +\text{cyclic permutations}=0,
\]
or in a~meaningful form reflecting that $\ad_x$ is a~derivation
\[
[x, [y,z]] =[[x, y],z]]+[y, [x,z]] \text{~~for any $x,y,z\in \fg$}.
\]
Any of these two forms of the Jacobi identity is routinely superized via the Sign Rule.

If $p=3$, the Jacobi identity requires one more condition, automatically followed from \eqref{JI} if $p\neq 3$:
\[
[x,[x,x]]=0\text{~for any~} x\in\fg_\od.
\]

By an \textit{ideal} $\fii$ of a~Lie superalgebra $\fg$ we always mean the ideal \textit{homogeneous} with respect to parity, i.e., $\fii={(\fii\cap\fg_\ev)\bigoplus (\fii\cap\fg_\od)}$; for $p=2$, the ideal should be closed with respect to squaring.

A~given Lie (super)algebra $\fg$ is said to be
\textit{simple} if $\dim\fg>1$ and $\fg$ has no proper ideals.

For any Lie \textbf{super}algebra $\fg$, its \textit{derived algebras} are
defined to be (for $i\geq 0$)
\[
\fg^{(0)}: =\fg, \quad
\fg^{(i+1)}=\begin{cases}[\fg^{(i)},\fg^{(i)}]&\text{for $p\neq
2$,}\\
[\fg^{(i)},\fg^{(i)}]+\Span\{g^2\mid g\in\fg^{(i)}_\od\}&\text{for
$p=2$}.\end{cases}
\]

The above were ``naive" definition of Lie superalgebras.
For a~functorial definition, working for any~$p$, see~\cite{Lo} and the Oberwolfach version of~\cite{KLLS}.

\ssec{Transitive Lie superalgebras. Cartan prolongations}\label{ssTrans}
Let $\fg_-=\oplus_{-d\leq i\leq -1}\ \fg_i$ be a~$\Zee$-graded nilpotent Lie (super)algebra, $\fg_0$ a~Lie (super)algebra of grading-preserving derivations of~$\fg_-$ (not necessarily maximal such algebra); additionally, let~$\fg_0$ faithfully act on~$\fg_-$.

Recall that the graded Lie superalgebra 
$\fb=\oplus_{i \geq -d}\ \fb_i$
is said to be \textit{transitive}
if for all $i\geq 0$ we have
\begin{equation*}
\label{12}
\{x\in \fb_i \mid [x,\fb_{-}]=0\}=0, \text{~~where $\fb_{-}:=\oplus _{i <0}\,\fb_i$}.
\end{equation*}

The maximal transitive $\Zee$-graded Lie (super)algebra whose non-positive part is $\fg_-\oplus \fg_0$ is called the \textit{Cartan prolong} of the pair $(\fg_-,\fg_0)$ and is denoted by $(\fg_-,\fg_0)_*$. This prolong is called \textit{complete} if $\fg_0$ is maximal possible in which case the prolong is determined by $\fg_-$. A \textit{partial Cartan prolong} of the pair $(\fg_-,\fg_0)$ is the one determined by a triple $(\fg_-,\fg_0, \tilde\fg_1)$, where $\tilde\fg_1$ is a $\fg_0$-submodule of $\fg_1$ such that $[\fg_{-1}, \tilde\fg_1]=\fg_0$.

\ssec{Basics from supergeometry} For details, see~\cite{MaG,SoS}.
Let $U$ be a~domain in~$\Ree^n$ and $V$ be an~$m$-dimensional vector space. An $n|m$-dimensional \textit{superdomain}~$\cU^{n|m}$
(or just $\cU$) is a~ringed space determined by a~pair $(U,\Cinf(\cU))$, where $\Cinf(\cU) = \Cinf(U)\otimes E^{\bcdot}(V)$ and $E^{\bcdot}(V)$ is the exterior algebra of $V$. Morphisms of such pairs are pairs $(f, \Phi)$, where $f: U\to U'$ is a~diffeomorphism and 
\[
\Phi: \Cinf(U')\otimes E^{\bcdot}(V')\tto \Cinf(U)\otimes E^{\bcdot}(V)
\]
is a~homomorphism of superalgebras (preserving parity and sending 1 to 1).

In what follows, we consider versions of the above with $\Cinf(\cU)$ replaced by $\cF:=\cF(\cU)$ which is the algebra of polynomials, or formal series, or divided powers --- depending on the context.

The elements of the Lie superalgebra $\fder\ \cF$ 
are called
\textit{vector fields}. We denote the left 
$\cF$-module 
$\fder\ \cF$ by~$\fvect(\cU)$.
The elements of the dual (right) 
$\cF$-module
\[
 \Covect(\cU) := (\fvect(\cU))^\ast = \Hom_{\cF}(\fvect(\cU),\cF) 
\]
are called \textit{covector fields}.
 The \textit{exterior differential} $d\colon \cF \to\Pi\Covect(\cU)$ 
is defined by the formula
\[
\textstyle df:=\sum dx_i\frac{\partial f}{\partial x_i}
\]
for any $f\in\cF$ 
and requirement $d^2=0$. 
The exterior algebra of $\Covect(\cU)$ over $\cF$ 
is the supercommutative superalgebra $\Omega^{\bcdot}(\cU) := E_{\cF}^{\bcdot}(\Covect(\cU))$
called the superalgebra of
\textit{differential forms} on~$\cU$. The algebra $\Omega^{\bcdot}(\cU) $ is naturally bigraded by degree of the forms (the superscript) and their parity. We will identity $\Omega^0(\cU)$ with $\cF$ and $\Omega^1(\cU)$ with
$\Pi(\Covect(\cU))$, where $\Pi$ is the reversal of parity functor.
Let the indeterminates $u_i$ be generators of the algebra $\cF$ of functions; denote by the $du_i$, where $p(du_i)=p(u_i)+\od$, the generators of the polynomial algebra of differential forms over $\cF$.

The \textit{interior product} with respect to~$X\in\fvect(\cU)$ is
the derivation~$\iota_X\colon\Omega(\cU)\to\Omega(\cU)$ such that
 for any $\alpha\in\Covect(\cU)$ and hence $\pi(\alpha)\in\Omega^1(\cU)$, we set 
 \[
\text{$\iota_X(\pi(\alpha)) =
(-1)^{p(X)}\langle X, \alpha \rangle$, where~}
\langle \sum f_i\partial_i, \sum dx_j \cdot g_j \rangle =\sum f_ig_i \text{~~for any $f_i, g_j\in \Cinf(\cU)$.}
\]
Clearly, $p(\iota_X) = p(X) + \od$ and $\deg \iota_X = -1$.

The differential form can be considered as a~multilinear function on~$\Pi(\fvect(\cU))$. In the literature, there are several formulas of such interpretation; in computations we use the following. Let $\tilde{X} =
\Pi(X)$ for any $X\in\fvect(\cU)$.
For any $\omega^l \in \Omega^l(\cU)$ and $X_1,\ldots,X_l\in\fvect(\cU)$, set
\[
 \omega^l(\tilde{X}_1,\ldots,\tilde{X}_l) =
 (-1)^{p(\omega^l)(p(X_1)+\ldots+p(X_l)+l)} \iota_{X_1}\ldots\iota_{X_l}\omega^l.
\]

Thanks to the possibility to turn any one-sided module over a~supercommutative superalgebra into a~two-sided one (for details, see \cite[Ch.1]{SoS}), we usually consider the differential forms as if $\Omega$ were a~left $\cF$-module; we always do so writing contact forms.

\ssec{Divided powers}\label{ssDivP}
Over $\Cee$, let $x_i$ be coordinates of the linear $m|n$-dimensional supermanifold, the first~$m$ of them
even, the rest being odd.
For an $(m+n)$-tuple of non-negative integers $\underline{r}=(r_1, \ldots , r_{m+n})$, where
$r_i=0$ or 1 for $i>m$, we set
\begin{equation}
\label{u^r} 
\textstyle u_i^{(r_{i})} := \frac{x_i^{r_{i}}}{r_i!}\quad
\text{and}\quad u^{(\underline{r})} := \prod\limits_{1\leq i\leq m+n}
u_i^{(r_{i})}.
\end{equation}
Over an~arbitrary field $\Kee$, we consider inseparable symbols $u_i^{(r_{i})}$. Clearly,
\begin{equation}
\label{divp}
\renewcommand{\arraystretch}{1.4}
\begin{array}{l} u^{(\underline{r})} \cdot u^{(\underline{s})} =
\left(\mathop{\prod}\limits_{m+1\leq i\leq m+n}
\min(1,2-r_i-s_i)\cdot(-1)^{\mathop{\sum}\limits_{m<i<j\leq m+n}
r_js_i}\right)\cdot \binom {\underline{r} + \underline{s}}
{\underline{r}}
u^{(\underline{r} + \underline{s})}, \\
\text{where}\quad\textstyle \binom {\underline{r} + \underline{s}}
{\underline{r}}:=\mathop{\prod}\limits_{1\leq i\leq m}\binom {r_{i}
+ s_{i}} {r_{i}}.
\end{array}
\end{equation}

Over any field~$\Kee$ of characteristic $p>0$, we consider the
supercommutative superalgebra
\begin{equation}\label{shear}
\cO(m; \un| n):=\Kee[u;
\un]:=\Span_{\Kee}\left(u^{(\underline{r})}\mid r_i
\begin{cases}< p^{N_{i}}&\text{for $i\leq m$}\\
=0\text{ or 1}&\text{for $i>m$}\end{cases}\right),
\end{equation}
with multiplication given by formula~\eqref{divp}
where $\un = (N_1,\dots, N_m)$ is the \textit{shearing
vector} whose coordinates $N_i\in\Zee_+\cup\infty$ (we assume that
$p^\infty=\infty$) are called the \textit{heights} of the respective indeterminates.
Set $\One:=(1,\dots,1)$. 

Over $\Cee$, consider the action of the partial derivative $\del_{x_i}$ of $\Cee[x_1, \dots, x_{m+n}]$ in the basis of divided powers. It is given by (recall the definition \eqref{u^r} of $u^{(r)}$)
\be\label{distPD}
\del_{x_i} u^{(\underline{r})} =
\begin{cases}
0 &\text{if~~} r_i = 0;\\
(-1)^{\max(0, i-m-1)} u^{((r_1,\dots,r_{i-1}, r_i-1, r_{i+1},\dots, r_a))} &\text{otherwise.}
\end{cases}
\ee
Since all the coefficients are integer, the map
$\del_{x_i}$ given by formula \eqref{distPD} is a~derivation of $\Kee[m;\un|n]$. We will denote this operator $\del_i$ and call the maps $\del_1, \dots, \del_a$ \textit{distinguished} partial derivatives each $\del_i$ serving as several
partial derivatives at once, for each of the generators $u_i$, $u_i^{(p)}$, $u_i^{(p^2)}$, \dots .

The \textit{general vectorial Lie superalgebra}
consists of the following derivations expressed in terms of \textit{distinguished} partial derivatives 
\[
\textstyle \fvect(m;\un|n)=\left\{\sum_{1\leq i\leq a} f_i\partial_i\mid f_i\in \cO(m;
\un|n)\right\}.
\]

Following Strade~\cite[vol.1, Index]{S}, we denote the module $\Vol^{\lambda}(m; \un |n)$ of
weighted densities of weight $\lambda\in\Kee$ (the rank-1 module over the algebra of functions $\cO(m; \un |n)$ on the $m|n$-dimensional supermanifold) by 
\be\label{volO}
\cO(m; \un |n)_{\lambda\Div}:=\cO(m; \un |n)\vvol(u)^{\lambda}, 
\ee
where ``div" stands for the ``divergence" defined by the formula
\be\label{divMod}
\Div (\sum f_i\partial_i):=\sum (-1)^{p(f_i)p(\partial_i)}\partial_i(f_i)
\ee
and the action of vector field~$D$ on~$\vvol(u)^{\lambda}$ is multiplication by~$\lambda\Div(D)$.

Following Strade, we denote the
subspace of volume forms with integral~0 by
\begin{equation}\label{O'}
\cO'(m; \un|n)_{\Div}:=\Span(u^{(a)}\vvol(u)\mid a_i<
p^{\un_i}-1\text{~~for all $i$}).
\end{equation}

\subsection{Critical coordinates and unconstrained shearing
vectors}\label{CritCoor} The coordinate of the shearing vector~$\un$
corresponding to an~even indeterminate of the $\Zee$-graded vectorial Lie (super)algebra~$\fg$ is said to
be \textit{critical} if it cannot take an arbitrarily big value.

The shearing vector~$\un$ without any imposed restrictions on its coordinates is said to be \textit{unconstrained}. 
\textbf{In what follows, the shearing vector of each Lie (super)algebra considered is unconstrained}, unless otherwise indicated.

\section{Lemmas on Weisfeiler regradings, see \cite{Sh14}}\label{SS:8.1} Let
$\fg=\mathop{\oplus}\limits_{i\geq -d}\fg_i$ be a~$\Zee$-graded
vectorial superalgebra. By a~\textit{Weisfeiler regrading} or just
W-regrading of $\fg$ we mean the passage from $\fg$ to the Lie
superalgebra $\fh=\mathop{\oplus}\limits_{i\geq -D}\fh_i$ isomorphic
to $\fg$ as abstract Lie superalgebra and such that
\begin{equation}
\begin{array}{ll}
a)&\text{$\fh$ is transitive, i.e., for any non-zero $X \in\fh_k$ where $k \geq 0$, there is a~$Y \in\fh_{-1}$}\\
&\text{such that $[X, Y]\neq 0$;}\\
b)&\text{$\fh_{\geq 0}:=\mathop{\oplus}\limits_{i\geq 0}\fh_i$ is
a~maximal
Lie subsuperalgebra of finite codimension;}\\
c)&\text{$\fh_0$-module $\fh_{-1}$ is irreducible};\\
d)&\text{$D<\infty$, i.e., $\fh$ is of finite depth.}
\end{array}
\end{equation}

In what follows in this section, $\fh$ denotes a~Lie superalgebra
obtained from $\fg$ by a~regrading. 
For our purposes it suffices to consider only \textit{compatible gradings} for which the decomposition
\begin{equation}
\label{7.1.1} \fh_{i}=\mathop{\oplus}\limits_{m(i)\leq j\leq M(i)}\fh_{i,
j},\text{ where $\fh_{i, j}=\fh_{i}\cap \fg_{j}$. }
\end{equation} is well-defined by definition: partly thanks to Skryabin, partly by hypothesis, we consider only gradings induced by the degrees of the indeterminates, see Subsection~\ref{gradings}.

We recall the proof of the next Lemmas from \cite{Sh14}.

\ssbegin[When $\fh_{-1}$ is homogeneous]{Lemma}[When $\fh_{-1}$ is homogeneous]\label{L4.1}
If $\fh_0\cap\fg_-=0$, then
$\fh_0\subset\fg_0$ and $\fh_{-1}$ is homogeneous.
\end{Lemma}

\begin{proof}
The Lie superalgebra $\fh_{0,0}$ transforms $\fh_{-1,j}$ into itself and the operators from $\fh_{0,k}$
send $\fh_{-1,j}$ into $\fh_{-1,j+k}$. Therefore, if the representation of $\fh_0$ on $\fh_{-1}$ is irreducible, then $\fh_{-1}$ is homogeneous with respect to the grading of $\fg$, i.e., $\fh_{-1} = \fh_{-1,j_0}$ for some $j_0$. But then, $\fh_{0,k}$ for $k > 0$ sends $\fh_{-1}$ to $0$. Since $\fh$ is transitive, $\fh_{0,k} = 0$ for all~$k > 0$, i.e., $\fh_0 = \fh_{0,0}$.\end{proof}

\ssbegin[When the gradings coincide]{Lemma}[When the gradings coincide]\label{L4.1.1}
If $\fh_0\cap\fg_-=0$ and
there exists a~nonzero ${X\in\fg_{-1}\cap \fh_-}$, then the gradings
of $\fh$ and $\fg$ coincide.
\end{Lemma}

\begin{proof} Since $\fh_{-1} =\fh_{-1,j_0}\subset\fg_{j_0}$, it follows that $\fh_{-k}\subset\fg_{k\cdot j_0}$. If $X\in\fh_{-k}$, then $-1=k\cdot j_0$, implying $k = 1$ and $j_0 = -1$, i.e., $\fh_{-1}\subset\fg_{-1}$, and therefore $\fh_{-1} = \fg_{-1}$ and $\fh_i = \fg_i$ for all~$i$.
\end{proof}

\ssec{How to describe all the $W$-regradings of the Lie (super)algebra $\fg$}\label{L3.3}\label{sec:HowTo}
Let~$\fh$ be a~Lie (super)algebra isomorphic to~$\fg$ as an~abstract algebra; namely, let~$\fh$ be a~Lie (super)algebra obtained from~$\fg$ by a~regrading. First, we should list all regradings not forbidden; our task is to select the ones that can be realized. If $p=0$, the even indeterminates can not be of degree $\leq 0$; if $p>0$, the indeterminates with unconstrained heights can not be of degree $\leq 0$. Steps

(i) Determine $\fh_{0}\cap\fg_{-}$.

(ii) Construct a~``minimal'' (better, though informally saying, ``most
tightly compressed'') regrading with the given intersection, i.e.,
such that preserves in $\fh_{-1}$ all the elements of $\fg_{-1}$
except for those that have to go away in view of the condition on the
intersection.

(iii) If the ``minimal'' regrading is a~Weisfeiler one, then with the
help of Lemmas \ref{L4.1} and \ref{L4.1.1} we prove that any other
$W$-regrading with the given intersection coincides with the minimal
one.

\ssec{Maurer--\/Cartan equations and regradings (after \cite{Shch})}\label{ss3.4} 
The above steps (i)--(iii) constitute a general method of describing $\Zee$-gradings of any Lie (super)algebra $\fg$. Here, we are interested in vectorial Lie (super)algebras.
Let $\fn = \bigoplus_{k= -d}^{-1}\fn_k$ be a $\Zee$-graded $m|n$-dimensional Lie (super)algebra such that
the superspace $\fn_{-1}$ generates~$\fn$. We denote by $\fv =\bigoplus_{k=-d}^{\infty}\fv_k$ the Lie (super)algebra~$\fvect(m|n)$ considered with the
nonstandard {depth-$d$} grading~\eqref{GradingMk}.
Let $f\colon \fn \to \fv_-:=\bigoplus_{k<0}\fv_k$ be an embedding
compatible with the grading~\eqref{GradingMk}. 

Let $\fg_{-}$ be the image of~$\fn$ in~$\fv$ and let $\{X_i\}_{1\leq i\leq m+n}$ be vector fields compatible with grading \eqref{GradingMk} and forming a basis
of $\fg_{-}$ with structure constants
\[
 [X_i,X_j] = \sum_{k}c^k_{ij} X_k,\qquad \text{where}~~c^k_{ij}\in\Kee.
\]
Let the $\omega^i$ be the dual 
differential 1-forms, i.e., $\omega^i(X_j) =
\delta^i_j$. Then the $\omega^i$ satisfy the following form of the Maurer\/--\/Cartan equation
\[
 d \omega^k = \nfrac12\sum_{i,j} (-1)^{p(X_i)p(X_j)+p(X_k)+p(X_j)} c^k_{ij}\omega^i\wedge\omega^j.
\]
We denote by the $Y_i$ the graded basis vectors of the centralizer $\mathfrak{cent}_{\fv_{-}}(\fg_{-})$ in~$\fv_{-}$, let $\Theta^i$ be the
1-form dual to $Y_i$. The elements $X\in\fg_s$ of the complete prolong~$\fg=\bigoplus_i\fg_i$ of $\fg_{-}$ can be written in the form
\[
 X = \sum_i\Theta^i(X) Y_i.
\]

For most of the Lie (super)algebras we consider, the initial grading $\fg=\oplus \fg_i$ embeds $\fg$ into
a~contact Lie (super)algebra of type $\fk$ or $\fm$. Accordingly, our first step is to single out~$\fg$ as
a~subalgebra in~$\fk$ or~$\fm$ by means of differential equations expressed in terms of $Y$-vectors.
As a~result, we get a~description of the elements of $\fg$ in terms of the functions generating $\fk$ or $\fm$,
and further one we work with this description of~$\fg$. 

The original descriptions of Skryabin algebras and of $\fvle$ do not embed them into $\fk$ or $\fm$.
Therefore, first of all we have to embed the components of non-positive degree into the Lie (super)algebra
$\fvect$ by means of the algorithm from~\cite{Shch}.
Seeking W-regradings we consider only regrading obtained by changing the degreed of
indeterminates~\eqref{GradingMk}.
So, the vectors $X_i$ remain homogeneous in any regrading, together with the vectors $Y_i$ and forms $\Theta^i$ dual to them.


\section{For $p=0$, mainly}\label{p=0} 

In this Section, we give some definitions and results for any $p$, but formulated for $p=0$ as a~point of reference. These results are known, we give them only for completeness of the picture.

\ssec{The two analogs of contact forms}\label{2contact} There is only one equivalence class of normal shapes
of the contact form.
(For the proof of this classical fact generalized to the super case over~$\Cee$,
see~\cite{BGLS};
for the case of modular Lie algebras, see~\cite{Sk1}; to consider the case of modular Lie superalgebras preserving a~(peri)contact structure is an \textbf{Open problem}; conjecturally, there is just one equivalence class.) Let the indeterminates~$p_i$, $q_i$ and~$t$ be even, the~$\xi_j$, $\eta_j$ and~$\theta_s$ be 
odd. The contact Lie superalgebra $\fk (2n+1|m)$ preserves
the distribution singled out by the 1-form (we mostly consider the cases $\ell=0$ or $1$)
\begin{equation}
\label{2.2.6} \alpha_1=dt+\mathop{\sum}\limits_{1\leq i\leq
n}(p_idq_i-q_idp_i)+ \mathop{\sum}\limits_{1\leq j\leq
k}(\xi_jd\eta_j+\eta_jd\xi_j)+
\begin{cases}0&
\text{ if }\ m=2k\\
\mathop{\sum}\limits_{1\leq s\leq
\ell}\theta_sd\theta_s&\text{ if }\ m=2k+\ell.\end{cases} 
\end{equation}

It is subject to a~direct verification that $\fk (2n+1|m)=\Span(K_f\mid f\in\Cee [t, p,
q, \theta])$, where
\begin{equation}
\label{K_f} 
\textstyle K_f=(2-E)(f) \nfrac{\partial}{\partial t}-H_f + \nfrac{\partial f}{\partial t} E,
\end{equation}
where $E=\mathop{\sum}\limits_i y_i \pder{y_{i}}$ (here the $y_{i}$
are all the coordinates except $t$) is the \textit{Euler operator}, and $H_f$ is the hamiltonian field with Hamiltonian $f$ that
preserves $d\alpha_1$:
\begin{equation}
\label{H_f} \textstyle H_f=\mathop{\sum}\limits_{i\leq n}\left(\pderf{f}{p_i}
\pder{q_i}-\pderf{f}{q_i} \pder{p_i}\right )
-(-1)^{p(f)}\left(\mathop{\sum}\limits_{j\leq
k}\left(\pderf{f}{\xi_j} \pder{\eta_j}+ \pderf{f}{\eta_j}
\pder{\xi_j}\right)+ \mathop{\sum}\limits_{j\leq l}\pderf{ f}{\theta_j}
\pder{\theta_j}\right).
\end{equation}
Indeed, $L_{K_f}\alpha_1=2\frac{\partial f}{\partial t}\alpha_1$.

In~\cite{L1}, another analog of the contact Lie superalgebra was discovered: it is the Lie superalgebra~$\fm(n)$ preserving the distribution singled out by the \textit{pericontact} form $\alpha_0$ given here in its normal shape~\eqref{alp0} also of just one equivalence class, see~\cite{BGLS} (to prove this in the modular case is an \textbf{Open problem}), where~$\tau$ and the~$\xi_i$ are odd indeterminates whereas the~$q_i$ are even:
\begin{equation}
\label{alp0}
\begin{array}{c}
\alpha_0=d\tau+\mathop{\sum}\limits_{1\leq i\leq n}(\xi_idq_i+q_id\xi_i).\end{array}
\end{equation}

It is subject to a~direct verification that $\fm(n)=\Span(M_f\mid f\in\Cee [\tau, 
q, \xi])$, where
\begin{equation}
\label{M_f} M_f=(2-E)(f) \nfrac{\partial}{\partial \tau}- Le_f -(-1)^{p(f)}
 \nfrac{\partial f}{\partial \tau} E, 
\end{equation}
where $E=\mathop{\sum}\limits_iy_i \pder{y_i}$ (here the $y_i$ are
all the coordinates except $\tau$), and
\begin{equation}
\label{Le_f} 
\textstyle Le_f=\mathop{\sum}\limits_{i\leq n}\left(
\pderf{f}{q_i}\ \pder{\xi_i}+(-1)^{p(f)} \pderf{f}{\xi_i}\
\pder{q_i}\right).
\end{equation}
Indeed, $L_{M_f}\alpha_0=-(-1)^{p(f)}2\frac{\partial f}{\partial \tau}\alpha_0$.

The two super analogs of the Lie algebra of Hamiltonian vector fields are 
\[
\begin{array}{l}
\fh(2n|m):=\{D\in\fvect(2n|m)\mid L_D(d\alpha_1)=0\}=\Span(H_f\mid f\in\Cee [p,
q, \theta]);\\
\fle(n):=\{D\in\fvect(n|n)\mid L_D(d\alpha_0)=0\}=\Span(Le_f\mid f\in\Cee [q, \xi]).
\end{array}
\]

\subsection{Weisfeiler regradings}\label{WeisFReg} If $p=0$, we can identify the elements of the nilpotent Lie algebra $\fg_-:=\oplus_{-d\leq i<0}\ \fg_{i}$ with elements of negative degree in $\fvect(n|m; \vec r)$ and the elements of $\fg_0$ with elements of 0th degree of $\fvect(m|n; \vec r)$ in a~non-standard (see Subsection~\ref{WeisReg}) grading $\vec r$ (the vector of degrees of indeterminates) of $\fvect(m|n)$, where $m|n=\sdim \fg_-$. Then, the \textit{Cartan prolong} is $(\fg_-, \fg_0)_{*}:=\oplus_{k\geq -d}\ \fg_{k}$, where the terms of Cartan prolongation for any $k>0$ are defined as
\[
\fg_k:=\{D\in \fvect(n|m; \vec r)_k\mid [D, \fg_i]\subset \fg_{k+i}\text{~~for any $i<0$}\}.
\]
Lately, the above-described procedure was often called \textit{generalized} prolongation because the Cartan prolongation was recalled in the literature of the 1960s for $d=1$ only, when the fact that Cartan himself used the prolongation for any finite $d$ was successfully forgotten.

If $p>0$, any Cartan prolong $(\fg_-, \fg_0)_{*}$ can be likewise embedded into $\fvect(m;\un|n)$, perhaps in a~non-standard grading of the latter; for details, see \cite{BGLLS2}.

We do not indicate $\vec r$ in the case of the \textit{standard grading} where --- for the \textbf{serial} Lie superalgebras --- the degrees of all indeterminates are equal to 1, except for the contact (resp. pericontact) Lie superalgebra $\fk$ (resp. $\fm$) where the degree of the \lq\lq time" indeterminate $t$ (resp.~$\tau$) is equal to 2; generally, $\vec r$ is shorthanded to the number $r$ of indeterminates of degree 0. The $\Zee$-grading of the Lie superalgebra of series $\fk$ or $\fm$ induces the grading of its simple subalgebra considered in \cite[Table (17)]{LS}; the W-graded Lie superalgebras $\fg$ of depth~2 listed in \cite[Table (17)]{LS} are not needed in this paper since they 
preserve the same distributions as the appropriately regraded the Lie superalgebras of series $\fk$ and $\fm$. 

All the W-gradings are listed in table~\eqref{nonstandgr}.
The standard realizations are marked by an asterisk 
$(*)$. Observe that the codimension of ${\cal L}_0$ attains its
minimum in the standard realization.
\footnotesize
\begin{equation}\label{nonstandgr}
\renewcommand{\arraystretch}{1.3}
\begin{tabular}{|c|c|}
\hline Lie superalgebra & its $\Zee$-grading \\ \hline

$\fvect (n|m; r)$, & $\deg u_i=\deg \xi_j=1$ for any $i, j$
\hfill $(*)$\\

\cline{2-2} $ 0\leq r\leq m$ & $\deg \xi_j=0$ for $1\leq j\leq r;$\\
&$\deg u_i=\deg \xi_{r+s}=1$ for any $i, s$
\\

\hline $\fm(n; r),$ & $\deg \tau=2$, $\deg q_i=\deg \xi_i=1$ for any
$i$ \hfill $(*)$\\ \cline{2-2} $\; 0\leq r< n-1$& $\deg \tau=\deg
q_i=2$, $\deg \xi_i=0$ for $1\leq i\leq r <n$;\\

& $\deg q_{r+j}=\deg \xi_{r+j}=1$ for any $j$\\
\hline

$\fm(n; n)$ & $\deg \tau=\deg q_i=1$, $\deg \xi_i=0$ for $1\leq
i\leq n$ \\ \hline

\cline{2-2} $\fk (2n+1|m; r)$, & $\deg t=2$, whereas, for any $i,
j, s$\hfill $(*)$\\
see eq.~\eqref{2.2.6};& $\deg
p_i=\deg q_i= \deg \xi_j=\deg \eta_j=\deg \theta_s=1$ \\

\cline{2-2} $0\leq r\leq [\frac{m}{2}]$& $\deg t=\deg \xi_i=2$,
$\deg \eta_{i}=0$ for $1\leq i\leq r\leq [\frac{m}{2}]$; \\
$r\neq k-1$ if $m=2k$ and $n=0$&$\deg p_i=\deg q_i=\deg
\theta_{j}=1$
for $j\geq 1$ and all $i$\\

\hline $\fk(1|2m; m)$ & $\deg t =\deg \xi_i=1$, $\deg \eta_{i}=0$
for $1\leq i\leq m$ \\ \hline
\end{tabular}
\end{equation}
\normalsize

\sssec{Notation in tables \eqref{1114} 
and \eqref{table32} and further on}\label{Other} The ground field is $\Kee$. \textbf{Let $p=0$ or $p>3$ till the end of the subsection}. If $p=3$ (resp. $2$), some modifications are needed in the tables, see (resp.~\cite{BGLLS2} and~\S\ref{Sp=2}).

The \textit{general linear} Lie
superalgebra of all supermatrices of size $\text{Size}=(p_1, \dots, p_{|\text{Size}|})$, is the ordered
collection of parities of the rows identical to that of the columns and equal to the parities of the basis vectors of the superspace $V$, is denoted by
$\fgl(\text{Size})$;
usually, for the \textit{standard} format (one of the simplest formats), $\fgl(\ev,
\dots, \ev, \od, \dots, \od)$ is abbreviated to $\fgl(\dim V_{\bar
0}|\dim V_{\bar 1})$. 
The \textit{supertrace} is the linear mapping $\fgl (\text{Size})\tto \Kee$ vanishing on the \textit{special linear} Lie subsuperalgebra
$\fsl(\text{Size})$, the commutant of $\fgl (\text{Size})$.

For $X\in \fgl_{\cC}(\Size)$ over a~supercommutative superalgebra $\cC$ having only the $ij$th non-zero matrix element $X_{ij}$, we have $p(X)=p_i+p_j+p(X_{ij})$, where $X_{ij}\in\cC$. The \textit{supertransposed} matrix $X^{st}$ is defined by the formula
\[
(X^{st})_{ij}=(-1)^{(p_{i}+p_{j})(p_{i}+p(X))}X_{ji}.
\]
In the standard format, this means
\[
X=\mmat A,B,C,D,\longmapsto X^{st}:=\begin{cases}\begin{pmatrix}A^t&C^t\\-B^t&D^t\end{pmatrix}&\text{if $p(X)=\ev$},\\
\begin{pmatrix}A^t&-C^t\\B^t&D^t\end{pmatrix}&\text{if $p(X)=\od$}.\\
\end{cases}
\]
The supermatrices $X\in\fgl(\text{Size})$ such that
\[
X^{st}B+(-1)^{p(X)p(B)}BX=0\quad \text{for a~homogeneous matrix
$B\in\fgl(\text{Size})$}
\]
constitute the Lie superalgebra $\faut (B)$ that preserves the
bilinear form $\cB$ on $V$ with the Gram matrix $B$. In order to identify the form $\cB\in \Bil(V, W)$ with an element of $\Hom(V, W^*)$, the Gram matrix $B=(B_{ij})$ of $\cB$ is given by the formula
\be\label{martBil}
B_{ij}=(-1)^{p(\cB)p(v_i)}\cB(v_{i}, v_{j})\text{~~for the basis vectors $v_{i}\in V$}.
\ee

Consider the \textit{upsetting} of bilinear forms
$u\colon\Bil (V, W)\tto\Bil(W, V)$ 
given by the formula \be\label{susyB}
u(\cB)(w, v):=(-1)^{p(v)p(w)}\cB(v,w)\text{~~for any $v \in V$ and $w\in W$.}
\ee
In terms of the Gram matrix $B$ of $\cB$: the form
$\cB$ is \textit{symmetric} if and only if 
\be\label{BilSy}
u(B)=B,\;\text{ where $u(B)=
\mmat R^{t},(-1)^{p(B)}T^{t},(-1)^{p(B)}S^{t},-U^{t},$ for $B=\mmat R,S,T,U,$.}
\ee
Similarly, $\cB$ is \textit{anti-symmetric} if $u(B)=-B$. Observe
that \textbf{the passage from $V$ to $\Pi (V)$ turns every symmetric
form $\cB$ on $V$ into an anti-symmetric form $\cB^\Pi$ on $\Pi (V)$ and anti-symmetric $\cB$ on $V$ into symmetric $\cB^\Pi$ on $\Pi (V)$ by setting}
\[
\text{$\cB^\Pi(\Pi(x), \Pi(y)):=(-1)^{p(B)+p(x)+p(x)p(y)}\cB(x,y)$ for any $x,y\in V$}.
\]

The Lie superalgebra preserving a~non-degenerate anti-symmetric even (resp. odd) bilinear form on $m|2n$-dimensional (resp. $n|n$-dimensional) superspace is denoted $\fosp(m|2n)$ (resp. $\fpe(n)$); let $\fspe(n):=\fpe(n)\cap\fsl(n|n)$.


$\fc(\fg)$ or $\fc\fg$ is the trivial central extension of the Lie (super)algebra $\fg$ by 1-dimensional even center.

$\fii\inplus \fa$ is the semidirect sum of the Lie (super)algebras, where $\fii$ is an ideal and $\fa$ is the quotient.

$\Kee[k]$ is the 1-dimensional $\Kee$-module over a~1-dimensional commutative Lie algebra $\Kee z$ spanned by an element $z$ which acts multiplying $1\in\Kee\simeq\Kee[k]$ by $k\in\Zee$.

$\Lambda(n)=\Lambda(\xi)$ is the Grassmann algebra with $n$ generators $\xi=(\xi_1, \dots, \xi_n)$.

$\fsl(m;\Lambda(n)):=\fsl(m)\otimes\Lambda(n)$ which is sometimes more convenient notation than $\fsl_{\Lambda(n)}(m)$.

The divergence of a~vector field in $\fvect(m|n)$ in coordinates $(x, \xi)$ is defined by eq.~\eqref{divMod}. The subalgebra of divergence-free vectorial Lie (super)algebra $\fg$ is called \textit{special} and denoted $\fs\fg$, e.g.,

$\fsvect(m|n):=\{D\in\fvect(m|n)\mid \Div D=0\}$.

$\fsvect_{a, b}(0|n):=\fsvect(0|n)\inplus \Kee(az+bD)$, where $a,b\in\Kee$ and
$D$ is
the operator that determines the standard $\Zee$-grading of $\fsvect(0|n)$,
whereas $z$ is a~central element. 

$\id_{\fg}$ denotes the tautological module $V$ over $\fg\subset\fgl(V)$. 

The \textit{Heisenberg} Lie superalgebra $\fhei(2n|m)$ is isomorphic to the sum of components of negative degree of the Lie superalgebra $\fk(2n+1|m)$ in its standard $\Zee$-grading.

In $\fk\fas(; 3\xi)_0$ and $\fk\fle(9|6; 2)_0$, see Table~\eqref{table32}, the subalgebra $\fsl(1|3)$ is considered embedded in $\fvect(0|3)$, it acts on $\Lambda(3)$ via the $\fvect(0|3)$-action.

For any $\fg_0$-module $V$ with \textit{lowest} weight
$\lambda$ and \textit{even} lowest weight vector, we make $V$ into a~$\fg_{\geq}$-module by setting $\fg_{+}(V)=0$, where $\fg_{+}=\mathop{\oplus}\limits_{i>
0}\fg_{i}$ and $\fg_{\geq}=\mathop{\oplus}\limits_{i\geq
0}\fg_{i}$; we realize $\fg=\fg_-\oplus \fg_0\oplus \fg_+$ by vector fields on the linear
supermanifold $\cK^{m|n}$ corresponding to the superspace $\fg_-:=\mathop{\oplus}\limits_{i<0}\fg_{i}$. Let $x=(u, \xi)$ be coordinates on $\cK^{m|n}$; then, the
superspace $T(V):=\Hom_{U(\fg_{\geq})}(U(\fg), V)$ is isomorphic, due
to the Poincar\'e--Birkhoff--Witt theorem, to ${\Cee}[[x]]\otimes
V$ and the elements of $T(V)$ have a~natural interpretation as formal
\textit{tensor fields of type} $V$; when the highest/lowest weight of $V$ is $\lambda=(a, \dots , a)$ we
simply write $T(\vec a)$ instead of $T(\lambda)$. 

We denote the generator of the $\fvect(m|n)$-module of volume forms $\Vol(m|n)$, as $\cF$-module,
corresponding to the ordered set of coordinates $x=(u,\theta)$ by $\vvol(x)$ or just $\vvol$;
the space of weighted $a$-densities for any weight $a\in\Kee$ is (compare with definition~\eqref{volO})
\[
\Vol^{a}(m|n)=T(\underbrace{a,\dots , a}_m; \underbrace{-a, \dots , -a}_n). 
\]
In particular, 
$\Vol^{a}(m|0)=T(\vec a)$, but
$\Vol^{a}(0|n)=T(\overrightarrow{-a})$. 

In particular, $\fvect(0|n)$ acts in the space of semi-densities for any $n$ and $D\in
\fvect(0|n)$: 
\begin{equation}\label{T12}
\textstyle T^{1/2}_D(f\vvol^{1/2})=(D(f)-\nfrac12 (-1)^{p(f)p(D)}f\cdot \Div D)\vvol^{1/2}.
\end{equation}

We set (do not confuse it with the space of functions $\cF:=\Vol^{0}(m|n)$):
\begin{equation}
\label{vol}
\textstyle \Vol_0(0|n):=\{v\in \Vol(0|n)\mid\int v=0\} \text{ and
$T_0(\vec 0):=\Lambda(n)/\Cee\cdot 1$.}
\end{equation}

As $\fsvect(m|n)$-modules, $\Vol\cong T(\vec 0)$. So, over
$\fsvect(0|n)$ we can set
\begin{equation}
\label{T00} T^0_{0}(\vec
0):=\Vol_0(0|n)/\Cee\vvol(\xi).\end{equation}

We represent an arbitrary
element $\widehat X\in\fas$ as a~pair $\widehat X=X+a\cdot 1_{4|4}$, where $X\in\fspe(4)$ and $a\in{\Cee}$. The bracket in $\fas$ is
\begin{equation}
\label{2.1.4}
\left[\mat {A &B\\ C&-A^t} +a\cdot 1_{4|4},
\mat{A' & B' \cr C' & -A'{}^t} +a'\cdot
1_{4|4}\right]= \left[\mat{ A& B \cr C & -A^t},
\mat{ A' & B' \cr C' & -A'{}^t
}\right]+\tr(C\widetilde C')\cdot 1_{4|4},
\end{equation}
where $\ \widetilde {}\ $ is the result of any
even permutation $(1234)\longmapsto(ijkl)$ extended via linearity from matrices
$C_{ij}:=E_{ij}-E_{ji}$ on which $\widetilde C_{ij}=C_{kl}$.

The Lie superalgebra $\fas$ can also be described with the help of
the spinor representation. For this, we need several vectorial
superalgebras. Recall that the Poisson Lie superalgebra is a~Lie subsuperalgebra of $\fk(2n+1|2m)$ (recall the expression of $K_f$, see \eqref{K_f}):
\[
\fpo(2n|2m):=\Span(K_f\mid f\in\Cee[p,q, \xi, \eta])
\]
with the
Poisson bracket (cf. with the contact bracket \eqref{kb1}, \eqref{kb2})
\begin{equation}
\label{pb}
\renewcommand{\arraystretch}{1.4}
\textstyle 
 \{f, g\}_{P.b.}=\sum \left(\pderf{f}{p_i}\ \pderf{ g}{q_i}-\pderf{f}{q_i}\
\pderf{ g}{p_i}\right)-(-1)^{p(f)} \sum \left(\pderf{f}{\xi_j}\ \pderf{ g}{\eta_j}+\pderf{f}{\eta_j}\
\pderf{ g}{\xi_j}\right).
\end{equation}

Recall that 
\[
\fh(2n|2m)=\fpo(2n|2m)/\mathfrak{center}, \text{~~where $\mathfrak{center}=\Span(K_1)$}. 
\]

Since
$\fsl(4)\cong\fo(6)$, we can identify $\fspe(4)_0$ with
$\fh(0|6)_0$.
It is not difficult to see that the elements of degree $-1$ in the
standard gradings of $\fspe(4)$ and $\fh(0|6)$ constitute isomorphic modules over
$\fsl(4)\cong\fo(6)$, so it is possible to embed $\fspe(4)_1$ into $\fh(0|6)_1$.

The extension $\fas$ discovered by A.~Sergeev is the result of the restriction to
$\fspe(4)\subset\fh(0|6)$ of the cocycle that turns $\fh(0|6)$ into
$\fpo(0|6)$.
The quantization deforms $\fpo(0|6)$ into
$\fgl(\Lambda(\xi))$; the through maps 
\[
T_\lambda:
\fas\tto\fpo(0|6)\tto\fgl(\Lambda (\xi))
\]
are representations of
$\fas$ in the $4|4$-dimensional modules $\spin_\lambda$ isomorphic
to each other for all~${\lambda\neq 0}$, so we assume $\lambda=1$ and do not indicate it. The explicit form of
$T_\lambda$ is as follows:
\begin{equation}
\label{spinLambda} 
\textstyle T_\lambda\colon \mat{A~& B \cr C & -A^t }+a\cdot
1_{4|4}\longmapsto \mat{A~& B-\lambda \widetilde C \cr C & -A^t
}+\lambda a\cdot 1_{4|4}, 
\end{equation}
where 
$\widetilde C$ is defined in
the line under eq.~\eqref{2.1.4}. Clearly, $T_\lambda$ is 
irreducible for any $\lambda$.

\underline{For $N=11$} in Table \eqref{1114}, the non-positive components of W-graded contact algebras are as follows (borrowed from \cite{LS}); the notation 
$\id_{\fosp}\boxtimes\Cee[-1]$ means the tensor product of the tautological $\fosp_{\Lambda(r)}(m-2r|2n)$-module $\id_{\fosp_{\Lambda(r)}(m-2r|2n)}$ by the $\fc$-module $\Cee[-1]$, where 
\[
\fc\subset \fc\fosp(m-2r|2n)\subset \fc\fosp_{\Lambda(r)}(m-2r|2n):=\fc\fosp(m-2r|2n)\otimes \Lambda(r)
\]
is the center and where $\Lambda(r)$ is a~natural module over $\fvect(0|r):=\fder\ \Lambda(r)$. The line 14 of Table~\eqref{1114} is analogous; excluded are the values of $r$ corresponding to non-Weisfeiler gradings. 

Set (recall~\eqref{M_f}) 
\begin{equation}
\label{2.7.2} 
\fb_{a, b}(n) :=\{M_f\in \fm (n)\mid a\; \Div
M_f=(-1)^{p(f)}2(a-bn)\nfrac{\partial f}{\partial\tau}\}. 
\end{equation}
Explicitly, in coordinates $q, \xi,
\tau$, see \eqref{M_f}, let $\Delta:=\mathop{\sum}_{i\leq
n}\frac{\partial^2 }{\partial q_i
\partial\xi_i}$ and $E=\mathop{\sum} (q_i\partial_{q_i}+\xi_i\partial_{\xi_i})$; then,
\begin{equation}
\label{2.7.4}
\begin{array}{ll}
\fb_{a, b}(n) &=\{M_f\in \fm (n)\mid (bn-aE)\nfrac{\partial f}{\partial\tau} =
a\Delta
f\}=\\
&\{D\in\fvect(n|n+1) \mid L_{D}(\vvol^a\alpha_{0}^{a-bn})=0\}. 
\end{array}
\end{equation}
 
For $\lambda =\frac{2a}{n(a-b)}\in\Cee \Pee^1$, we can abbreviate $\fb_\lambda(n):=\fb_{a, b}(n)$. Since $\fsm(n)\simeq \fb_{2/(n-1)}(n)$, the family $\fb_\lambda(n)$ is a~deformation of $\fsm(n)$, preserving the same distribution for all $\lambda\neq 0$.

\underline{For $N=21$, 22} in Table \eqref{1114}: the terms ``$\fg_{-i}$'' denote the
superspace isomorphic to the one in quotation marks, but with the
action given by the rules~\eqref{1.3.5} and~\eqref{1.3.6}.

\underline{For $N=21-23$}: let $D:=\diag
(1_{k}, -1_{k})\in \fpe(k)$, let $V$ denote the tautological
$k|k$-dimensional $\fspe(k)$-module. Let $W=V\otimes \Lambda(r)$ and $X\in
\fvect(0|r)$. Let 
\be\label{top}
\Xi=\xi_1\cdots\xi_n\in \Lambda(\xi_1, \dots
,\xi_n).
\ee
Denote by $T^r$, where $r$ is integer, the representations of $\fvect(0|r)$ in
$W$ \index{$T^r$, the representations of $\fvect(0\vert r)$} in $\fspe(n-r)\otimes\Lambda(r)$ and $\fspe(n-r)\otimes
\Lambda(r)$ given by the formula
\begin{equation}
\label{1.3.4} T^r(X)=1\otimes X+D\otimes \nfrac{1}{n-r}\Div X.
\end{equation}
\underline{$\fg=\fsle'(n; r)_{0}$ for $r\neq n-2$}. For $\fg_{0}$,
we have:
\begin{equation}
\label{1.3.5}
\begin{array}{l}
\text{$\fvect(0|r)$ acts on the ideal $\fspe(n-r)\otimes
``\text{$\Lambda(r)$}"$ via
$T^r$, see eq.~\eqref{1.3.4}};\\
\text{any $X\otimes f\in\fspe(n-r)\otimes \Lambda(r)$ acts in
$\fg_{-1}$ as
$\id\otimes f$ and in $\fg_{-2}$ as $0$};\\
\text{any $X\in\fvect(0|r)$ acts in $\fg_{-1}$ via $T^r$ and in
$\fg_{-2}$ as $X$.}
\end{array}
\end{equation}
\underline{$\fg=\fsle'(n; n-2)$}. For $\fg_{0}$, we note that
\[
\fspe(2)\simeq \Cee(\Le_{q_1\xi_1-q_2\xi_2})\inplus
\Cee\Le_{\xi_1\xi_2},
\]
whereas $\fg_{-2}$ and $\fg_{-1}$ are as above, for $r< n-2$. Set
$\fh=\Cee(\Le_{q_1\xi_1-q_2\xi_2})$. In this case
\begin{equation}
\label{1.3.6} \fg_0\simeq \underline{``\text{$\left (\fh\otimes
\Lambda(n-2)\inplus \Cee\Le_{\xi_1\xi_2}\otimes
(\Lambda(n-2)\setminus \Cee\xi_3\cdots\xi_n)\right )$}"}\inplus
T^1(\fvect(0|n-2)).
\end{equation}
The action of $\fvect(0|n-2)$, the quotient of $\fg_{0}$ modulo the
underlined ideal in quotation marks in \eqref{1.3.6} is performed via~\eqref{1.3.5}. In the subspace
$\xi_1\xi_2\otimes\Lambda(n-2)\subset\fg_0$, this action is the same
as in the space of volume forms. So we can consider everything, except for
the terms proportional to $\Xi$, see eq.~\eqref{top}, or, speaking correctly, take
the irreducible submodule of functions with integral 0.

\underline{For $N=27$, 32, 33}: the terms ``in quotation marks'' denote the
superspace isomorphic to the one in quotation marks, but with the
action given by eqs.~\eqref{1.3.5} and~\eqref{1.3.7}. For $N=32$, we have
\be\label{32}
\fg_0=((\fspe(n{-}r)_{n, n-2})\otimes \text{``$\Lambda(r)$''})\inplus T^r(\fvect(0|r)).
\ee

In the exceptional case $ar=bn$, \textit{i.e.},
$\lambda= \frac{2}{n-r}$, we see that the
$\fvect(0|r)$-action on the ideal
$(\fc\fspe(n-r)\otimes\Lambda(r))\inplus \fvect(0|r)$ of $\fg_0$, and
on $\fg_-$, is the same as for $\fsle'$, see eq.~\eqref{1.3.5}.

In table~\eqref{1114}, to save space in line 18 we denote the $\fpe(n-r;\Lambda
(r))\inplus \fvect(0|r)$-module $\id_{\fpe(n-r;\Lambda
(r))}$ by $\id_{\fpe}$, and we use similar abbreviations in other lines.
\index{$\fh_{\lambda}(2\vert 2)$}

\underline{For $N=24-34$}: 1) We assume that
$\lambda=\nfrac{2a}{n(a-b)}\neq 0$, 1, $\infty$; these three
exceptional cases (corresponding to the ``drop-outs'' $\fle(n)$,
$\fb_{1}'(n)$ and $\fb'_{\infty}(n)$, respectively) are considered
se\-pa\-rately.

2) The irreducibility condition of the $\fg_0$-module $\fg_{-1}$ for
$\fg=\fb'_{\infty}$ excludes $r=n-1$.

3) The case where $r=n-2$ is extra exceptional, so in table~\eqref{1114} we assume that
\begin{equation}
\label{*1} 0<r<n-2;\; \; \text{ additionally $a\neq b$ and $(a,
b)\neq \alpha(n, n-2)$ for any $\alpha\in\Cee$.} 
\end{equation}

\underline{For $N=24-27$}: we consider $\fb_{a,b}(n; r)$ for
$0<r<n-2$ and $ar-bn\neq 0$; in particular, this excludes
$\fb'_{\infty}(n; n)=\fb'_{a, a}(n; n)$ and $\fb'_{1}(n;
n-2)=\fb'_{n, n-2}(n; n-2)$.

If $z$ is the central element of $\fc\fspe(n-r)$ that acts on
$\fg_{-1}$ as $-\id$, then
\begin{equation}
\label{1.3.7} z\otimes \psi\text{ acts on $\fg_{-1}$ as
$-\id\otimes \psi$, and on $\fg_{-2}$ as $-2\id\otimes\psi$.} 
\end{equation}
Set
\begin{equation}
\label{c} c:=\nfrac{a}{ar-bn}.
\end{equation}

Let $x\cdot \str\otimes\id$ be the representation of 
$\fpe(n-r)=\fs\fpe(n-r)\inplus \Cee D$, which is $\id_{\fs\fpe(n-r)}$
on $\fs\fpe(n-r)$ and sends $D=\diag(1_{n-r}, -1_{n-r})\in \fpe(n-r)$ to $2x\cdot\id$. 

\subsection{The exceptional simple vectorial Lie superalgebras 
 as Cartan prolongs and their W-gradings, see \cite{BGLLS2}}
The five families of such Lie superalgebras are
given below in their W-gradings as Cartan prolongs $(\fg_{-1},\fg_{0})_{*}$
or generalized Cartan prolongs $(\fg_{-},\fg_{0})_{*}$ 

For 
$\fg_{-}=\mathop{\oplus}\limits_{-2\leq i\leq -1}\fg_{i}$, we write
$(\fg_{-2}, \fg_{-1}, \fg_{0})_{*}$ 
instead of $(\fg_{-}, \fg_{0})_{*}$.
In Table~\eqref{table32}, indicated is also one of
the Lie superalgebras from the list of series~\eqref{nonstandgr} as an
ambient, in a~grading described in eq.~\eqref{clsfeq10}, which contains the exceptional one as a~maximal subalgebra.
The $W$-graded superalgebras of \underline{depth $3$} appear as
regradings of the listed ones at certain values of $\vec r$; for the
corresponding terms $\fg_{i}$ for $i\leq 0$, see Table~\eqref{table32} and eqs.~\eqref{a1}, \eqref{a2}. 

Denote by $\fb\fa(n)$ the
\textit{antibracket} Lie superalgebra ($\fb\fa$ is ``Anti-Bracket'' read
backwards) whose space is $W\oplus \Cee\cdot z$, where $W$ is an
$n|n$-dimensional superspace endowed with a~non-degenerate
anti-symmetric \textbf{odd} bilinear form $B$ and $z$ is odd; the bracket in $\fb\fa(n)$ is
given by the following relations:
\begin{equation}\label{2.5.4}
 \text{$z$ spans the center; $[v,w]=B(v, w)\cdot z$ for any $v, w\in W$.}
\end{equation}
In the exceptional cases, the vector of degrees $\vec r$ after a~semi-colon is often shorthanded by a~symbol
$r$, and in the case $r=0$, where no indeterminate is of degree $0$, the sign $r=0$ is omitted, as in serial
cases. For example, for $\fk\fa\fs^\xi$ which is the one of two isomorphic subalgebras of $\fk(1|6)$, namely
the one containing $K_{\xi_1\xi_2\xi_3}$, the meaning of notation $r=0$, $1\xi$, $3\xi$ is clear: 
none, or one or three of the $\xi$'s have degree 0 (and the
corresponding $\eta$'s acquire degree 2), same applies to $3\eta$.
The letter $K$ denotes a~W-grading compatible with parity, the exceptional vector $\vec r=CK$ denotes the
grading discovered by Cheng and Kac~\cite{CK2}, see~\eqref{table32}.
Tables \eqref{clsfeq10} 
describe these regradings. 
\tiny
\be\label{1114}
\renewcommand{\arraystretch}{1.3}
\begin{tabular}{|c|c|c|c|c|}
\hline $N$&$\fg$&$\fg_{-2}$&$\fg_{-1}$&$\fg_0$\cr \hline \hline
$7$&$\fh(2n|m; r),\ 2r<m$&$T_0(\vec
0)$&$\id_{\fosp}$&$\fosp(m-2r|2n;\Lambda (r))\inplus \fvect(0|r)$\cr
\hline

$11$&\parbox{35mm}{\tiny{\mbox{}\\ $\fk(2n+1|m;
r)$ for $r\neq k-1$ if $m=2k$ and $n=0$ \\ 
}}&$\Lambda(r)\boxtimes\Cee[-2]$&$\id_{\fosp}\boxtimes\Cee[-1]$&$\fc\fosp_{\Lambda
(r)}(m-2r|2n)\inplus \fvect(0|r)$\cr \hline


$13$&$\fk(1|2m+1; m)$&$\Lambda(m)$&$\Pi(\Lambda (m))$&$\Lambda
(m)\inplus \fvect(0|m)$\cr \hline
\hline

$14$&$\fm(n;
r),\text{~~where~~}r\neq n$&$\Pi(\Lambda(r))\boxtimes\Cee[-2]$&$\id_{\fpe}$&
$\fc\fpe_{\Lambda
(r)}(n-r)\inplus 
\fvect(0|r)$\cr \hline

$18$&$\fle(n; r)$&$\Pi(T_0(\vec 0))$&$\id_{\fpe}$&$\fpe(n-r;\Lambda
(r))\inplus \fvect(0|r)$\cr \hline


$21$&$\fsle'(n; r),\ r\neq n-2, n$&``$\Pi(T_0(\vec 0))$'' & $\id_{\fspe}\boxtimes``\Lambda
(r)$'' &$(\fspe(n-r)\otimes ``\Lambda (r)")\inplus
T^r(\fvect(0|r))$\cr \hline

$22$&$\fsle'(n; n-2)$&``$\Pi(T_0(\vec 0))$'' &$\id\boxtimes``\Lambda
(r)$''&see eq.~\eqref{1.3.6}\cr 
\hline
$24$&$\fb_{\lambda}(n)$&$\Pi(\Cee[-2])$&$\id$&
 $\fspe(n)\inplus\Cee(az+bD)$\cr
 \hline

$25$&$ \fb_{\lambda}(n;r),\ \ r<n-2$
 & \tiny{$\Pi(-c\,\str)\boxtimes
 \Vol^{2c}(0|r)$} &
 $(-\frac{c}{2}\,\str\otimes\id)
 \boxtimes\Vol^{c}(0|r)$&\tiny{
 $(\fpe(n{-}r)\otimes\Lambda(r))\inplus \fvect(0|r)$}\\
 \hline


 $27$ &
 $ \fb_{2/(n-r)}(n; r),\ \
 r<n-2$
 & $\Pi(\Cee)\boxtimes \Lambda
 (r)$&$\id\boxtimes ``\Lambda (r)$''& \tiny{$\fc\fpe(n-r)\otimes ``\Lambda
 (r)$'' $\inplus T^r(\fvect(0|r))$}\cr
 \hline

 $28$&$\fb'_{\infty}(n)$&$\Pi(\Cee)$&$\id_{\fspe}$& $\fspe(n)_{a, a}$\cr
 \hline

 $29$&
 $
 \fb'_{\infty}(n; r),\ \
 r<n-2$ & $\Pi(\Cee)\boxtimes \Lambda(r)$
 &$\id_{\fspe}\boxtimes\Lambda(r)$& \tiny{$((\fspe(n-r)_{a,a})\otimes\Lambda (r))\inplus \fvect(0|r)$}\cr
 \hline


 $31$&$\fb'_{1}(n)$&$\Pi(\Cee)$&$\id_{\fspe}$&$\fspe(n)_{n, n-2}$\cr
 \hline

 $32$&
 $\fb'_{1}(n; r),\ \
 r<n-2$
 &``$\Pi(\Vol_{0}(0|r))$''& $\id\boxtimes ``\Lambda(r)$''
 & see eq.~\eqref{32}\cr 
 \hline

 $33$&$\fb'_{1}(n; n-2)$&``$\Pi(T_0(\vec 0))$'' & $\id_{\fspe}\boxtimes``\Lambda (r)$'' &
 \parbox{45mm}{\tiny{\mbox{}\\ $((\fspe(n{-}r)_{n, n-2})
 \otimes ``\Lambda(r)$'') see ${N=32}$
 for $\fc\fspe(2)$ instead of $\fspe(2)$}}
 \cr \hline
\hline
\end{tabular}
\ee
\normalsize
\sssbegin{Theorem}[\cite{LS1}]\label{K&M} The non-isomorphic W-gradings of simple Lie superalgebras over $\Cee$ preserving non-integrable distributions are listed in Table \eqref{1114} borrowed from the complete description of W-graded simple vectorial Lie superalgebras in \cite{LS1}. For the growth vectors of the corresponding non-integrable distributions, see eqs.~\eqref{grVecK} and \eqref{grVecM}. The non-simple prolongs $\fg$ preserve the same distributions as their simple derived algebras $\fg'$. 

For any $p>2$, the same description, mutatis mutandis, provides with examples of symmetries of non-integrable distributions over $\Kee$. 
\end{Theorem} 

\subsection{The $W$-gradings of the exceptional Lie
superalgebras}\label{SS:1.3.4}
In Table~\eqref{table32}, there are given the terms $\fg_{i}$
for $-2\leq i\leq 0$ of the 15 exceptional $W$-graded algebras; ``un" is the number of parameters of unconstrained coordinates of $\un$, see Subsection~\ref{CritCoor}.
\begin{equation}\label{table32}\tiny
\renewcommand{\arraystretch}{1.3}
\begin{tabular}{|c|c|c|c|c|c|}
\hline $\fg$&$\fg_{-2}$&$\fg_{-1}$&$\fg_0$&$\sdim\fg_{-}$&un\cr

\hline

\hline $\fv\fle(4|3)$&$-$&$\Pi(\Lambda(3)/\Cee
1)$&$\fc(\fvect(0|3))$&$4|3$&3\cr

\hline

$\fv\fle(4|3; 1)$&$\Cee[-2]$&$\id_{\fsl(2;\Lambda(2))}$&
$\fc(\fsl(2;\Lambda(2)))\inplus T^{1/2}(\fvect(0|2))$&$5|4$&3\cr

\hline

$\fv\fle(4|3; K)$&$\id_{\fsl(3)}\boxtimes
\Cee[-2]$&$\Pi(\id^*_{\fsl(3)}\boxtimes \id_{\fsl(2)}\boxtimes
\Cee[-1])$&$\fsl(3)\oplus\fsl(2)\oplus\Cee z$&$3|6$&3\cr

\hline \hline

$\fv\fas(4|4)$&$-$&$\spin$&$\fas$&$4|4$&1\cr \hline \hline

$\fk\fas$&$\Cee[-2]$&$\Pi(\id_{\fo(6)})$& $\fc\fo(6)$&$1|6$&1\cr

\hline

$\fk\fas(; 1\xi)$
&$\Lambda(1)$&$\id_{\fsl(2)}\boxtimes\id_{\fgl(2;\Lambda(1))}$&
$\fsl(2)\bigoplus(\fgl(2;\Lambda(1))\inplus \fvect(0|1))$&$5|5$&1\cr

\hline

$\fk\fas(; 3\xi)$&$-$&
$\Lambda(3)$&$\Lambda(3)\inplus \fsl(1|3)$&$4|4$&1\cr 
\hline $\fk\fas(;
3\eta)$&$-$&$\Vol_{0}(0|3)$& $\fc(\fvect(0|3))$&$4|3$&1\cr

\hline \hline

$\fm\fb(4|5)$&$\Pi(\Cee[-2])$&$\Vol
(0|3)$&$\fc(\fvect(0|3))$&$4|5$&3\cr \hline

$\fm\fb(4|5; 1)$&$\Lambda(2)/\Cee 1$ &$\id_{\fsl(2;\Lambda(2))}$
&$\fc(\fsl(2;\Lambda(2))\inplus T^{1/2}(\fvect(0|2))$&$5|6$&3\cr

\hline

$\fm\fb(4|5; K)^{\star}$&$\id_{\fsl(3)}\boxtimes
\Cee[-2]$&$\Pi(\id^*_{\fsl(3)}\boxtimes \id_{\fsl(2)}\boxtimes
\Cee[-1])$&$\fsl(3)\oplus\fsl(2)\oplus\Cee z$&$3|8$&3\cr

\hline \hline $\fk\fle(9|6)$&$\Cee[-2]$&$\Pi(T^0_{0}(\vec
0))$&$\fsvect(0|4)_{3, 4}$&$9|6$&5\cr

\hline

$\fk\fle(9|6; 2)$&$\Pi(\id_{\fsl(1|3)})$&$\id_{\fsl(2;\Lambda(3))}$&
$\fsl(2;\Lambda(3)) \inplus \fsl(1|3)$&$11|9$&5\cr

\hline

$\fk\fle(9|6; K)$&$\id_{\fsl(5)}$&$\Pi(\Lambda^2(\id_{\fsl(5)}^*))$&$\fsl(5)$&$5|10$&5\cr

\hline

$\fk\fle(9|6; CK)^{\star}$&$\id_{\fsl(3;\Lambda(1))}^*$&$\id_{\fsl(2)}\boxtimes
\id_{\fsl(3;\Lambda(1))}$&$\fsl(2)\bigoplus\left (\fsl(3;\Lambda(1))
\inplus \fvect(0|1)\right)$&$9|11$&5\cr

\hline
\end{tabular}
\end{equation}
\normalsize

The star ${\star}$ marks the Lie superalgebras of depth 3 whose
components $\fg_{-3}$ are described in eqs.~\eqref{a1} and~\eqref{a2}. Observe that none of the simple $W$-graded vectorial Lie
superalgebras is of depth $>3$ and only two algebras are of depth 3:
$\fm\fb(4|5; K)$, for which we have
\begin{equation}\label{a1}
\fm\fb(4|5; K)_{-3}\cong \Pi(\id_{\fsl(2)}),
\end{equation}
and another one, $\fk\fle(9|6; CK)$, for which we have
\begin{equation}\label{a2}
\fk\fle(9|6; CK)_{-3}\simeq \Pi(\id_{\fsl(2)}\boxtimes\Cee[-3]).
\end{equation}

In tables \eqref{clsfeq10}, 
there are listed the exceptional superalgebras, their indeterminates and their respective
degrees in the regrading $R(r)$ of the ambient superalgebra, see \cite{BGLLS2}.
{\tiny 
\begin{table}[ht]\centering
%
\begin{equation}\label{clsfeq10}
\setcounter{MaxMatrixCols}{16}
\begin{tabular}{|l|l|}
\hlx{hv}
$\fv\fle(4|3)$ & $ R(0)=\begin{pmatrix}
x_1&x_2&x_3&y&|&\xi_1&\xi_2&\xi_3\\ 1&1&1&1&|&1&1&1\\ \end{pmatrix}$\cr \hlx{vv}

$\fv\fle(5|4; 1)$ &$R(1)=\begin{pmatrix}
x_1&x_2&x_3&y&|&\xi_1&\xi_2&\xi_3\\ 2&1&1&0&|&0&1&1\\ \end{pmatrix}$\cr \hlx{vv}

$\fv\fle(3|6; K)$ &$R(K)=\begin{pmatrix} 
x_1&x_2&x_3&y&|&\xi_1&\xi_2&\xi_3\\
2&2&2&0&|&1&1&1\\\end{pmatrix}$\cr \hlx{vhv}

\hlx{hv}
$\fm\fb(4|5)$ &$R(0)=\begin{pmatrix} 
x_0&x_1&x_2&x_3&|&\xi_0&\xi_1&\xi_2&\xi_3&\tau \\
1&1&1&1&|&1&1&1&1;&2\\\end{pmatrix}$\cr
\hlx{vv}
$\fm\fb(5|6; 1)$ & $R(1)=\begin{pmatrix} 
x_0&x_1&x_2&x_3&|&\xi_0&\xi_1&\xi_2&\xi_3&\tau \\
0&2&1&1&|&2&0&1&1;&2\\\end{pmatrix}$\cr
\hlx{vv}
$\fm\fb(3|8; K)$ &$R(K)=\begin{pmatrix}
x_0&x_1&x_2&x_3&|&\xi_0&\xi_1&\xi_2&\xi_3&\tau \\
 0&2&2&2&|&3&1&1&1;&3\\\end{pmatrix}$\cr
\hlx{hv} 

\hlx{vhv}
$\fk\fas(1|6)$ &$R(0)=\begin{pmatrix}
t&|&\xi_1&\xi_2&\xi_3&\eta_1&\eta_2&\eta_3 \\
 2&|&1&1&1&1&1&1\\\end{pmatrix}$\cr
\hlx{vv}
$\fk\fas(5|5; 1\xi)$ &$R(1\xi)=\begin{pmatrix}
t&|&\xi_1&\xi_2&\xi_3&\eta_1&\eta_2&\eta_3 \\
 2&|&0&1&1&2&1&1\\\end{pmatrix}$\cr
\hlx{vv}
$\fk\fas(4|4; 3\xi)$ &$R(3\xi)=\begin{pmatrix} 
t&|&\xi_1&\xi_2&\xi_3&\eta_1&\eta_2&\eta_3 \\
1&|&0&0&0&1&1&1\\\end{pmatrix}$\cr
\hlx{vv}
$\fk\fas(4|3; 3\eta)$ &$R(3\eta)=\begin{pmatrix}
t&|&\xi_1&\xi_2&\xi_3&\eta_1&\eta_2&\eta_3 \\
 1&|&1&1&1&0&0&0\\\end{pmatrix}$\cr
 \hlx{vh}


\hlx{vhv}
$\fk\fle(9|6)$ &$R(0)=\begin{pmatrix} 
q_1&q_2&q_3&q_4&p_1&p_2&p_3&p_4&t&|&\xi_1&\xi_2&\xi_3&
\eta_1&\eta_2&\eta_3\\
1&1&1&1&1&1&1&1;&2&|&1&1&1&1&1&1\\
\end{pmatrix}$\cr
\hlx{vv}
$\fk\fle(11|9; 2)$ &$R(2)=\begin{pmatrix} 
q_1&q_2&q_3&q_4&p_1&p_2&p_3&p_4&t&|&\xi_1&\xi_2&\xi_3&
\eta_1&\eta_2&\eta_3\\
1&1&2&2&1&1&0&0;&2&|&0&1&1&2&1&1\\
\end{pmatrix}$\cr
\hlx{vv}
$\fk\fle(5|10; K)$ &$R(K)=\begin{pmatrix} 
q_1&q_2&q_3&q_4&p_1&p_2&p_3&p_4&t&|&\xi_1&\xi_2&\xi_3&
\eta_1&\eta_2&\eta_3\\
2&2&2&2&1&1&1&1;&2&|&1&1&1&1&1&1\\
\end{pmatrix}$\cr
\hlx{vv}
$\fk\fle(9|11; CK)$&$R(CK)=\begin{pmatrix} 
q_1&q_2&q_3&q_4&p_1&p_2&p_3&p_4&t&|&\xi_1&\xi_2&\xi_3&
\eta_1&\eta_2&\eta_3\\
3&2&2&2&0&1&1&1;&3&|&2&2&2&1&1&1\\
\end{pmatrix}$\cr
\hlx{vh}
\end{tabular}
\end{equation}
\end{table}
}

\normalsize

\[
\textbf{Hereafter, the contact and pericontact distributions are marked by a~\lq\lq C".}
\]

\sssbegin{Theorem}[\cite{Sh5, Sh14,LS}]\label{Exc} The non-isomorphic W-gradings of exceptional vectorial Lie superalgebras are classified in Table~\eqref{table32}, see \cite{LS}. For the growth vectors on non-integrable distributions, see Table~\eqref{grVex}. The same is true, mutatis mutandis, over $\Kee$ for any $p>2$.
\begin{equation}\label{grVex}\footnotesize
\renewcommand{\arraystretch}{1.3}
\begin{array}{ll}
\begin{tabular}{|c|c|}
\hline 
$\fg$&$v_\fg$\cr
\hline
\hline
$\fv\fle(4|3; 1)$&$(4|4, 5|4)C$\cr
\hline
$\fv\fle(4|3; K)$&$(0|6, 3|6)$\cr
\hline \hline
$\fk\fle(9|6)$&$(8|6,9|6)C$\cr
\hline
$\fk\fle(9|6; 2)$&$(8|8,11|9)$\cr
\hline
$\fk\fle(9|6; K)$&$(0|10, 5|10)$\cr
\hline
$\fk\fle(9|6;
CK)$&$(6|6, 9|9,9|11)$\cr
\hline
\end{tabular} &

\begin{tabular}{|c|c|}
\hline $\fg$&$v_\fg$\cr
\hline
\hline
$\fk\fas$&$(0|6, 1|6)C$\cr
\hline
$\fk\fas(; 1\xi)$
&$(4|4,5|5)$\cr
\hline \hline
$\fm\fb(4|5)$&$(4|4, 4|5)C$\cr \hline
$\fm\fb(4|5; 1)$&$(4|4,5|6)$\cr
\hline
$\fm\fb(4|5; K)$&$(0|6, 3|6, 3|8)$\cr
\hline \hline 
\end{tabular}
\end{array}
\end{equation}\end{Theorem}

\ssec{Weisfeiler regradings of contact algebras}\label{WeisReg} The $\Zee$-gradings of vectorial Lie superalgebras are defined by the vector $\vec r:=(\deg x_1, \dots, )$ of degrees of the indeterminates. \textbf{For W-gradings}, this
vector can often be shorthanded by the number $r$ of indeterminates of degree 0, so instead of $\fvect(n|m;\vec r)$ we write $\fvect(n|m; r)$ and similarly for the subalgebras of $\fvect(n|m; r)$, bar exceptional examples, see Table~\eqref{table32}. 

Recall that from the $K$-functor point of view, the \textit{superdimension} of the superspace $V$ is $a+b\eps$, where $a=\dim V_\ev$, $b=\dim V_\od$ and $\eps^2=1$, often shorthanded as $a|b$. The \textit{growth vector} $v_\fk$ of the non-integrable distribution preserved by the contact Lie superalgebra 
$\fk (2n+1|m;r)$, with the case where $n=0$, $m=2k$ and $r=k$ excluded (because it does not correspond to any W-grading), see eq.~\eqref{nonstandgr}, is as follows:
\be\label{grVecK}
\begin{array}{ll}
r=0&v_\fk=(2n|m, 2n+1|m)\\
1\leq r\leq [\frac{m}{2}] &v_\fk=(D_\fk, D_\fk+2^{r-1}(1+\eps)),\text{~~where $D_\fk=(2n|m-2r)2^{r-1}(1+\eps)$}.
\end{array}
\ee

The growth vector $v_\fm$ of the non-integrable distribution preserved by the pericontact Lie superalgebra 
$\fm (k|k+1;r)$ --- usually abbreviated to $\fm (k;r)$, see eq.~\eqref{nonstandgr}, --- is as follows:
\be\label{grVecM}
\begin{array}{ll}
r=0&v_\fm=(k|k, k|k+1)\\
1\leq r<k-1&v_\fm=(D_\fm, D_\fm+2^{r-1}(1+\eps)),\text{~~where $D_\fm=(k-r|k-r)2^{r-1}(1+\eps)$}.
\end{array}
\ee

Observe that
\be\label{obs}
\text{If $k=m-r=2n+r$ and $r>0$, then $v_\fk=v_\fm$. }
\ee

In~\S\ref{Sequiv}, we show that the two distributions satisfying condition \eqref{obs} are inequivalent.

\ssec{Peculiarities of $p$ small}\label{ssPsmall}
Table~\eqref{table32} clearly demonstrates possible values of~$p$ for which the incarnations of the algebra might differ from its incarnations for $p$ large or $p=0$. These are almost all cases for $p=2$ carefully studied in~\cite{BGLLS2}; we did not consider these cases here, except for what is described in Section~\ref{Sp=2}. Let us list the other cases:

$p=5$: $\fk\fle(9|6; K)_0\simeq \fsl(5)$ has a~center now. Nevertheless, this does not affect the negative part we are interested in; and the components $\fk\fle(9|6; K)_i$ for $i<0$ remain irreducible. 

$p=3$: the component of $\fkas_1$ containing $\xi_1\xi_2\xi_3$ remains irreducible, as for any other~$p$. The component $\fk\fle(9|6)_0$ is obtained by reduction of parameters $(3,4)$ modulo 3 (direct verification).

\section{$p=5$, Melikyan algebra, see \cite{Me, GL}}\label{SMel}

\ssbegin{Theorem}\label{ThMe} The two known W-gradings of the Melikyan algebra $\fme$ exhaust all its W-gradings: \textup{$1)$} preserving the contact structure, \textup{$2)$} preserving the same distribution as is preserved by the exceptional simple Lie algebra $\fg(2)$ in one of its $\Zee$-gradings; the growth vector being $(2,3,5)$. \end{Theorem} 

This section is devoted to a~proof of this theorem.

\ssec{Melikyan's first description of $\fme(\un)$}\label{ssMel1} The following expressions of the basis elements of non--positive degree in terms of the indeterminate $u$ of
$\fvect(1; \underline{1})$ and the $\fvect(1; \underline{1})$-module $\cO(1;1)$ rewritten in terms of the generating functions of elements of
$\fk(5;\un)$, i.e., in terms of the indeterminates
$t; p_1, p_2, q_1, q_2$ constituted the initial Melikyan's construction of $\fme(\un)$ as the Cartan prolong of the following components. Since table (40) in the journal version of \cite{GL} is incorrect, let us correct it and give the details for verification (this will be used in the latest arXiv version of \cite{GL}). First of all, recall the formula for the contact bracket (with a~minus in front of the Poisson bracket); in what follows we skip subscript k.b.:
\be\label{kb1}
\{f,g\}_{k.b.}=(2-E)(f)\del_t(g)-\del_t(f)(2-E)(g)-\sum(\del_{p_i}(f)\del_{q_i}(g)- \del_{q_i}(f)\del_{p_i}(g)),
\ee
where $E:=\sum y_i\del_{y_i}$ and the $y_i$ are all the indeterminates except $t$.
Hence, $\{p_i,q_i\}=-1$. The formula 
\[
(f,g):=\int fdg =\text{the coefficient of $u^{(4)}$ in the expansion of $fdg$}
\]
defines a~non-degenerate anti-symmetric bilinear form on $\cO(1;1)$:
\[
\begin{array}{l}
(u^{(4)},u)=\int u^{(4)}du=1,\\
(u^{(3)},u^{(2)})=\int u^{(3)}udu=4=-1.
\end{array}
\]
This implies the description of a~basis of $\fme_{-1}$:
\[
p_1\leftrightarrow -u^{(4)}, \ \ p_2\leftrightarrow u^{(3)}, \ \ q_2\leftrightarrow u^{(2)}, \ \ q_1\leftrightarrow u.
\]
Therefore,
 $\del_u \leftrightarrow -p_2q_1+q_2^{(2)}$. The action of $u\del_u$ on $\cO(1;1)$ in terms of $u$ and $p, q$ is as follows:
\[
u\del_u: \begin{cases}-u^{(4)}\mapsto -u\cdot u^{(3)}=u^{(4)},\; u^{(3)}\mapsto u\cdot u^{(2)}=3u^{(3)}=-2u^{(3)},\; u^{(2)}\mapsto u\cdot u=2u^{(2)},\; u\mapsto u;\\
p_1\mapsto -p_1,\; p_2\mapsto -2p_2,\; q_2\mapsto 2q_2, \; q_1\mapsto q_1.\end{cases}
\]
Hence,
\[
u\del_u \leftrightarrow -p_1q_1-2p_2q_2.
\]

Similarly,
\[
\begin{array}{l}
u^{(2)}\del_u \leftrightarrow -p_1q_2+2p_2^{(2)},\\
u^{(3)}\del_u\leftrightarrow -p_1p_2, \\
u^{(4)}\del_u\leftrightarrow p_1^{(2)}.\\
\end{array}
\]
Clearly, $\fme_{-2}=\Kee 1$, and the center of $\fme_{0}$ is $z\leftrightarrow t$.

In this realization, all elements of degree $-1$ and the corresponding indeterminates look as if they are on equal footing and it is absolutely unclear why $\un$ depends only of 2 parameters, the height of the other 3 indeterminates being constrained. However, they are NOT on equal footing! They would have been if they were interchangeable under the action of $\fme_{0}$, i.e., if $\fme_{0}$ had acted on them symmetrically.

\sssec{Description of $\fme_1$}\label{ssDescrMe1} The generating functions independent of $t$ are as follows
\[
p_1^{(3)}, \; p_1^{(2)}p_2, \; p_1^{(2)}q_2-2p_1p_2^{(2)}, \; p_1^{(2)}q_1+2p_1p_2q_2+p_2^{(3)}, \;
p_1q_2^{(2)}-p_1p_2q_1-2p_2^{(2)}q_2.
\]

The fact that these functions belong to $\fme_1$ is subject to a~direct verification. The fact that there are no other functions independent of $t$ follows from the dimension considerations since, as is not difficult to see,
\be\label{fG}
\text{the functions independent of $t$ generate the central extension $\fG$ of $\fh(2; r)$ for $r=(1,0)$}.
\ee

Observe that different $\fG$ is not isomorphic to the known central extension --- the Poisson Lie algebra $\fpo(2; r)$ with the height of the indeterminate of degree 0 being constrained, see Subsection~\ref{H^2}. 
In more detail: first, the quotient of the space of functions independent of $t$ modulo $\Kee\ 1$ generates the Hamiltonian Lie algebra $\fh(2; r)$. However, the generating function 1 in $\fme$ is not the same as generating function 1 in the Poisson Lie algebra. In the regrading $\fpo(2; r)$ of the Poisson Lie algebra on the space of functions in $p$ and $q$ the depth is equal to 1 and all functions depending only on $p$ --- they span $\fpo(2; r)_{-1}$ --- commute. Whereas the depth of the Melikyan algebra $\fme$ is equal to 2, $\fme_{-1}$ is spanned by functions in $p$, except 1, and they do not commute but yield $1\in\fme_{-2}$.

Let $G:=-p_2q_1+q_2^{(2)}$ be the generating function corresponding to $\del_u$. Then, the functions depending on $t$ are as follows:
 \[
 \begin{array}{l}
 F_{p_1}=tp_1+2p_1p_2q_2-p_2^{(3)}, \\
 F_{p_2}=\{F_{p_1},G\}=tp_2+p_1q_2^{(2)}-2p_1p_2q_1, \\
 F_{q_2}=\{F_{p_2},G\}=-tq_2+p_2q_2^{(2)}+p_2^{(2)}q_1+p_1 q_1q_2,\\
 F_{q_1}=\{F_{q_2},G\}=tq_1+2q_2^{(3)}-p_2q_1q_2-2p_1q_1^{(2)}.
 \end{array}
 \]

 The operators corresponding to the generating functions in $\fme_-$ are as follows:
 \begin{equation}\label{K-}
 K_{p_i}:=p_i\del_t-\del_{q_i} \text{~and~} K_{q_i}:=q_i\del_t+\del_{p_i} \text{~for~} i=1,2, \text{~and~}K_1:=2\del_t.
 \end{equation}

As is explained in \cite{Shch}, to select the generating functions belonging to $\fme$ from all generating
functions of $\fk(5)$, we should consider the differential equations expressed in terms of the operators
commuting with the operators~\eqref{K-}; these operators, called in \cite{Shch} and in what follows by
``$Y$-vectors", are 
\begin{equation}\label{hatK-}
 \widetilde K_{p_i}=p_i\del_t+\del_{q_i} \text{~and~} \widetilde K_{q_i}=q_i\del_t-\del_{p_i}\text{~for~} i=1,2, \text{~and~}
 \widetilde K_1=\del_t.
 \end{equation}

Observe that all generating functions in $\fme_0$ satisfy, in particular, the following equations:
\begin{equation}\label{eqme}
(\widetilde K_{p_1})^2 (f)=0, \quad \widetilde K_{p_1}\widetilde K_{p_2}(f)=0.
\end{equation}

Observe that to single out the whole algebra $\fme\subset \fk(5)$ we need more equations. 
Observe that the operators $\widetilde K_{p_i}$ act on the functions $f$ independent of $t$ just as
$\del_{q_1}$. Hence, on such functions, the equations~\eqref{eqme} turn into
\begin{equation}\label{fpq}
\textstyle \frac{\del^2 f}{(\del_{q_1})^2}=0, \quad \frac{\del^2 f}{\del_{q_1}\del_{q_2}}=0 \Lra f=q_1\vf(p_1,p_2)+\psi(p_1,p_2,q_2).
\end{equation}

If $\deg_t f=k>0$, i.e., $f=t^k g(p_1,p_2,q_2,q_1)+h$, where $\deg_t h<k$, then taking $k$ times the contact bracket of $f$ with $1$, we see that $g\in\fme$, and hence is also of the form~\eqref{fpq}.

\paragraph{What is $\fG$, see definition~\eqref{fG}?}\label{H^2}
Since the algebra acts by zero on its cohomology, it is easy to see that the weight 0 (modulo 5) of the cocycle might occur only in products of elements of negative degree, so the central extension can only increase the depth of the algebra $\fh(2; \un)$, so it suffices to calculate the extensions setting $\un=\One$ to simplify computations.
The space $H^2(\fh(2;\One))$ is spanned by the following cocycles whose indexes are their degrees in the grading $(1,-1)$ and $\what{x}$ denotes the element of $\Pi(\fg^*)$ dual to $x\in\fg$:
\begin{align*}
 c_{-5} = {} & \what{p^{(2)}}\wedge \what{p^{(3)}} - \what{p}\wedge\what{p^{(4)}},
 & c_0 = {}& \what{p}\wedge\what{q},
 & c_{5} = {} & \what{q^{(2)}}\wedge \what{q^{(3)}} - \what{q}\wedge\what{q^{(4)}}. 
\end{align*}
 \normalsize
The cocycle $c_0$ corresponds to the Poisson algebra, the other two cocycles correspond to algebras isomorphic
to~$\fG$.
The whole picture with the center~$\fc$ corresponding to the 3 cocycles and the centrally extended algebra denoted by ``?" is as follows: 
\[
 0\tto \fc\tto ? \tto \fh(2;\One)\tto 0.
\]
It resembles, especially since $\fpsl(2|2)\simeq \fh'(0|4)$, the deform of $\fosp_0(4|2)$, see~\cite[p.37]{BGLLK}:
\[
 0\tto \fsl(2)\tto \fosp_0(4|2)\tto \fpsl(2|2)\tto 0.
\]

\ssec{A W-regrading $3\fme(\un)$ of $\fme(\un)$}\label{ssKuzMel}
From Yamaguchi's theorem cited in~\cite{GL} and~\cite{Shch}
we know that the Cartan prolong of $(\fg(2)_-,\fg(2)_0)$ (in any $\Zee$-grading of $\fg(2)$)
is isomorphic to $\fg(2)$, at least,over $\Cee$.
Kuznetsov, see~\cite{Ku1}\footnote{
 Melikyan informed us that he deposited this description to VINITI depositions, but failed to publish it.} observed that, for $p=5$, the non-positive parts
of~$\fg(2)$ in the grading of depth~3 are isomorphic to the
respective non-positive parts of the Melikyan algebra in one of
its $\Zee$-gradings --- let us denote it $3\fme(\un)$. The Melikyan algebra $3\fme(\un)$ is the Cartan prolong of this non-positive part. In more detail: let $V[k]$ be the $\fgl(V)$-module which is $V$ as $\fsl(V)$-module, and let a~fixed central element $z\in
\fgl(V)$ act on $V[k]$ as $k\cdot \id$. Then,
\begin{equation}\label{melg2}
\renewcommand{\arraystretch}{1.4}
\footnotesize{
\begin{tabular}{|c|c|c|c|}
\hline
$\fg_0$&$\fg_{-1}$&$\fg_{-2}$&$\fg_{-3}$\cr \hline
$\fgl(2)\simeq\fgl(V)$&$V=V[-1]$&$E^2(V)$&$V[-3]$\cr
\hline
\end{tabular}
}
\end{equation}
For a~basis of these components in terms of vector fields, see~\cite{Shch}, where the ground field is $\Cee$, but the formulas remain valid in characteristic~$5$ as well if we consider indeterminates as divided powers:
\footnotesize
\begin{equation}\label{realme}
\renewcommand{\arraystretch}{1.4}
\begin{tabular}{|l|l|}
\hline $\fg_{0}$&$X_- = u_2 \del_1 + u_1 u_2^{(2)} \del_4 + u_2^{(3)} \del_5
+ u_5
\del_4$,\\
&$X_+ = u_1 \del_2 + 2 u_1^{(3)} \del_4 + u_4 \del_5$,\\
&$h_1 = u_1 \del_1 + u_3 \del_3 + 2 u_4 \del_4 + u_5 \del_5$,\\
&$h_2 = u_2 \del_2 + u_3 \del_3 + u_4 \del_4 + 2 u_5 \del_5$\\
\hline\end{tabular}\quad \renewcommand{\arraystretch}{1.4}
\begin{tabular}{|l|l|}
\hline $\fg_{-1}$&$\del_1 -u_2 \del_3- u_1 u_2 \del_4 - u_2^{(2)}
\del_5$;\quad
$\del_2$\\
\hline
$\fg_{-2}$&$\del_3 + u_1 \del_4 + u_2 \del_5$\\
\hline
$\fg_{-3}$&$\del_4$;\quad $\del_5$\\
\hline
\end{tabular}
\end{equation}
\normalsize 

Unlike $p\neq 5$ case, the complete prolong of
$\fg(2)_-\oplus \fg(2)_0$ in the depth 3 grading of $\fg(2)$
strictly contains $\fg(2)$, see~\cite{GL}. 
In this realization, we still do not see why the shearing vector $\un$ depends only of 2 parameters, but at least this realization corroborates the following empirical fact formulated with the word ``usually'', because sometimes it is violated: \textbf{usually, the heights of the indeterminates of non-maximal degree in a~certain W-grading are constrained}, see~\cite{GL}.


\begin{proof}[Proof\nopoint] of Theorem~\ref{ThMe} follows from a~series of lemmas.
 Let $\fg$ be a~W-regrading of~$\fme$. Since the generating functions~$p_1$ and~$t$ enter~$\fme$ in however high power,
 whereas bracketing with~$q_1$ acts on the functions independent of~$t$ as~$\del_{p_1}$,
 it follows that~$q_1$ and~$1$ must belong to~$\fg_-$.

 In this situation, Lemma~\ref{L4.1.1} tells us that if~$\fg_0\cap \fme_-=0$, then the grading of~$\fg$ coincides with the grading of~$\fme$.

\textit{A priori} $\dim (\fme_-\cap \fg_0)\leq 3$. We have $3$ elements that can belong to this intersection: $p_1$, $p_2$ and $q_2$ whose weights with respect to $p_1q_1+2p_2q_2\in\fme_0$ are different. They do not form a~basis as we will show. 

If $\dim (\fme_-\cap \fg_0)=2$, there are lots of 2-dimensional planes these 3 elements span, but we are interested only in the planes spanned by some two of these 3 elements. 

If $\dim (\fme_-\cap \fg_0)=1$, the situation is similar.

Consider the W-regraded algebra $\fg$ with the following Lie grading $\gr$ (\textit{grade}):
\be\label{gr}
 \begin{array}{l}
\deg p_1=3,\ \ \deg q_1=0, \ \ \deg p_2=2, \ \ \deg q_2=1, \ \ \deg t=3\\
\gr(f)=\deg f-3.
 \end{array}
 \ee
Clearly, the depth of $\fg$ is equal to 3. Only functions of degree 0, i.e., only functions depending on $q_1$
can belong to $\fg_{-3}$.
Taking eq.~\eqref{fpq} into account we see that $\fg_{-3}=\Span(1,q_1)$.
The component $\fg_{-2}$ is manned by the functions of degree~1, i.e., of the form $q_2\vf(q_1)$, hence, with
eq.~\eqref{fpq} being taken into account, $\fg_{-2}=\Span(q_2)$.
Finally, $\fg_{-1}$ is manned by the functions of degree~2, which with eq.~\eqref{fpq} being again taken into
account, means that $\fg_{-1}=\Span(p_2,G=-q_1p_2+q_2^{(2)})$.
Besides, $\fg_0=\Span(t,\ p_1q_1+2p_2q_2, \ p_1, \ F_{q_1})$.
In this way, we get the depth-3 W-grading of the Melikyan algebra, earlier denoted by $3\fme$.

\sssbegin{Lemma}\label{L1}
If $\fg$ is a~W-regrading of $\fme$ such that $p_1\in\fg_0$, then $\fg=3\fme$.
\end{Lemma}

\begin{proof} Observe that since $\{p_1, G\}=-p_2$, it follows that $p_2$ and $G$ must belong to the same component $\fg_k$. Therefore, $q_2=\{G,p_2\}\in\fg_{2k}$ and $q_1=\{G,q_2\}\in\fg_{3k}$. But $q_1\in\fg_-$, so $k<0$, i.e., $p_2, G, q_2\in\fg_-$. Thus, $\fg_0\cap 3\fme_-=0$ and $p_2\in\fg_-\cap 3\fme_{-1}$; both conditions of Lemma~\ref{L4.1.1} are satisfied, so $\fg=3\fme$.
\end{proof}

\sssec{The W-regrading $\fme(; p_2)$ of $\fme$}\label{ssMeP2}
Consider a~W-regrading of $\fg=\fme$ which we denote $\fme(;p_2)$: we set
\be\label{gr1}
 \begin{array}{l}
\deg p_1=2, \ \ \deg q_1=-1, \ \ \deg p_2=1, \ \ \deg q_2=0, \ \ \deg t=1\\
\gr(f)=\deg f-1.
 \end{array}
 \ee

If $\deg_t f=k>0$ and $f=t^{(k)} g(p_1,p_2,q_2,q_1)+h$, where $\deg_t h<k$, we call $t^{(k)}g$ the \textit{leading term} of $f$.

Conditions~\eqref{fpq} imply that if~$f$ is independent of~$t$, then $\deg_{q_1}f\le 1$,
whereas if~$f$ depends on~$t$, and~$t^kg$ is the leading term of~$f$, then $\deg_{q_1}g\le 1$.
Therefore, only functions of the form $q_1f(q_1,q_2)$ have negative degree relative this grading.
Thanks to~\eqref{fpq} there is only one (up to a~factor) such function: $q_1$.
Therefore, the depth of this grading is equal to 2 and $\fme(;p_2)_{-2}=\Span(q_1)$.

To describe the other components we need all 5 equations that single out $\fme\subset \fk(5)$.
The 3 equations that complement~\ref{eqme} are as follows:
\begin{equation}\label{eqme2}
 \begin{array}{l}
 ((\widetilde K_{p_2})^2-\widetilde K_{p_1}\widetilde K_{q_2})(f)=0,\\
 ((\widetilde K_{q_2})^2-2\widetilde K_{q_1}\widetilde K_{p_2})(f)=0,\\
 (2\widetilde K_{p_1}\widetilde K_{q_1}-\widetilde K_{p_2}\widetilde K_{q_2}-\widetilde K_1)(f)=0.
 \end{array}
\end{equation}

On the space of functions independent of~$t$, these conditions are as follows:
\begin{equation}\label{fpq2}
 \begin{array}{l}
 \frac{\del^2 f}{(\del q_2)^2}+\frac{\del^2 f}{\del q_1\del p_2}=0, \\
 \frac{\del^2 f}{(\del p_2)^2}+2\frac{\del^2 f}{\del p_1\del q_2}=0, \\
 2\frac{\del^2 f}{\del p_1 \del q_1}-\frac{\del^2 f}{\del p_2\del q_2}=0.
 \end{array}
\end{equation}

Let us describe $\fme(;p_2)_{-1}$.
This component consists of degree-0 functions in grading~\eqref{gr1}.
Taking~\eqref{fpq} into account we see that the space of such functions independent of~$t$
is spanned by the functions of the form $\vf(q_2)$ and $p_2q_1$.
The first of equations~\eqref{fpq2} brings us $q_2^{(2)}-p_2q_1$ and equality
\[
 \textstyle \frac{\del^3 f}{(\del q_2)^3}=\frac{\del^3 f}{\del q_2\del q_1\del p_2}=0. 
\]
Thus, of functions independent of~$t$, only three function: $1$, $q_2$ and $q_2^{(2)}-p_2q_1$ lie in $\fme(;p_2)_{-1}$.

If $f$ depends on~$t$ and has degree 0 relative grading~\eqref{gr1}, then its leading coefficient should ---
taking~\eqref{fpq} into account --- be proportional to $tq_1$. 

As a~result,
\[
 \fme(;p_2)_{-1}=\Span(1,q_2, G, F_{q_1}),
\]
where $\{1,F_{q_1}\}=2q_1$, and $\{G,q_2\}=q_1$.

Let us now describe $\fme(;p_2)_0$.
It contains functions of degree $1$ relative grading~\eqref{gr1}.
Of functions independent of~$t$, it can contain --- taking~\eqref{fpq} into account --- only $p_2f(q_2)$, $p_1q_1$ and $p_2^{(2)}q_1$.
The 3rd of equations~\eqref{fpq2} yields $p_1q_1+2p_2q_2$ and equality
\[
 \textstyle \frac{\del^3 f}{\del p_2(\del q_2)^2}=2\frac{\del^3 f}{\del p_1\del q_1\del q_2}=0.
\]
The 2nd of equations~\eqref{fpq2} means that
\[
 \textstyle \frac{\del^3 f}{(\del p_2)^2\del q_1}=-2\frac{\del^3 f}{\del p_1\del q_2\del q_1}=0.
\]

Thus, $\fme(;p_2)_0$ contains 2 functions independent of~$t$. namely, $p_2$ and $p_1q_1+2p_2q_2$,
whereas the leading elements of~$f$ depending on~$t$ are of the form
\[
 t, \quad tq_2, \quad t(q_2^{(2)}-p_2q_1)\quad t^{(2)}q_1.
\]

The functions having the first 3 leading elements are
\[
 t, \; F_{q_2}, \; F_G=t(q_2^{(2)}-p_2q_1)+2p_1p_2q_1^{(2)}-p_1q_1q_2^{(2)}+2p_2^{(2)}q_1q_2
\]
Finally, the leading element of $\{ F_{q_2}, F_G\}$ is $-2t^{(2)}q_1$, as is easy to verify.
Therefore,
\[
 \fme(;p_2)_0=\Span(p_2, p_1q_1+2p_2q_2, t, F_{q_2}, F_G, \{ F_{q_2}, F_G\}).
\]

The bracketing with $p_2\in \fme(;p_2)_0$ is as follows:
\begin{equation}\label{adp2}
 \ad_{p_2}: F_{q_1}\mapsto 2D\mapsto -2q_2\mapsto 2.
\end{equation}

The bracketing with $F_{q_2}\in \fme(;p_2)_0$ is (up to coefficients) as follows:
\[
 \ad_{F_{q_2}}: 1\mapsto q_2\mapsto D\mapsto F_{q_1}.
\]

The bracketing with $F_G$ and $\{ F_{q_2}, F_G\}$ is as follows:
\[
 \begin{array}{l}
 \ad_{F_G}: 1\mapsto 2D\mapsto 0, \; q_2\mapsto -F_{q_1}\mapsto 0; \\ 
 \ad_{\{F_{q_2}, F_G\}}:
 1\mapsto F_{q_1}, \ \ q_2\mapsto 0, \ \ D\mapsto 0, \ \ F_{q_1}\mapsto 0.
 \end{array}
\]

Finally, the element $-2t-p_1q_1-2p_2q_2$ acts as a~grading operator on $\fme(;p_2)$.

Thus, the algebras~$\fme(;p_2)$ and~$\fme$ are isomorphic as graded algebras, and this regrading is an automorphism of $\fme$.

\parbegin{Lemma}\label{L2}
If $\fg$ is a~W-regrading of $\fme$ such that $p_2\in\fg_0$, then $\fg=\fme(;p_2)$.
\end{Lemma}

\begin{proof}
 Since $\{G,p_1\}=p_2$ and $p_1\notin \fg_0$ as follows from Lemma~\ref{L1}, it follows that $G\notin \fg_0$.
 Since $\{F_{q_1},p_2\}=2G$, the elements $F_{q_1}$ and~$G$ belong to the same component of~$\fg$, i.e., $F_{q_1}\notin \fg_0$.
 Besides, $\{p_2, q_2\}=-1$, and since $1\in \fg_-$, it follows that $q_2\in \fg_-$ as well.
 Therefore, $\fg_0\cap \fme(;p_2)_-=0$ and $q_2\in \fme(;p_2)_{-1}\cap \fg_-$, i.e.,
 the conditions of Lemma~\ref{L4.1.1} are satisfied, and hence $\fg=\fme(;p_2)$.
\end{proof}

\parbegin{Lemma}\label{L3}
There are no W-regradings $\fg$ of $\fme$ such that $q_2\in\fg_0$.
\end{Lemma}

\begin{proof}
 As we have already observed, if~$\fg$ is a~$W$-regrading of $\fme$, then $q_1\in\fg_-$ and $1\in\fg_-$.
 If $q_2\in \fg_0$, then~$q_1$ and~$G$ must belong to the same component $\fg_k$, where $k<0$, because $\{q_2,G\}=-q_1$.
 The elements $p_2$ and $1$ must also lie in one component~$\fg_m$, where $m<0$, because $\{q_2,p_2\}=1$.
 But $\{G,p_2\}=q_2$, implying $k+m=0$.
 This is a~contradiction.
\end{proof}
Theorem~\ref{ThMe} is proved.
\end{proof} 

\section{$p=3$, the Frank and Ermolaev algebras, see \cite{Sk,GL}}\label{SFrandEr}

\ssec{The Frank algebra}\label{Fr} The Lie algebra~$\fk(3;\un)$ is considered preserving the distribution singled out by the form
\begin{equation}\label{cont}
 dt -pdq+qdp.
\end{equation}
The Frank algebra~$\ffr(n)$ has the same non-positive part as $\fk(3;\un)$, where $\un=(n,1,1)$,
but its degree-1 component is smaller.
All its elements of positive degree are explicitly written in eq.~(83) in~\cite{GL} in terms of the functions
generating contact vector fields;
here, we reproduce them having corrected typos. For $i>0$, we have:
\begin{equation}\label{fr} 
\begin{array}{ll}
\text{for }\ffr_{2i-1}:&\left\{p{}^{(2)}q t^{(i-1)}-pt^{(i)},\ \
q{}^{(2)}p t^{(i-1)}+q t^{(i)}\right\},\\

\text{for }\ffr_{2i}:&\left\{p{}^{(2)}t^{(i)}, \ \
pqt^{(i)}, \ \ q{}^{(2)}t^{(i)},\ \ p{}^{(2)}q{}^{(2)}t^{(i-1)}-t^{(i+1)} \right\}
\text{~~if $1\leq i\leq n-1$}\\

\text{for }\ffr_{2(n+1)}:&\left\{p{}^{(2)}t^{(n+1)}, \ \
pqt^{(n+1)}, \ \ q{}^{(2)}t^{(n+1)} \right\}.\\\end{array}
\end{equation}
For $\fk(3;\un)$, the heights of the indeterminates are not constrained, for $\ffr(n)$, the heights of $p$ and $q$ are constrained.

\sssec{The W-regrading $\widetilde \ffr$ of the Frank algebra $\ffr$} 
\label{WgradFr} From eq.~\eqref{cont} we see that there are two symmetric possibilities for a~W-regrading: either $\deg p=0$ or $\deg q=0$.
Let, for definiteness, $\deg p=0$.
Then, from eq.~\eqref{cont} we see that $\deg t=\deg q=1$, and we get the following basis of non-positive components of the W-regraded algebra of depth 1 we denote $\widetilde \ffr$:
\begin{equation}\label{tilfr}
 \begin{array}{ll}
 \widetilde \ffr_{-1}:&\left\{1,\ \ p, \ \ p^{(2)}\right\},\\ 
 \widetilde \ffr_{0}:&\left\{q,\ \ qp,\ \ qp{}^{(2)}-pt, \ \ t,\ \ p{}^{(2)}t \right\}.\\
 \end{array}
\end{equation}

\paragraph{Why $\widetilde \ffr_{-1}$ is an irreducible $\widetilde \ffr_{0}$-module}
Let us look how~$\widetilde\ffr_0$ acts on~$\widetilde\ffr_{-1}$:
\[
 \begin{array}{lll}
 \{q,f(p)\}=\frac{\del f}{\del_p}&\Lra &q \lra \del_p, \\
 \{qp,f(p)\}=p\frac{\del f}{\del_p}&\Lra &qp \lra p\del_p.
 \end{array}
\]
The elements~$qp$ and~$t$ diagonally act on~$\widetilde\ffr_{-1}$ and in the basis $(1,p,p^{(2)})$ have the form
\[
 \begin{array}{l}
 qp=\diag(0,1,2), \\ t=\diag(-2,-1,0)\Lra t-qp=\diag(1,1,1)\Lra t\lra p\del_p+1=p\del_p+\Div(p\del_p).
 \end{array}
\]
The element $pt-p^{(2)}q$ acts as follows:
\[
 \ad_{pt-p^{(2)}q}: 1\mapsto p,\; p\mapsto 0,\; p^{(2)} \mapsto 0 \Lra pt-p^{(2)}q \lra p^{(2)}\del_p +\Div(p^{(2)}\del_p).
\]
Finally, $p^{(2)}t$ sends $1\mapsto p^{(2)}$, it sends~$p$ and~$p^{(2)}$ to~$0$, so it works as the
multiplication by~$p^{(2)}$.

How could this be that we have operators of multiplication by~$1$ and by~$p^{(2)}$, but no multiplication by~$p$?
Let us look at the bracket of~$p^{(2)}t$ with~$q$,
and how this bracket can be described in terms of a~differential operator~$D$ corresponding to~$p^{(2)}t$:
\[
 \{q,p^{(2)}t\}=-qp^{(2)}+pt\lra [\del_p,D]=p^{(2)}\del_p+p\Lra D=p^{(3)}\del_p+p^{(2)}=p^{(3)}\del_p+\Div(p^{(3)}\del_p).
\]
Thus, $\widetilde \ffr_{-1}$ is an~irreducible $\widetilde \ffr_{0}$-module, where
\[
 \widetilde \ffr_{0}\simeq \fvect(1; 1)\oplus \Span(p^{(3)}\del_p , z).
\]
Here~$z$ is the center of~$\widetilde \ffr_{0}$ and grading operator on the whole~$\widetilde \ffr$
while the vector field~$D$ generated by the ``virtual" or ``ghost"
(non-existing for $\un=\One$, but nevertheless influencing the existing elements, see~\cite{LeP} and~\cite[Remark 2.8.1]{BGLLS2})
element $p^{(3)}\del_p$ establishes irreducibility of the $\widetilde \ffr_{0}$-module $\widetilde \ffr_{-1}$
considered by the $\fvect(1; 1)$-part of $\widetilde \ffr_{0}$ as the space of volume forms $\cO(1; 1)_{\Div}$.

\sssbegin{Theorem}\label{TmFr}
All W-gradings of Frank algebra are $1)$ $\widetilde \ffr$ of depth $1$ and $2)$ $\ffr$ which preserves a~contact distribution.
\end{Theorem}

\ssec{Ermolaev algebra}\label{Er}
Recall that as $\Zee/2$-graded Lie
algebra, it is defined as follows in terms of functions in the $u_i$, and the volume element $v:=\vvol$:
\begin{equation}
\label{er1}
\renewcommand{\arraystretch}{1.4}
\begin{array}{l}
\fer(\un):=\fvect(2; \un)\oplus
\cO(2;\un)_{\Div},\\
\fer(\un)^{(1)}= \fvect(2; \un)\oplus
\cO'(2;\un)_{\Div}. 
\end{array}\end{equation}

To define the bracket, recall that $v^2=v^{-1}$ and observe that 
\begin{equation}\label{dudu/v2}
(du_i)v^{-1}=\sign(ij)\del_j\;\text{ for any permutation $(ij)$ of
$(12)$}.
\end{equation}
For any $f v, g v \in\cO(2; \un)_{\Div}$, set
\begin{equation}
\label{er2} {}[fv, gv]:=(fdg-gdf)v^{-1}= (f\del_2(g)-g\del
_2(f))\del_1+(g\del_1(f)-f\del _1(g))\del_2 \text{};
\end{equation}
define the other products canonically. Define the $\Zee$-grading of
$\fer(\un)$ as induced by the standard $\Zee$-grading of
$\cO(2;\un)$:
\begin{equation}
\label{er3} \fer(\un)_i:=\fvect(2; \un)_i\oplus
\left(\cO(2;\un)_{\Div}\right)_{i+1}\quad \text{for any
$i\geq -1$}.
\end{equation}
Then, $\fer(\widetilde\un)$ is defined as the Cartan prolong of the
following non-positive part:
\begin{equation}
\label{er}
\renewcommand{\arraystretch}{1.4}
\begin{array}{lll}
\fer_{-1}&=&\Span_\Kee(\del_1, \del_2;\; \vvol);\\
\fer_{0}&=&\Span_\Kee(u_iD_j;\; u_k\vvol\mid i, j, k=1, 2).
\end{array}
\end{equation}

Constructing terms of higher degree we see that only the indeterminate $x_3$ corresponding to $\vvol$ is constrained and no W-regradings are possible.
We have arrived at the following statement. 

\sssbegin{Theorem}\label{TmEr}
The only W-grading of Ermolaev algebra is of depth $1$; hence, no non-integrable distributions.
\end{Theorem}

\section{$p=3$, the Skryabin algebras, see \cite{Sk29, Sk,GL}}\label{SSkr}
In~\cite{GL}, all Skryabin algebras are realized as (complete or partial) Cartan prolongs of the non-positive part of certain simple Lie algebras with Cartan matrix considered with a~certain $\Zee$-grading of their Chevalley generators; by symmetry only degrees of positive generators are indicated. 

\ssec{The deep Skryabin algebra $\fg=\mathfrak{dy}(10)$}\label{dy}
Consider the Brown algebra $2\fbr(3)$, see~\cite{GL},
\begin{equation}\label{001}
\text{with Cartan matrix \footnotesize $\begin{pmatrix}
 2 & -1 & 0 \\ -2 & 2 & -1 \\ 0 & -1 & \ev\\
 \end{pmatrix}$ \normalsize and grading $(0,0,1)$.}
\end{equation}

For $\fg=\mathfrak{dy}(10)$ or its special subalgebra $\fs\mathfrak{dy}(10)$,
let $\fg_-\simeq 2\fbr(3)_-$ be realized by vector fields as follows:
\footnotesize
\begin{equation}\label{realdy2}
\renewcommand{\arraystretch}{1.4}
\begin{tabular}{|l|l|}
 \hline $\fg_{-1}$
 &$X_1 = 2 \del_1 + 2u_2\del_4 + 2u_3\del_5 + 2u_2^{(2)}u_3\del_9 + (u_2u_3 + 2u_6)\del_7 + (2u_3u_4 + u_7)\del_8$,\\
 &$X_2 = 2 \del_2 + u_5\del_7 + u_7\del_9,\qquad
 X_3 = 2 \del_3 + u_2\del_6 + 2u_4\del_7 + u_2u_4\del_9 + (u_2 u_5+u_7)\del_{10}$\\
 \hline $\fg_{-2}$
 &$X_4:= \del_{4} + u_{5} \del_{8} + u_{6} \del_{9}$,\qquad
 $X_5:= \del_{5} + u_{6}\del_{10}$,\qquad
 $X_6:= \del_{6}$\\
 \hline
 $\fg_{-3}$
 &$X_7:=\del_{7}$\\
 \hline $\fg_{-4}$
 &$X_8:=\del_{8},\quad X_9:=\del_{9},\quad
 X_{10}:=\del_{10}$\\
\hline
\end{tabular}
\end{equation}
\normalsize

Then, see~\cite{GL}, we have
$2\fbr(3)_0=\mathfrak{dy}(10)_0=\fgl(3)$ and $\fs\mathfrak{dy}(10)_0=\fsl(3)$:
\footnotesize
\begin{equation}\label{realdy}
\renewcommand{\arraystretch}{1.4}
\begin{tabular}{|ll|}
\hline
$x_1 =$&$ u_3\del_2 + 2u_3^{(2)}\del_6 + (2u_4u_3^{(2)} + u_{10})\del_9 + u_5\del_4 + 2u_3^{(2)}u_5\del_{10} + u_5^{(2)}\del_8$,\\
$x_2 =$&$ 2u_2\del_1 + 2u_2^{(2)}\del_4 + 2u_2^{(2)}u_3\del_7 + (2u_2u_3 + 2u_6)\del_5 + (2u_2u_3u_4 + 2u_9)\del_8 + 2u_6^{(2)}\del_{10}$,\\
$x_3 =$&$ [x_1,x_2] = 2u_3\del_1 + 2u_3^{(2)}\del_5 + u_2u_3^{(2)}\del_7 + (2u_2^{(2)}u_3^{(2)} + u_6^{(2)})\del_9 + (2u_4u_3^{(2)} + u_5u_6 + 2u_{10})\del_8 + u_6\del_4$,\\
$y_1 =$&$ u_2\del_3 + 2u_2^{(2)}\del_6 + u_4\del_5 + 2u_2^{(2)}u_4\del_9 + u_4^{(2)}\del_8 + (2u_5u_2^{(2)} + u_9)\del_{10}$,\\
$y_2 =$&$ 2u_1\del_2 + 2u_1^{(2)}\del_4 + 2u_1^{(2)}u_3\del_7 + (u_1u_3 + 2u_5)\del_6 + (u_2u_3u_1^{(2)} + u_3u_4u_1 + 2u_8)\del_9$\\
 &$ + (u_1^{(2)}u_3^{(2)} + u_1u_5u_3 + 2u_5^{(2)})\del_{10}$,\\
$y_3 =$&$ [y_1,y_2] = u_1\del_3 + u_1^{(2)}\del_5 + 2u_1^{(2)}u_2\del_7 + (u_1^{(2)}u_2^{(2)} + 2u_4^{(2)})\del_9 + 2u_4\del_6
 + u_1^{(2)}u_4)\del_8 + (2 u_4 u_5+u_8)\del_{10}$,\\
$h_1 =$&$ [x_1,y_1] = 2u_2\del_2 + u_3\del_3 + 2u_4\del_4 + u_5\del_5 + 2u_9\del_9 + u_{10}\del_{10}$,\\
$h_2 =$&$ [x_2,y_2] = 2u_1\del_1 + u_2\del_2 + 2u_5\del_5 + u_6\del_6 + 2u_8\del_8 + u_9\del_9$,\\
$h_3 =$&$2u_2\del_2 + 2u_3\del_3 + 2u_4\del_4 + 2u_5\del_5 + u_6\del_6+u_7\del_7+u_8\del_8$\\
\hline
\end{tabular}
\end{equation}
\normalsize

The vector fields~$Y_i$ commuting with~$\fg_-$ and their weights are as follows
\footnotesize{
\[\begin{tabular}{ll}
\hline
$\fg_{-1}$&$\begin{array}{ll}
 (0,1,0) & \del_1+u_3\del_5+ u_2u_3\del_7+2 u_2u_3^{(2)} \del_{10}+2 u_3 u_4 \del_8+ u_2 u_6 \del_9,\\
 (1,-1,1)& \del_2+2 u_1 \del_4+u_3 \del_6+2 u_1 u_3 \del_7+2 u_1^{(2)} u_3 \del_8+ u_1 u_2 u_3 \del_9 + u_1u_6 \del_9,\\
 (-1,0,1)& \del_3+2 u_1 u_2^{(2)} \del_9 + u_2u_4\del_9 + u_2u_5\del_{10},\\
 \end{array}$\\
 \hline
$\fg_{-2}$&$\begin{array}{ll}
 (1,0,1)& \del_4+2 u_3 \del_7+2 u_1 u_3 \del_8+2 u_3^{(2)} \del_{10},\\
 (-1,1,1)& \del_5+u_2 \del_7+u_2^{(2)} \del_9+u_2 u_3 \del_{10}+2 u_4 \del_8,\\
 (0,-1,2)& \del_6+2 u_1 \del_7+2 u_1^{(2)} \del_8+2 u_4 \del_9+2 u_5 \del_{10},\\
\end{array}$\\
 \hline
$\fg_{-3}$&$\begin{array}{ll} (0,0,2)& \del_7+u_1 \del_8+u_2 \del_9+u_3 \del_{10},\\
 \end{array}$\\
 \hline
$\fg_{-4}$&$\begin{array}{ll} (0,1,2)& \del_8,\\
 (1,-1,3)& \del_9,\\
 (-1,0,3)& \del_{10}. \end{array}$\\
 \hline
\end{tabular}\]
}
\normalsize

The corresponding dual forms $\Theta^i$ are as follows
\begin{gather*}
\Theta^1 = du_1,\qquad
\Theta^2 = du_2,\qquad
\Theta^3 = du_3,\\
\Theta^4 = du_4+u_1du_2,\qquad
\Theta^5 = du_5+2 u_3du_1,\qquad
\Theta^6 = du_6+2 u_3du_2,\\
\Theta^7 = du_7+u_1du_6+2 u_2du_5+u_3du_4+ u_1 u_3 du_2,\\
\Theta^8 = du_8+2 u_1du_7+2 u_1^{(2)}du_6+\left(u_1 u_2\right)du_5+u_4du_5,\\
\Theta^9 = du_9+ 2u_2 du_7 + (2u_1u_2 + u_4)du_6 + u_2^{(2)}du_5 + 2u_2u_3du_4 + (u_1u_2^{(2)} + 2u_2u_4)du_3\\
{} + (u_1u_2u_3 + 2u_3u_4 + 2u_1u_6)du_2 + (u_2^{(2)} u_3 + 2 u_2u_6)du_1, \\
\Theta^{10} = du_{10}+2 u_3du_7+2 u_1 u_3du_6+2 u_3^{(2)}du_4+2 u_1 u_3^{(2)}du_2+u_5du_6+2 u_2 u_5du_3+2 u_3 u_5 du_2.
\end{gather*}

\sssec{W-regradings}\label{WgradDY}
Let $d_i:=\deg u_i$. The indeterminates with constrained heights are the $u_i$ for $i=1, \dots, 7$, see~\cite{GL}.
Since the summands of the $\Theta^i$ must be of the same degree, we get 
\[
\begin{array}{lll}
d_4=d_1+d_2,&
d_5=d_1+d_3,&
d_6=d_2+d_3,\\
d_7=d_1+d_6=d_2+d_5=d_3+d_4,&
d_8=d_1+d_7=d_4+d_5, &
d_9=d_2+d_7=d_4+d_6,\\
d_{10}=d_3+d_7=d_5+d_6
\end{array}
\]
or equivalently
\begin{equation}\label{solDY}
\begin{array}{lll}
d_4=d_1+d_2,&
d_5=d_1+d_3,&
d_6=d_2+d_3,\\
d_7=d_1+d_2+d_3,&
d_8=2d_1+d_2+d_3,&
d_9=d_1+2d_2+d_3,\\
d_{10}=d_1+d_2+2d_3.
\end{array}
\end{equation}

\paragraph{The W-regradings~$\fdy(9)$ of~$\fdy(10)$.}
Let $d_1= d_2=0$. Denote this regrading of~$\fdy(10)$ by~$\fdy(9)$.
Then, relations~\eqref{solDY} imply
\[
\begin{array}{l}
d_1= d_2= d_4=0,\\
d_3=d_5=d_6=d_7=d_8=d_9\text{~(we set:) $=1$};\ \
d_{10}=2.\\
\end{array}
\]
For a~basis of the negative part of $\fg=\fdy(9)$ in terms of functions generating contact vector fields we take
(the degree in~$\fdy(10)$, the function~$f$ generating a~contact vector field corresponding to the vector field in~$\fdy(10)$)
\tiny
\be\label{pAndQ}
 \begin{array}{c|l|l}
 \deg & f &\text{vector field in~$\fdy(10)$}\\
 \hline
 -4 & 1 & \del_{10} \\
 -4 & p_1 & \del_8 \\
 -4 & p_2 & \del_9 \\
 -3 & p_3 & \del_7 \\
 -2 & p_4 & \del_5 + u_6 \del_{10} \\
 0 & q_1 & 2u_1\del_3 + 2u_1^{(2)}\del_5 + u_1^{(2)}u_2\del_7 + (2 u_1^{(2)} u_2^{(2)} + u_4^{(2)})\del_9 + u_4\del_6 + 2u_1^{(2)}u_4\del_8 + (u_4 u_5+2 u_8)\del_{10} \\
 0 & q_2 & 2u_2\del_3 + u_2^{(2)}\del_6 + 2u_4\del_5 + u_2^{(2)}u_4\del_9 + 2u_4^{(2)}\del_8 + (u_5 u_2^{(2)}+2 u_9)\del_{10} \\
 -1 & q_3 & \del_3+ 2u_2\del_6 + u_4\del_7 + 2u_2u_4\del_9 + (2 u_2 u_5+2 u_7)\del_{10} \\
 -2 & q_4 & \del_6 \\
\end{array} 
\ee
\normalsize

For a~basis of $\fdy(9)_0$ we take the elements whose leading parts are listed in table~\eqref{dy(9)_01}, 
where $\deg$ is the degree in~$\fdy(10)$ with respect to $-(s_8+s_9)$, the new name $s_i$, the generating
function, the vector field in~$\fdy(10)$.
The component $\fs\fdy(9)_0$ does not contain $h_3 = p_1q_1 + p_3q_3 + 2p_4q_4$. 
\tiny
\be\label{dy(9)_01}
 \begin{array}{c|l|l|l}
 \deg & s_i & f &\text{vector field in terms of the }u_i\\
 \hline
 4 & s_1 & q_2^{(2)}
 & 2u_2^{(2)} u_4^{(2)}\del_4 + 2u_2^{(2)}u_3u_4^{(2)}\del_7 + 2u_1u_4u_2^{(2)} + 2u_4^{(2)}u_2\del_1 + (2u_3u_4u_2^{(2)} +\dots\\ 
 2 & s_2 & p_4 q_2
 & (2 u_3 u_2^{(2)}+u_6 u_2)\del_3 + 2u_2^{(2)} u_4\del_4 + 2u_2^{(2)}u_3u_4\del_7 + (2u_1u_2^{(2)} + 2u_4u_2)\del_1+\dots\\
 0 & s_3 & 2 p_1 q_2+2 p_4^{(2)}
 & 2 u_2\del_1 + 2u_2^{(2)}\del_4 + 2u_2^{(2)}u_3\del_7 + (2u_2u_3 + 2u_6)\del_5 + (2u_2u_3u_4 + 2u_9)\del_8 + 2u_6^{(2)}\del_{10} \\
 1 & s_4 & 2 p_3 q_2+p_4 q_3
 & (u_1u_2 + u_4)\del_1 + (u_2u_3 + 2u_6)\del_3 + u_2u_4\del_4 +2u_3u_4^{(2)}\del_8 +\dots\\
 -1 & s_5 & p_1 q_3+2 p_3 p_4
 & 2 \del_1 + 2u_2\del_4 + 2u_3\del_5 + 2u_2^{(2)}u_3\del_9 + (u_2 u_3+2 u_6)\del_7 + (2 u_3 u_4+u_7)\del_8 \\
 3 & s_6 & q_2 q_3
 & u_2u_4^{(2)}\del_4 + (u_4^{(2)} + u_1u_2u_4)\del_1 + (u_2u_3u_4 + u_2u_7 + 2u_9)\del_3 + (2u_2u_3u_4^{(2)} +\dots \\ 
 -2 & s_7 & p_3^{(2)} - p_1 q_4 - p_2 p_4
 & 2 \del_4 + 2 u_5\del_8 + 2 u_6\del_9 \\
 0 & s_8 & p_1 q_1-p_2 q_2+p_4 q_4
 & u_1\del_1 + 2u_2\del_2 + u_5\del_5 + 2u_6\del_6 + u_8\del_8 + 2u_9\del_9 \\
 0 & s_9 & p_1 q_1+p_3 q_3-p_4 q_4
 & 2 u_2\del_2 + 2 u_3\del_3 + 2 u_4\del_4 + 2 u_5\del_5 + u_6\del_6 + u_7\del_7 + u_8\del_8 \\
 4 & s_{10} & q_1 q_2
 & u_1^{(2)}u_2^{(2)} u_4\del_4 + 2u_1u_4^{(2)}\del_1 + (u_1 u_4 u_2^{(2)}+u_4^{(2)} u_2)\del_2 + (u_1 u_2 u_3 u_4+u_2 u_8+u_1 u_9)\del_3+\dots\\
 2 & s_{11} & 2 q_3^{(2)}+q_2 q_4
 & 2u_1u_4\del_1 + 2u_4^{(2)}\del_4 + (2 u_3 u_4+u_2 u_5+u_7)\del_3 + u_4 u_5\del_5 + u_4^{(2)}u_5\del_8+\dots\\

2 & s_{12} & p_4 q_1+q_2 q_4
 & u_1^{(2)} u_2^{(2)}\del_4 + (u_1 u_2 u_3+u_2 u_5+u_1 u_6)\del_3 + (u_2 u_3 u_1^{(2)}+u_6 u_1^{(2)}+u_4 u_5+u_8)\del_5+\dots\\
 1 & s_{13} & p_3 q_1+2 q_3 q_4
 & u_1^{(2)}u_2\del_4+(u_1 u_2+u_4)\del_2 + u_1^{(2)} u_3\del_5+(u_1 u_3+u_5)\del_3+(2 u_1 u_2 u_3+2 u_2 u_5+u_7)\del_6+\dots\\
 -1 & s_{14} & p_2 q_3+p_3 q_4
 & 2 \del_2+u_5\del_7+u_7\del_9 \\
 3 & s_{15} & q_1 q_3
 & 2u_1^{(2)} u_2 u_4\del_4 + (u_1 u_2 u_3 u_4+2 u_7 u_4+u_2 u_8)\del_6+\dots\\

 0 & s_{16} & 2 p_2 q_1+q_4^{(2)}
 & 2 u_1\del_2 + 2u_1^{(2)}\del_4+ 2u_1^{(2)}u_3\del_7 + (u_1 u_3+2 u_5)\del_6+(u_2 u_3 u_1^{(2)}+u_3 u_4 u_1+2 u_8)\del_9+\dots\\
 2 & s_{17} & 2 q_1 q_4
 & 2u_1u_4\del_2 + 2u_1^{(2)}u_4\del_4 + (u_1 u_3 u_4+u_5 u_4+u_8)\del_6+\dots\\

 
 4 & s_{18} & 2 q_1^{(2)}
 & 2u_1u_4^{(2)}\del_2 + 2u_1^{(2)}u_4^{(2)}\del_4+(u_1 u_3 u_4^{(2)}+u_8 u_4)\del_6+\dots\\
 0 & s_{19} & t & u_3\del_3+u_5\del_5+u_6\del_6+u_7\del_7+u_8\del_8+u_9\del_9+2 u_{10}\del_{10} \\
\end{array}
\ee
\normalsize

\paragraph{Description of $\fdy(9)_0$}
Actually, we need the fourth column in table~\eqref{dy(9)_01} obtained with the help of \textit{SuperLie}
only to compose columns one and three. All calculations are then performed manually in terms of the contact bracket.
Clearly, $t$ is a~central element of~$\fdy(9)_0$.
Let $\fh$ be the subalgebra of $\fdy(9)_0$ spanned by the generating functions independent of~$t$.
Then, the grading of~$\fdy(10)$ induces the grading $\fh=\oplus_{i=-2}^4 \ \fh_i$:
\be\label{g}
\begin{array}{|l|l|l|l||l|l|l|l|}
\hline
 & \fh_{-2} & \fh_{-1} & \fh_0 & \fh_1 & \fh_2 & \fh_3 & \fh_4 \\
 \hline\hline
 \dim(\fh_i) & 1 & 2 & 4 & 2 & 4 & 2 & 3 \\
 \hline
 \text{Basis of }\fh_i & s_7 & s_5, s_{14} & s_3,s_8,s_9,s_{16} & s_4, s_{13} & s_2,s_{11},s_{12},s_{17} & s_6,s_{15} & s_1, s_{10},s_{18} \\
 \hline
\end{array} 
\ee

The element $s_8+s_9\in \fh_0$ is central in~$\fh_0$ and acts on $\fh$ as a~2(grading operator), whereas
$\Span(s_3,s_8,s_{16})\simeq\fsp(2)$.
Let $\fk(3;\One)$ be realized in terms of the contact form $dT-PdQ+QdP$.
We see that $\fh$ is a~subalgebra of $\fk(3;\One)$.

To describe the positive part of $\fh$, observe that $\dim\fk(3;\One)_4=4$ and looking at table~\eqref{g} we deduce that, as $\fk(3;\One)_0$-module, $\fk(3;\One)_4\simeq \Kee TP^{(2)}Q^{(2)}\oplus \fh_4$. Brackets of $\fh_{-1}$ with $\fh_{4}$ yield
\[
\begin{array}{l}
\fh_4\simeq\Span(T^{(2)}P^{(2)}, \ T^{(2)}PQ, \ T^{(2)}Q^{(2)}),\\
\fh_3\simeq\Span(T^{(2)}P-P^{(2)}Q, \ T^{(2)}Q+Q^{(2)}P),\\
\fh_2\simeq 
\Span(TP^{(2)},\ TPQ, \ TQ^{(2)})\oplus \Kee (T^{(2)}-P^2Q^2),\\
\fh_1\simeq\Span(TP-P^2Q, TQ+PQ^2).\\
\end{array}
\]
Thus, $\fh\simeq\ffr(1)$, the Frank subalgebra of $\fk(3;\One)$, see its basis~\eqref{fr}; hence, $\fdy(9)_0=\ffr(1)\oplus \Kee t$.

\paragraph{The $W$-regrading~$\fdy(11)$ of~$\fdy(10)$} Set $d_1=0$, $d_2=1$, $d_3=1$. Then, relations~\eqref{solDY} yield 
\[
 \text{$d_4= d_5=1$,\ \
 $d_6= d_7= d_8=2$, \ \
 $d_9= d_{10}=3$.}
\] 
Let us call this regrading $\fg=\fdy(11)$.
\tiny
\[
\renewcommand{\arraystretch}{1.4}
\begin{tabular}{|l|l|}
 \hline
 $\fg_{-3}$ & $\del_9$, \quad $\del_{10}$ \\
 \hline
 $\fg_{-2}$ & $\del_6$, \quad $\del_7$, \quad $\del_{8}$ \\
 \hline
 $\fg_{-1}$
 & $2\del_4 + u_5\del_8 + u_6\del_9$,\quad
 $ \del_5 + u_6\del_{10}$,\quad
 $2\del_2 + u_5\del_7 + u_7\del_9$,\\
 & $2 \del_3 + u_2\del_6 - u_4\del_7 + u_2u_4\del_9 + (u_2u_5 + u_7)\del_{10}$,\\
 & $2u_1\del_2 + 2u_1^{(2)}\del_4 + 2u_1^{(2)}u_3\del_7 + (u_1 u_3+2 u_5)\del_6 + (u_2u_3u_1^{(2)}+u_3u_4u_1+2u_8)\del_9$\\
 &$ + (u_1^{(2)}u_3^{(2)} + u_1u_5u_3 + 2u_5^{(2)})\del_{10}$,\\
 & $u_1\del_3 + u_1^{(2)}\del_5 + 2u_1^{(2)}u_2\del_7 + (u_1^{(2)} u_2^{(2)}+2 u_4^{(2)})\del_9 + 2u_4\del_6 + u_1^{(2)}u_4\del_8 + (2 u_4 u_5+u_8)\del_{10}$\\
 \hline
$\fg_{0}$ & $X_+ := u_2\del_3 + 2u_2^{(2)}\del_6 + u_4\del_5 + 2u_2^{(2)}u_4\del_9 + u_4^{(2)}\del_8 +(2u_5u_2^{(2)}+u_9)\del_{10}$, \\
 & $D := 2u_1u_4\del_2 + 2u_1^{(2)}u_4\del_4 + (u_1 u_3 u_4+u_5 u_4+u_8)\del_6 + (u_2 u_3 u_4 u_1^{(2)}+2 u_2^{(2)} u_5 u_1^{(2)}+2 u_3 u_4^{(2)} u_1+u_4^{(2)} u_5)\del_9$\\
 & ${} + (2u_3 u_4 u_1^{(2)}+u_2 u_5 u_1^{(2)})\del_7 + 2u_1u_5\del_3 + 2u_1^{(2)}u_5\del_5 + 2u_1^{(2)}u_4u_5\del_8 + (u_1^{(2)} u_4 u_3^{(2)}+u_1 u_4 u_5 u_3+2 u_4 u_5^{(2)})\del_{10}$, \\
 & $A_- := u_1^{(2)}u_2\del_4 + (u_1 u_2+u_4)\del_2 + u_1^{(2)}u_3\del_5 + (u_1 u_3+u_5)\del_3 + (2u_1u_2u_3+2u_2u_5+u_7)\del_6 + u_1^{(2)} u_3 u_4\del_8$\\
 & ${}+ (2u_1^{(2)} u_3 u_2^{(2)}+2 u_1 u_3 u_4 u_2+2 u_4 u_5 u_2)\del_9+(u_4 u_5+u_8)\del_7+(2 u_1^{(2)} u_2 u_3^{(2)}+2 u_1 u_2 u_5 u_3+u_2 u_5^{(2)})\del_{10}$, \\
 & $T := 2 u_1\del_1 + u_2\del_2 + 2 u_5\del_5 + u_6\del_6 + 2u_8\del_8 + u_9\del_9$, \\
 & $H:= 2u_2\del_2 + u_3\del_3 + 2u_4\del_4 + u_5\del_5 + 2u_9\del_9 + u_{10}\del_{10}$ \\
 & $E:= 2u_2\del_2 + 2u_3\del_3 + 2u_4\del_4 + 2u_5\del_5 + u_6\del_6 + u_7\del_7 + u_8\del_8$ (Center),\\
 & $A_+ := 2\del_1 + 2u_2\del_4 + 2u_3\del_5 +2u_2^{(2)}u_3\del_9 + (u_2 u_3+2 u_6)\del_7 + (2 u_3 u_4+u_7)\del_8$, \\
 & $X_- := u_3\del_2 + 2u_3^{(2)}\del_6 + (2 u_4 u_3^{(2)}+u_{10})\del_9 + u_5\del_4 + 2u_3^{(2)}u_5\del_{10} + u_5^{(2)}\del_8$. \\
 \hline
\end{tabular}
\]
\normalsize

The multiplication table for $\fg_0$ is as follows
\[
\tiny
\begin{array}{l|llllllll}
 & X_+ & D & A_- & T & H & E & A_+ & X_- \\
\hline
 X_+ & 0 & 0 & 0 & 2 X_+ & 2 X_+ & 0 & 0 & 2 H \\
 D & 0 & 0 & 0 & 2 D & 0 & 0 & 2 A_- & 0 \\
 A_- & 0 & 0 & 0 & A_- & 0 & 0 & E & 0 \\
 T & X_+ & D & 2 A_- & 0 & 0 & 0 & A_+ & 2 X_- \\
 H & X_+ & 0 & 0 & 0 & 0 & 0 & 0 & 2 X_- \\
 E & 0 & 0 & 0 & 0 & 0 & 0 & 0 & 0 \\
 A_+ & 0 & A_- & 2 E & 2 A_+ & 0 & 0 & 0 & 0 \\
 X_- & H & 0 & 0 & X_- & X_- & 0 & 0 & 0 
\end{array}
\]
\normalsize

We see that $S:=\Span(X_-, H, X_+)\simeq \fsl(2)$ and $\Span(A_-, E, A_+)\simeq \fhei(2|0)$.
Recall that $\fhei(2|0)\simeq\fk(3|0)_-$, so the Lie algebra of outer derivations of $\fhei(2|0)$ is $\fk(3|0)_0\simeq\fgl(2)$.
The operator $T$ multiplies $X_\pm$ and $A_\pm$ by $\pm 1$.
Let $\fr:=\fhei(2|0)\inplus \Kee D$.
This $\fr$ is the radical of $\fg_0$, and $\fg_0\simeq(\fr\oplus \fsl(2))\inplus \Kee T$.

\parbegin{Lemma}[On W-gradings of $\fs\fdy(10)$]\label{Lsdy}
 The algebra $\fs\fdy(10)$ has only one W-grading.
\end{Lemma}

\begin{proof}
1) Since $\fs\fdy(9)_-=\fdy(9)_-$, both $\fdy(9)$ and $\fs\fdy(9)$ preserve the same distribution, while (recall notation $\fhei$, see Subsection~\ref{Other})
\begin{equation}\label{sdy_0}
\fs\fdy(9)_0=\Span(s_7; s_5, s_{14}; s_3,s_8,s_{16})=\fhei(2|0)\inplus \fsl(2).
\end{equation}
However, the grading of~$\fs\fdy(9)$ is not a~Weisfeiler one.
Indeed, $\fhei(2|0)=\Span (s_7)\inplus \Span(s_5,s_{14})$, so that
$\Span(s_5,s_{14})\subset \fs\fdy(10)_{-1}$ while
$\Span (s_7) )\subset \fs\fdy(10)_{-2}$, see table~\eqref{dy(9)_01}.

2) The operators $H, T, E\in \fdy(11)_0$ correspond to $h_1, h_2, h_3\in \fdy(10)_0$, respectively.
Having deleted~$E$ (divergence-freeing $\fdy(11)_0$) we must delete~$D$ and~$A_-$ as well.
Then, the action of $\fdy(11)_0\cap\fs\fdy(10)$ is reducible on $\fdy(11)_{-1}$ which is the sum of three 2-dimensional spaces whose bidegrees,
induced by $\fdy(11)$ and $\fdy(10)$, are $(-1,-2)$; $(-1,-1)$ and $(-1,0)$ while $\fsl(2)$, see eq.~\eqref{sdy_0},
has bidegree $(0,0)$, and hence preserves each of these three subspaces.
Since $\text{bideg}(A_+)=(0,1)$, $\text{bideg}(A_-)=(0,-1)$, $\text{bideg}(D)=(0,2)$,
then, without~$D$ and~$A_-$, the subspace of $\text{bideg}(-1,0)$ is invariant.
\end{proof}

\sssbegin{Theorem}\label{TmDy}
The deep Skryabin algebra has $3$ W-gradings, its special subalgebra has one W-grading.
Each of these W-regraded algebras preserves a~non-integrable distribution: 
\begin{equation}\label{GrVdy}
\renewcommand{\arraystretch}{1.4}
\footnotesize
\begin{tabular}{|c|c|}
\hline $\fg$&$v_\fg$\\
\hline
$\fdy(9)$ 
&$(8, 9)$C\\
$\fdy(10)$ and $\fs\fdy(10)$&$(3, 6,7,10)$\\
$\fdy(11)$ &$(6,9,11)$\\
\hline
\end{tabular}
\end{equation}
\end{Theorem}

\begin{proof}
The conditions in Section~\ref{sec:HowTo} imply that any minimal regrading~$\fh=\bigoplus_i\fh_i$ of
$\fg=\fdy(10)$ is defined by two subspaces $\fg_{-1}^0$ and~$\fg_{-1}^{-1}$ such that
\[
 \fg_{-1} = \fg_{-1}^0 \oplus \fg_{-1}^{-1},\quad
 \fg_{-1}^0 \subset \fh_0,\quad
 \fg_{-1}^{-1} \subset \fh_{-1}.
\] 
First of all, note that the elements $s_8$ and $s_9$ diagonally act on~$\fdy(9)$ and should remain in
0th component in any regrading.
Thus, any new regrading should be compatible with the grading defined by
weights with respect to~$s_8$ and~$s_9$. In particular, the subspaces $\fg_{-1}^0$ and~$\fg_{-1}^{-1}$ should
be (sum) of eigenspaces for~$s_8$ and~$s_9$.
The $\fg_0$-module $\fg_{-1}$ splits into three one-dimensional subspaces spanned by $X_1$, $X_2$, and $X_3$,
respectively, with distinct weights.
Therefore, the subspaces~$\fg_{-1}^0$ and~$\fg_{-1}^{-1}$ are defined by subsets $S_0$ and $S_{-1}$ of the set $I := \{1,2,3\}$ such that $S_0\cap S_{-1}=\emptyset$ and
$S_0\cup S_{-1} = I$. Set $d_i = \deg X_i = 0$ for $i\in S_0$ and $d_j = \deg X_j = 1$ for $j\in S_{-1}$.

Note, that $d_1$, $d_2$, and~$d_3$ can not simultaneously be equal to zero. Therefore, taking into account that the indeterminates $u_i$ are on equal footing for $i=1,2,3$ (as well as for $i=4,5,6$ and $i=8,9.10$) we arrive at the following cases: 
\begin{enumerate}[label=\arabic*)]
\item $d_1=1$, $d_2=1$, $d_3=1$: defines the original grading $\fdy(10)$.

\item $d_1=0$, $d_2=0$, $d_3=1$: defines $\fdy(9)$.

\item $d_1=0$, $d_2=1$, $d_3=1$: defines $\fdy(11)$.
\end{enumerate}
For the special subalgebra, see Lemma~\ref{Lsdy}; the theorem is proved.
\end{proof}

\subsection{The middle Skryabin algebra $\fg=\fm\fy(6)$}
Consider $\fo(7)$
\[
\text{with CM \footnotesize $\begin{pmatrix}
 2 & -1 & 0 \\ -1 & 2 & -1 \\ 0 & -2 & 2
 \end{pmatrix}$ \normalsize and grading $(0,0,1)$.}
\]
Then, for a~basis of non-negative elements of $\fg_{\leq 0}\simeq\fo(7)_{\leq 0}$ we can take
\tiny
\begin{equation}\label{MY6}
\renewcommand{\arraystretch}{1.4}
\begin{tabular}{|l|l|}
 \hline
 $\fg_0$
 & $(2,-1,0)\ x_1:= u_3\del_2 + u_5\del_4 + u_3^{(2)}\del_6$,\\
 & $(-1,2,-2)\ x_2:= u_2\del_1 + u_2^{(2)}\del_4 + u_6\del_5$,\\
 & $(1,1,-2)\ x_3:= [x_1,x_2] = u_3\del_1 + (u_2u_3 + 2u_6)\del_4 + u_3^{(2)}\del_5$,\\
 & $(-2,1,0)\ y_1:= u_2\del_3 + u_2^{(2)}\del_6 + u_4\del_5$\\
 & $(1,-2,2)\ y_2:= u_1\del_2 + u_1^{(2)}\del_4 + u_5\del_6$,\\
 & $(-1,-1,2)\ y_3:= [y_1,y_2] = 2u_1\del_3 + 2u_1^{(2)}\del_5 + (2u_1u_2 + u_4)\del_6$\\
& $h_1 := [x_1,y_1] = 2u_2\del_2 + u_3\del_3 + 2 u_4\del_4 + u_5\del_5$,\\
 & $h_2 := [x_2,y_2] = 2u_1\del_1 + u_2\del_2 + 2 u_5\del_5 + u_6\del_6$,\\
 & $h_3 := 2u_1\del_1 + 2 u_4\del_4 + 2 u_5\del_5$,\\
 \hline
 $\fg_{-1}$
 & $(0,1,-2): \del_1 -u_2\del_4 + u_3\del_5,\quad
 (1,-1,0): \del_2 - u_3\del_6,\quad
 (-1,0,0): \del_3$\\ 
 \hline
 $\fg_{-2}$
 & $(1,0,-2): \del_4,\quad (-1,1,-2): \del_5,\qquad (0,-1,0): \del_6$ \\
 \hline
\end{tabular}
\end{equation}
\normalsize
The element $h_1-h_2$ is central in~$\fg_0\simeq\fgl(3)$.

The vector fields~$Y_i$ that commute with~$\fo(7)_-$ and their weights are as follows
\footnotesize
\[\begin{tabular}{ll}
\hline
$\fg_{-1}$&$\begin{array}{ll}
 (0,1,-2)\to{}& \del_1,\\
(1,-1,0)\to{}& \del_2+2 u_1 \del_4,\\
(-1,0,0)\to{}& \del_3+2 u_1 \del_5+2 u_2 \del_6,\\
 \end{array}$\\
 \hline
$\fg_{-2}$&$\begin{array}{ll}
(1,0,-2)\to{}& \del_4,\\
(-1,1,-2)\to{}& \del_5,\\
(0,-1,0)\to{}& \del_6.
\end{array}$\\
 \hline
\end{tabular}\]

\normalsize

The corresponding dual forms $\Theta^i$ are as follows
\be\label{MYtheta}\footnotesize
\begin{array}{lll}
 \Theta^1 = du_1,& \Theta^2 = du_2,& \Theta^3 = du_3,\\
 \Theta^4 = u_1du_2 + du_4, & \Theta^5 = u_1du_3 + du_5, & \Theta^6 = u_2du_3 + du_6. 
\end{array}
\ee

\sssec{W-regradings}\label{WgradMY} Recall that $d_i:=\deg u_i$. The homogeneity of the $\Theta^i$ for $i=4,5,6$ implies 
\be\label{rk3}\footnotesize
\begin{array}{l}
 d_4 =d_1+d_2,\ \
 d_5 =d_1+d_3,\ \
 d_6 =d_2+d_3.
\end{array}
\ee


\paragraph{The W-regrading~$\fmy(7)$ of~$\fmy(6)$}\label{MY(7)_0}
Let $d_1=0$, and hence $d_4 =d_2= d_5 =d_3=1$, and $d_6 =2$, see conditions~\eqref{rk3}.
We denote this regrading by $\fmy(7)$, since $\dim\fmy(7)_-=7$ as we will see.

Let us show how $\fmy(7)\simeq \fmy(6)$ is embedded in $\fk(7;\un)$ and establish a~correspondence between the elements of $\fmy(6)$ and generating functions of $\fk(7;\un)$ in indeterminates $t, p_i, q_i$ for $i=1,2,3$:
\be\label{MY7geq}\tiny
\begin{array}{ll|l}
\hline
 \fmy(7)_{-2}:{}
 & \del_6 & 1\\
 \hline
 \fmy(7)_{-1}:{}
 & \del_2 + 2u_3\del_6 & p_1\\
 & \del_3 & -q_1\\
 & \del_4 & p_2\\
 & u_1\del_3 + 2u_1^{(2)}\del_5 + (2u_4+2u_1u_2)\del_6 & q_2\\
 & u_1\del_2+2u_1^{(2)}\del_4+u_5\del_6 & p_3\\
 & \del_5 & q_3\\
 \hline
 \fmy(7)_0:{}
 & \del_1 + 2 u_2 \del_4 + 2 u_3 \del_5 &p_1q_3+q_1p_2 \\
 & 2 u_2 \del_2 + u_3 \del_3 + 2 u_4 \del_4 + u_5 \del_5 & p_1q_1+p_2q_2+p_3q_3\\
 & 2 u_1 \del_1 + u_2 \del_2 + 2 u_5 \del_5 + u_6 \del_6 & t-p_1q_1+p_2q_2\\
 & 2 u_1 \del_1 + 2 u_4 \del_4 + 2 u_5 \del_5 & p_2q_2-p_3q_3\\
 & 2 u_3 \del_2 + 2 u_5 \del_4 + u_3^{(2)} \del_6 & p_1^{(2)}-p_2p_3\\
 & 2 u_2 \del_3 + 2 u_4 \del_5 + u_2^{(2)} \del_6 & q_2q_3-q_1^{(2)}\\
 & 2 u_1^{(2)} \del_1 + u_4 \del_2 + u_1u_3 \del_3 + 2 u_1^{(2)} u_3 \del_5 + 2 (u_1 u_2 u_3 + u_3 u_4) \del_6 & p_1q_2\\
 & 2 u_1^{(2)} \del_1 + (u_1 u_2) \del_2 + u_5 \del_3 + 2 u_1^{(2)} u_2 \del_4 & p_3q_1\\
 & (2 u_5 + u_1 u_3) \del_2 + 2 u_1^{(2)} u_3 \del_4 + u_3 u_5 \del_6 & p_1p_3\\
 & (2 u_1 u_2 + u_4) \del_3 + u_1^{(2)} u_2 \del_5 + 2 u_1 u_2^{(2)} \del_6 & q_1q_2\\
 & u_1 u_4 \del_2 + 2 u_1 u_5 \del_3 + 2 u_1^{(2)} u_4 \del_4 + u_1^{(2)} u_5 \del_5 + (u_1 u_2 u_5 + u_4 u_5) \del_6 & -q_2p_3\\
 & u_1 u_5 \del_2 + 2 u_1^{(2)} u_5 \del_4 + u_5^{(2)} \del_6 & p_3^{(2)}\\
 & u_1 u_4 \del_3 + 2 u_1^{(2)} u_4 \del_5 + 2 (u_1 u_2 u_4 + u_4^{(2)}) \del_6 & -q_2^{(2)}\\
 \hline
 \end{array}
\ee

Let us describe $\fmy(7)_0$. Its center is $\fc=\Kee t$. Further, $\fmy(7)_0$ contains two commuting with each other subalgebras isomorphic to $\fsl(2)$:
\[
\begin{array}{l}
S_1=\Span(p_1^{(2)}-p_2p_3,q_1^{(2)}-q_2q_3, p_1q_1+p_2q_2+p_3q_3), \\
S_2=\Span(p_1q_3+q_1p_2,q_1p_3+p_1q_2, p_2q_2-p_3q_3).
\end{array}
\]

Besides, $\fmy(7)_0$ contains two $S_1$-invariant subspaces on each of which $S_1$ acts via the adjoint representation:
\[
\begin{array}{l}
V_1=\Span(q_1 p_3-p_1q_2, \ \ p_1p_3, \ \ q_1q_2), \\
 V_2=\Span(p_3^{(2)}, \ \ q_2p_3,\ \ q_2^{(2)}).
 \end{array}
\]
Moreover, 
\[
\begin{array}{l}
[V_2,V_2]=0, \ \ [V_1,V_1]=V_2: \\ 
\{q_1 p_3-p_1q_2, \ p_1p_3\}=p_3^{(2)}, \; \{q_1 p_3-p_1q_2, \ q_1q_2\}=-q_2^{(2)}, \; \{p_1p_3,\ q_1q_2\}=q_2p_3.
\end{array}
\]

Therefore, $V_1\oplus V_2$ is a~subalgebra isomorphic to $\fmy(6)_-$, i.e., to a~negative part of $\fo(7)$ in the grading of $\fmy(6)$. 
The action of $S_2$ on $V_1\oplus V_2\oplus S_1$ is as follows:
\[
\begin{array}{l}
\ad_{p_1q_3+q_1p_2}:\begin{array}{lclclcl}
 q_2^{(2)} & \mapsto & q_1q_2 & \mapsto & -q_1^{(2)}+q_2q_3 & \mapsto & 0\\
 q_2p_3 & \mapsto & q_1p_3-p_1q_2 & \mapsto & p_1q_1+p_2q_2+p_3q_3 & \mapsto & 0\\
 p_3^{(2)} & \mapsto & p_1p_3 & \mapsto & p_1^{(2)}-p_2p_3 & \mapsto & 0\\
 \end{array}\\
\hline
\ad_{p_1q_2+q_1p_3}:\begin{array}{lclcl}
 q_1q_2 & \mapsto & -q_2^{(2)} & \mapsto & 0\\
 q_1p_3-p_1q_2 & \mapsto & -q_2p_3 & \mapsto & 0\\
 p_1p_3 & \mapsto & -p_3^{(2)} & \mapsto & 0\\
 \end{array}
\end{array}
\]
Introduce a~dummy indeterminate~$x$.
The commutation rules in~$\fmy(7)_0$ make it possible to identify~$S_1$ with $S_1 \otimes 1$,
$V_1$ with $S_1 \otimes x$, and $V_2$ with $S_1 \otimes x^{(2)}$.
Then, the action of $p_1q_3 +q_1p_2$ turns into that of $\del_x$, and the action of $p_1q_2 +q_1p_3$ into that of $x^{(2)}\del_x$.
Observe that $S_2\simeq\fvect(1;1)$. 
Since the divergence on Lie (super)algebras of series $\fk$ is given by $\partial_t$, it follows from~\eqref{MY7geq} that 
$\fsmy(7; \un)_-=\fmy(7; \un)_-$ and $\fsmy(7; \un)_0\simeq\fmy(7)_0/\fc=\fsl(2) \otimes \cO(1; 1) \inplus \fvect(1;1)$.

\paragraph{$\fmy(7)_{-1}$ as $\fmy(7)_0$-module}
The central element~$t$ acts as the multiplication by $-1$.
The direct calculations show that if we introduce a~space $W$ by identifying $\Span(p_2,q_3)$ with $W\otimes 1$,
then $\Span(p_1,q_1)$ becomes identified with $W\otimes x$, and $\Span(p_3,q_2)$ with $W\otimes x^{(2)}$, and hence 
$\fmy(7)_{-1}\simeq W\otimes \cO(1;\One)$, whereas the action of $S_1\otimes \cO(1;\One )$ is $\id_{\fsl(2)}\otimes \cO(1;\One )$.
Having considered the elements of~$S_2$ as vector fields in $\fvect(1;\One )$, we see that~$D$ acts in $\fmy(7)_{-1}$ via the formula:
\[
 \ad_D: w\otimes f\mapsto w\otimes (L_D(f)-\Div(D)f \text{~~for any~~} w\in W, f\in \cO(1;\One ).
\]

For $p=3$, we have $-1=2=\nfrac12$, so the $\fmy(7)_{0}$-module $\fmy(7)_{-1}$ is $\id\otimes \Vol^{1/2}$ with the natural bracket with values in $\fmy(7)_{-2}$ induced by the antisymmetric non-degenerate bilinear form.

\paragraph{Remark}\label{strange1}
Consider a~non-standard grading of $\fsvect(3; \un)$ for $\un=(\infty, \infty, 1)$ setting $\deg x_3=0$ and $\deg x_1= \deg x_2=1$.
Normally, we can not set $\deg x_3=0$ because the height of $x_3$ is unbounded,
and we forbid gradings with infinite dimensions of homogeneous components; but if $N_3=1$, we can set $\deg x_3=0$.
Observe that $(\fmy(7)_{0})/\fc\simeq\fsvect(3; \un)_{0}$ and $\fmy(7)_{-1}\simeq\fsvect(3; \un)_{-1}$: 
\be\label{my7}
\fmy(7)_0=\fc(\fsl(2)\otimes \cO(1;\One) \inplus \fvect(1;\One)), \ \ \fmy(7)_{-1}=\id\otimes \Vol^{1/2} \text{~~as a~$\fg_0$-module.}
\ee
Having added a~center to $(\fmy(7)_{0})/\fc$ and taking the $(\fmy(7)_{0})$-invariant anti-symmetric non-degenerate form on $\fmy(7)_{-1}$ into account we acquire $\fmy(7)_{-2}$ as the space of values of this form. 




\sssbegin{Theorem}\label{TmMy}
The above two examples exhaust all W-gradings of $\fmy(6)$.
The growth vectors $v_\fg$ of non-integrable distributions preserved by the W-regradings of the middle
Skryabin algebra $\fmy(6)$ and its divergence free subalgebra $\fs\fmy(6)$ are as follows
\begin{equation}\label{GrVdy1}
\renewcommand{\arraystretch}{1.4}
\begin{tabular}{|c|c|}
\hline $\fg$&$v_\fg$\\
\hline
$\fmy(6)$ and $\fs\fmy(6)$&$(3,6)$\\
$\fmy(7)$ and $\fs\fmy(7)$&$(6,7)$C\\
\hline
\end{tabular}
\end{equation}
\end{Theorem}

\begin{proof} First of all, the above-described correspondence between the indeterminates $u_i$ and $t, p_j,q_j$ allows us to express each regrading not only in terms of the degrees $d_i:=\deg u_i$, but also in terms of the degrees of the contact indeterminates. In particular, by setting
\[
\deg p_2=\deg q_3=0,\ \ \deg p_1=\deg q_1=1, \ \ \deg p_3=\deg q_2=\deg t=2,\ \ 
\]
we get a~back the grading of $\fmy(6)$, where
\[
\fg_{-2}=\Span(1,\ \ p_2,\ \ q_3), \quad \fg_{-1}=\Span(p_1,\ \ q_1, \ \ p_1q_3+q_1p_2).
\]

The conditions on regradings in contact indeterminates are as follows. Let 
\[
\deg t=k, \; \deg p_i=m_i, \; \deg q_i=n_i, \; i=1,2,3.
\]
Since the polynomials in the 2nd, 3rd and 4th lines in table~\eqref{MY7geq} must remain in the degree-0
component (they act diagonally)
whereas all other monomials-summands in $\fmy(7)_0$ are linearly independent,
it follows that all polynomials in table~\eqref{MY7geq} must remain homogeneous in any grading.
This implies the following constraints on $k, m_i,n_i$:
\begin{equation}\label{kmn}
\begin{array}{l}
k=m_1+n_1=m_2+n_2=m_3+n_3=N, \\ 
m_1+n_3=m_2+n_1,\ \ n_2+n_3=2n_1,\ \ m_2+m_3=2m_1.
\end{array}
\end{equation}
We have $\gr(f)=\deg(f)-N$.

The corank of the system~\eqref{kmn}, same as the corank of the system~\eqref{rk3} for~$d_i$, is equal to 3. 
The condition $d_1=0$ imposes one more condition: the element $p_1q_3+q_1p_2$ must remain in the 0th component.
This requirement yields one more equation:
\[
m_1+n_3=m_2+n_1=N,
\]
which, together with the system~\eqref{kmn}, imply:
\[
m_1=m_2=m_3 \text{~~and~~} n_1=n_2=n_3.
\]

Observe that the heights of~$p_3$ and~$q_2$ are unconstrained.
Therefore, $q_3$ and~$p_2$ which act on functions independent of~$t$ as~$\del_{p_3}$ and~$-\del_{q_2}$, respectively,
must remain in the negative part of any regrading.
But then so must $p_i,q_i$ for all $i$, and the generating function $1$.
By Lemmas~\ref{L4.1} and~\ref{L4.1.1} this implies that such a~regrading is not a~Weisfeiler one,
unless it is equivalent to the grading of~$\fmy(7)$ with $m_i=n_i=1$ for all $i$ and $k=2$. 

Thanks to Lemmas~\ref{L4.1} and~\ref{L4.1.1} we have to equate at least one $d_i$ to 0;
it could be only some of $i=1,2,3$ since the heights of other indeterminates are unconstrained;
moreover, we can not equate to 0 more than one degree, since the commutator of any two of the first three corresponding elements of $\fg_{-1}$
yields an~element corresponding to an~indeterminate with unconstrained height;
lastly, the~$u_i$ for $i=1,2,3$ are on equal footing, as well as the $u_i$ for $i=4,5,6$, so it suffices to
set $d_1=0$.
Finally, incompressibility of the grading and fixing which direction of the grading is positive imply that
there exists only one, up to an isomorphism, $W$-grading of $\fmy(6)$ for which $d_1=0$: it is $\fmy(7)$.
\end{proof}

\paragraph{The little Skryabin algebras}\label{SkLittl}
In~\cite{GL}, the divergence-free subalgebras of~$\fmy(6; \un)$ and its filtered deformations are called \textit{little Skraybin algebras}. 
There are known 5 little Skryabin algebras~\cite{Sk}.
We consider only the $\Zee$-graded one, $\fsmy(6; \un)$, the special (divergence-free) subalgebra of $\fmy(6;\un)$,
originally denoted $\text{Y}_{(1)}(\un)$;
the other little Skryabin algebras are filtered deforms of the graded one, see~\cite{LaY, LaD, LaZ}. 
The algebra $\fsmy(6; \un)$ is the Cartan prolong whose
non-positive part is given by divergence-free part of~\eqref{MY6}.
Hence, the non-integrable distributions preserved by $\fsmy(6; \un)$ and $\fmy(6; \un)$ coincide. 

\subsection{The big Skryabin algebra $\fg=\fby(7)$}\label{fby}
Then, $\fg_0=\fgl(3)$ and $\fg_-$ can be realized as follows
(indicated are basis vectors and their non-zero weights with respect to the $h_i$), see~\cite[Section~5.5]{GL}:
\tiny
\begin{equation}\label{realby}
\renewcommand{\arraystretch}{1.4}
\begin{tabular}{|l|l|}
\hline $\fg_{0}\simeq\fgl(V)$&$(-1,0,1)$\quad $u_3 \del_1 + u_3^{(2)}
\del_5 + \left(u_2 u_3-u_4\right) \del_6 -u_2
u_3^{(2)} \del_7$\\
&$(-1,1,0)$\quad $-u_2 \del_1 + u_4 \del_5$\\
&$(0,-1,1)$\quad $-u_3 \del_2 + u_3^{(2)} \del_4 + \left(u_1 u_3 +u_5
\right)\del_6 -u_1
u_3^{(2)} \del_7$;\\
&$h_1=u_1 \del_1 + u_5 \del_5 + u_6 \del_6 + u_7 \del_7,$\\
&$h_2=u_2 \del_2 + u_4 \del_4 + u_6 \del_6 + u_7 \del_7,$\\
&$h_3=u_3 \del_3 + u_4 \del_4 + u_5 \del_5 + u_7 \del_7$;\\
&$(0,1,-1)$\quad $-u_2 \del_3 + u_2^{(2)} \del_4 + \left(u_1 u_2 + u_6
\right)\del_5 -u_1
u_2^{(2)} \del_7$,\\
&$(1,-1,0)$\quad $-u_1 \del_2 + u_5 \del_4$,\\
&$(1,0,-1)$\quad $-u_1 \del_3 + \left(-u_1 u_2 + u_6
\right)\del_4-u_1^{(2)} \del_5 +
u_1^{(2)} u_2 \del_7$\\
\hline $\fg_{-1}\simeq V$&$X_1:=\del_1+ u_2\del_6 + u_3 \del_5 - u_4
\del_7,\quad
X_2:=\del_2 - u_3 \del_4 - u_1 \del_6-u_5 \del_7,\ \
X_3:=\del_3 - u_6 \del_7$\\
\hline $\fg_{-2}\simeq E^2(V)$&$X_4:=\del_4 + u_1 \del_7,\quad
X_5:= \del_5 + u_2 \del_7,\quad
X_6:= \del_6 + u_3 \del_7$\\
\hline
$\fg_{-3}\simeq E^3(V)$&$X_7:=\del_7$\\
\hline
\end{tabular}
\end{equation}
\normalsize 

Clearly, $\fs\fby(7)_-=\fby(7)_-$, and $\fs\fby(7)_0=\fsl(3)$.

The vector fields~$Y_i$ commuting with $\fg_-$ are 
\tiny
\[
\renewcommand{\arraystretch}{1.4}
\begin{tabular}{|l|l|}
 \hline
 Weight & $Y_i$ \\
 \hline 
 $(2,0,0)$ & $Y_1 = \del_1 + 2u_2\del_6 + u_2u_3\del_7 + u_4\del_7$,\\
 $(0,2,0)$ & $Y_2 = \del_2 + u_1\del_6 + 2u_1u_3\del_7 + u_5\del_7$,\\
 $(0,0,2)$ & $Y_3 = \del_3 + u_1\del_5 + 2u_2\del_4 + u_6\del_7$,\\
 $(0,2,2)$ & $Y_4 = \del_4 + 2 u_1\del_7$, \\
 $(2,0,2)$ & $Y_5 = \del_5+2 u_2\del_7$,\\
 $(2,2,0)$ & $Y_6 = \del_6+2 u_3\del_7$,\\
 $(2,2,2)$ & $Y_7 = \del_7$ \\
 \hline
\end{tabular}
\]
\normalsize

For the dual basis of the 
module of 1-forms we take 
\begin{gather*}
\Theta^1 = 2 du_1,\qquad
\Theta^2 = 2 du_2,\qquad
\Theta^3 = 2 du_3,\qquad
\Theta^4 = 2 du_4+2 u_2du_3 ,\\
\Theta^5 = 2 du_5+u_1du_3 ,\qquad
\Theta^6 = 2 du_6+2 u_2du_1 +u_1du_2 ,\\
\Theta^7 = 2 du_7+u_4du_1 +u_5du_2 +u_6du_3 +2 u_1du_4 +2 u_2du_5 +2 u_3du_6 .
\end{gather*}

We can equate to 0 the heights of only constrained indeterminates listed in Remark~\ref{strange2}.
Looking at $\Theta^4$, $\Theta^5$ and $\Theta^6$ we see that no two of the first three indeterminates can be
of degree 0 simultaneously.
The homogeneity of the~$\Theta^i$ implies that
\be
\label{rk37}
\begin{gathered}
d_4 =d_2+d_3,\qquad
d_5 =d_1+d_3,\qquad
d_6 =d_1+d_2,\\
d_7 =d_1+d_4=d_2+d_5=d_3+d_6=d_1+d_4
\end{gathered}
\ee
or equivalently
\be\label{rk3sol}
d_4 = d_2+d_3,\quad d_5 = d_1+d_3,\quad d_6 = d_1+d_2, \quad d_7 = d_1+d_2+d_3.
\ee

\sssec{The W-regrading~$\fby(8)$ of~$\fby(7)$} Let $d_1=0$. Then, $\Theta^5$, $\Theta^6$ and $\Theta^7$ imply, respectively $d_3=d_5$, $d_2=d_6$, and
$d_4=d_7=d_2+d_5$. Let $d_3=d_5=d_2=d_6:=1$, then $d_4=d_7=2$. We denote this regrading $\fby(8;\un)$. For a~basis of the negative part we take
\footnotesize
\[
\renewcommand{\arraystretch}{1.4}
\begin{tabular}{|l|l|}
 \hline
 $\fby(8)_{-2}$
 & $c_1:=\del_7$, \ \ $c_2:=\del_4 + u_1\del_7$ \\
 \hline
 $\fby(8)_{-1}$
 & 
 $b_2:=\del_2 - u_3\del_4 - u_1\del_6 - u_5\del_7$,\ \ 
 $b_3:=\del_3 - u_6\del_7$,\ \ $b_5:=\del_5 + u_2\del_7$,\ \ 
 $b_6:=\del_6 + u_3\del_7$,\\ 
 &
 $b_1:=u_1\del_2 - u_5\del_4$,\ \ 
 $b_4:=u_1\del_3 - (u_6 - u_1 u_2)\del_4 - u_1^{(2)}\del_5 + u_1^{(2)} u_2\del_7$\\
 \hline
\end{tabular}
\]
\normalsize

For a~basis of $\fby(8)_{0}$ we take
\footnotesize
\begin{align*}
 a_1 :={}& u_2\del_6 + u_3\del_5 + 2 u_4\del_7 + \del_1,\\
 a_2 :={}& 2 u_3\del_2 + u_1 u_3\del_6 + u_3^{(2)}\del_4 + 2 u_1 u_3^{(2)}\del_7 + u_5\del_6,\\
 a_3 :={}& u_1\del_1 + u_5\del_5 + u_6\del_6 + u_7\del_7,\\
 a_4 :={}& u_2\del_2 + 2 u_3\del_3 + 2 u_5\del_5 + u_6\del_6,\\
 a_5 :={}& 2 u_2\del_3 + u_1 u_2\del_5 + u_2^{(2)}\del_4 + 2 u_1 u_2^{(2)}\del_7 + u_6\del_5,\\
 a_6 :={}& u_1^{(2)}\del_1 + 2u_1 u_2\del_2 + 2 u_1^{(2)} u_2\del_6 + 2u_1 u_3\del_3 + 2 u_1^{(2)} u_3\del_5 
 + 2 u_1 u_2 u_3\del_4 + 2 u_1^{(2)} u_2 u_3\del_7 + 2u_1 u_4\del_4\\
 {}&{}+ 2 u_1^{(2)} u_4\del_7 + 2 u_2 u_5\del_4 + u_1 u_2 u_5\del_7 + 2u_3 u_6\del_4 + u_1 u_3 u_6\del_7 + u_7\del_4 + u_1u_7\del_7,\\
 a_7 :={}& u_1^{(2)}u_2\del_6 + u_1u_3\del_3 + u_1^{(2)}u_3\del_5 + u_1u_2u_3\del_4 + u_1^{(2)}u_2u_3\del_7 + u_5\del_3 + u_6\del_2 + 2u_1u_6\del_6 + 2u_3u_6\del_4 + u_1u_3u_6\del_7,\\
 a_8 :={}& u_1 u_3\del_2 + u_5\del_2 + 2u_1 u_5\del_6 + 2 u_3 u_5\del_4 + u_1 u_3 u_5\del_7,\\
 a_9 :={}& 2 u_1^{(2)} u_2\del_5 + 2 u_1 u_2^{(2)}\del_4 + 2 u_1^{(2)} u_2^{(2)}\del_7 + u_6\del_3,\\
 a_{10} :={}& u_1^{(2)} u_2\del_2 + 2 u_1 u_5\del_3 + 2 u_1^{(2)} u_5\del_5 + 2 u_1 u_2 u_5\del_4 + 2 u_1^{(2)} u_2 u_5\del_7 + 2 u_1 u_6\del_2 + u_5 u_6\del_4 + 2 u_1 u_5 u_6\del_7,\\
 a_{11} :={}& u_1 u_5\del_2 + 2 u_5^{(2)}\del_4 + u_1 u_5^{(2)}\del_7,\\
 a_{12} :={}& u_1^{(2)} u_2\del_3 + 2 u_1 u_6\del_3 + 2 u_1^{(2)} u_6\del_5 + 2 u_1 u_2 u_6\del_4 + 2 u_1^{(2)} u_2 u_6\del_7 + u_6^{(2)}\del_4 + 2 u_1 u_6^{(2)}\del_7,\\
 a_{13} :={}& u_2\del_2 + u_3\del_3 + 2 u_4\del_4 + u_5\del_5 + u_6\del_6 + 2 u_7\del_7.
\end{align*}
\normalsize

The multiplication table in $\fby(8)_0$ is as follows:
\tiny
\[
\begin{array}{c|ccccccccccccc}
 {} & a_1 & a_2 & a_3 & a_4 & a_5 & a_6 & a_7 & a_8 & a_9 & a_{10} & a_{11} & a_{12} & a_{13} \\
 \hline
 a_1 & 0 & 0 & a_1 & 0 & 0 & a_3+2 a_{13} & a_4 & a_2 & 2 a_5 & 2 a_7 & a_8 & 2 a_9 & 0 \\
 a_2 & 0 & 0 & 0 & 2 a_2 & 2 a_4 & 0 & 2 a_8 & 0 & a_7 & a_{11} & 0 & a_{10} & 0 \\
 a_3 & 2 a_1 & 0 & 0 & 0 & 0 & a_6 & a_7 & a_8 & a_9 & 2 a_{10} & 2 a_{11} & 2 a_{12} & 0 \\
 a_4 & 0 & a_2 & 0 & 0 & 2 a_5 & 0 & 0 & a_8 & 2 a_9 & 0 & a_{11} & 2 a_{12} & 0 \\
 a_5 & 0 & a_4 & 0 & a_5 & 0 & 0 & 2 a_9 & a_7 & 0 & 2 a_{12} & 2 a_{10} & 0 & 0 \\
 a_6 & 2 a_3+a_{13} & 0 & 2 a_6 & 0 & 0 & 0 & 2 a_{10} & a_{11} & 2 a_{12} & 0 & 0 & 0 & 0 \\
 a_7 & 2 a_4 & a_8 & 2 a_7 & 0 & a_9 & a_{10} & 0 & 2 a_{11} & 2 a_{12} & 0 & 0 & 0 & 0 \\
 a_8 & 2 a_2 & 0 & 2 a_8 & 2 a_8 & 2 a_7 & 2 a_{11} & a_{11} & 0 & a_{10} & 0 & 0 & 0 & 0 \\
 a_9 & a_5 & 2 a_7 & 2 a_9 & a_9 & 0 & a_{12} & a_{12} & 2 a_{10} & 0 & 0 & 0 & 0 & 0 \\
 a_{10} & a_7 & 2 a_{11} & a_{10} & 0 & a_{12} & 0 & 0 & 0 & 0 & 0 & 0 & 0 & 0 \\
 a_{11} & 2 a_8 & 0 & a_{11} & 2 a_{11} & a_{10} & 0 & 0 & 0 & 0 & 0 & 0 & 0 & 0 \\
 a_{12} & a_9 & 2 a_{10} & a_{12} & a_{12} & 0 & 0 & 0 & 0 & 0 & 0 & 0 & 0 & 0 \\
 a_{13} & 0 & 0 & 0 & 0 & 0 & 0 & 0 & 0 & 0 & 0 & 0 & 0 & 0 \\
\end{array}
\]
\normalsize

The algebra $\fby(8)_0$ has a~$1$-dimensional center $\fc=\Span(a_{13})$ and a~$9$-dimensional ideal 
\[
{I=\Span(a_2, a_4, a_5, a_7, a_8, a_9, a_{10}, a_{11}, a_{12})}.
\]
Clearly, $S_1:=\Span(a_2=X_-, \ a_4=H, \ a_5=X_+)\simeq\fsl(2)$; set 
\[
\text{$V_1=\Span(a_9, \ a_7, \ a_8)$
and $V_2=\Span(a_{12}, \ a_{10}, \ a_{11})$.}
\]
Clearly, $V_1$ and~$V_2$ are isomorphic adjoint $\fsl(2)$-modules. Moreover, 
\[
\begin{array}{l}
[V_2,V_2]=0, \ \ [V_1,V_1]=V_2. \\ 
\end{array}
\]
The multiplication table in $\fby(8)_0/(I\oplus\fc)\simeq \fsl(2)$ is as follows:
\[
\begin{array}{l}
[a_1,\ a_6]=a_3, \ \ [a_3, a_1]=-a_1, \ \ [a_3, a_6]=a_6. \\ 
\end{array}
\]
Introduce a~dummy indeterminate $x$.
The commutation rules in $\fby(8)_0$ make it possible to identify $S_1$ with $S_1 \otimes 1$, $V_1$ with $S_1 \otimes x$, and $V_2$ with $S_1 \otimes x^{(2)}$. Then, the action of $a_1$ turns into that of $\del_x$, and the action of $a_6$ into that of $x^{(2)}\del_x$. Since $\fsl(2)\simeq\fvect(1;1)$, 
we see that $\fby(8)_0$ is isomorphic to $\fc(\fsl(2) \otimes \cO(1; 1) \inplus \fvect(1;1))$, and since the center is a direct summand, the negative part of divergence-free subalgebra is the same while $\fs\fby(8)_0\simeq \fsl(2) \otimes \cO(1; 1) \inplus \fvect(1;1)$. 

We see that $\fby(8)_{-1}\simeq W\otimes \cO(1;\One )$, where $W=\id_{\fsl(2)}$ and $\fby(8)_{-2}\simeq \cO(1;\One )/\Kee 1$,
compare with expression~\eqref{my7}.
Let $\pi: \cO(1;\One )\to\cO(1;\One )/\Kee 1$.
Then, for any $v,w\in W$, and $f,g\in\cO(1;\One )$, the bracket $[\fby(8)_{-1}, \fby(8)_{-1}]\to \fby(8)_{-2}$ is given by the formula
\[
 [v\otimes f, w\otimes g]=(w,v)\pi(fg),
\]
where $(-,-)$ is the non-degenerate $\fsl(2)$-invariant anti-symmetric form on $W$.

\sssec{Remark}\label{strange2}
Compare with Remark~\ref{strange1}.
Consider a~non-standard grading of $\fh(4; \un)$ for $\un=(1,\ \infty, \infty, \infty)$ setting $\deg p_1=0$
the degrees of other indeterminates being equal to~1.
Normally, we can not set $\deg p_1=0$ because the height of $p_1$ is unbounded, and we forbid gradings with
infinite dimensions of homogeneous components; but if $N_1=1$, we can set $\deg p_1=0$.
Denote thus graded algebra by $\fh(4; \un; 1)$.
Observe that $\fby(8)_{0}/\fc\simeq\fh(4;\un; 1)_{0}$, where $\fh(4;\un; 1)$ is considered in the non-standard grading 
\be\label{byNst}
\deg p_1 =0, \ \ \deg q_1 =2, \ \ \deg p_2 = \deg q_2 =1, 
\ee
and then $\fby(8)_{-1}\simeq \fh(4; \un; 1)_{-1}$.
Having added a~center to $(\fby(8)_{0})/\fc$ and taking the Cartan prolongation we acquire $\fby(8)$ with $\dim\fby(8)_{-2}=2$.

Since $\fby(8)_-$ coincides with $\fh(4;1)_-$ while the bases of the 0th components of these algebras differ
by a~grading operator, it follows that $\fh(4;1)$ can be embedded to $\fby(8)$ as a~subalgebra.
Let $\fh(4;1)$ be realized in terms of generating functions in $p_1,p_2,q_1,q_2$ with the grading~\eqref{byNst}. 

Now, the form of $\fh(4;1)_{\leq 0}$ shows that the height of $p_1$ is equal to 1, and there are no
constraints on the heights of other indeterminates.
Since, on one hand, the element $q_1$ acts on the functions generating $\fh(4;1)$ as $-\del_{p_1}$, whereas,
on the other hand, its action on the negative part coincides with the action of $a_1\in \fby(8)$, we conclude
that in $\fby(8)$, and hence in the initial algebra $\fby$, the height of $u_1$ is constrained to be equal
to~1.
Besides, the elements $p_2,\ q_2\in\fh(4;1)$ correspond to $\del_5+u_2\del_7,\ \del_6+u_3\del_7\in \fby(8)$
implying that there are no constraints on the heights of $u_5$, $u_6$.
Finally, the element $p_1$ acts in $\fh(4;1)$ as $\del_{q_1}$ and corresponds to $\del_7\in\fby(8)$.
Therefore, there are no constraints on the height of $u_7$, either. 

\sssbegin{Theorem}\label{TmBy}
The above two pairs of examples exhaust all W-gradings of the big Skryabin algebra and its divergence-free
subalgebra.
The growth vectors $v_\fg$ of non-integrable distributions preserved by the W-regradings of these algebras are as follows
\begin{equation}\label{GrVdy2}
\renewcommand{\arraystretch}{1.4}
\begin{tabular}{|c|c|}
\hline $\fg$&$v_\fg$\\
\hline
$\fby(7)$ and $\fs\fby(7)$&$(3, 6,7)$\\
\hline
$\fby(8)$ and $\fs\fby(8)$&$(6,8)$\\
\hline
\end{tabular}
\end{equation}
\end{Theorem}

\begin{proof}
Turning to the initial depth-3 grading of $\fby$ we see that the indeterminates $u_1, u_2, u_3$ are on equal
footing, and similarly $u_4, u_5, u_6$ are also on equal footing.
Therefore, the heights of $u_2, u_3$ are constrained to be equal to 1, and the heights of $u_4, u_5, u_6$ are not constrained.

For notation $X_i$, see table~\eqref{realby}. 
Let $\fg$ be a~$W$-regrading of $\fby$ for which $d_1=0$.
Then, $X_1\in \fg_0$, and since $[X_1,X_3]=-X_5$, the elements $X_3$ and $X_5$ lie in the same component of
negative degree because the height of the indeterminate $u_5$ corresponding to $X_5$ is unconstrained.
Likewise, the elements $X_2$ and $X_6$ lie in the same component.
Finally, since $X_4=[X_2,X_3]$ whereas $X_7=-[X_2,X_5]$, it follows that $X_4, X_7\in\fg_-$.
Thus, $\fg_-\cap \fby(8)_0=0$ and Lemmas~\ref{L4.1} and~\ref{L4.1.1} imply that $\fg=\fby(8)$ as $\Zee$-graded
Lie algebras.
\end{proof}

\section{$p=3$, the three types of superized Melikyan algebras, see \cite{BL, BGL}}\label{SBrj}

\ssec{$\fg:=\fMe(3;\un|3)$}\label{Me33}
Let $\fMe(3;\un|3)$ denote a~super analog of~$\fme(\un)$.
Observe that $\fk(1;1|1)$ is not simple if $p=3$: namely, $\fk(1;1|1)'\simeq\fosp(1|2)$.
Its non-positive components are 
\be\label{vague}
\begin{array}{l}
\fg_0=\fk(1;1|1)'\oplus \Kee z,\\
\fg_{-1}=\Pi(\fk(1;1|1)')\stackrel{\text{as spaces}}{\simeq} \cO(1;1|1)/\Kee 1,\\
\fg_{-2}=\Kee 1.\\
\end{array}
\ee
We assume that $\fk(1;1|1)$ preserves the distribution singled out by the contact form $\alpha=dt-\theta
d\theta$; and consider the following $\fg_0$-invariant anti-symmetric non-degenerate bilinear form on~$\fg_{-1}$:
\[
 (f, g)=\int_\cB fdg \pmod\alpha= \int_\cB(f_0g_0+(f_1g_1-f_0g_0')\theta)d\theta=\int (f_1g_1-f_0g_0')dt 
\] 
for any $f=f_0+f_1\theta$ and $g=g_0+g_1\theta$, where $f_0,f_1,g_0,g_1\in\cO(1;1)$ and
where the Berezin integrals $\int_\cB$ on $1|1$-dimensional supermanifolds and $\int $ on $1$-dimensional
manifolds are defined for $p>2$ as follows 
\[
 \begin{array}{ll}
 \int_\cB:fd\theta\mapsto b_{p-1}&\text{for $f=\sum (a_it^{(i)} + b_it^{(i)}\theta)\in\cO(1;1|1)$,}\\
 \int :gdt\mapsto b_{p-1}&\text{for $g=\sum b_it^{(i)}\in\cO(1;1)$}.
 \end{array}
\]
Therefore, for $p=3$, there is an anti-symmetric non-degenerate form on $\fg_{-1}=\cO(1;1|1)/\Kee 1$:
\[
 \begin{array}{l}
 (t, t^{(2)})=-\int t\cdot tdt=-\int 2t^{(2)}dt =1\qquad (t^{(2)}, t)=-\int t^{(2)}dt=-1;\\
 (\theta, t^{(2)}\theta)= (t^{(2)}\theta, \theta)=1,\ \ 
 (t\theta, t\theta)=\int 2t^{(2)}dt=-1.
 \end{array}
\]
We realize $\fk(1; 1|1)'$ by generating functions in $t, \theta$. Then, 
\[
 \fg_{-1}= \Span(t, t^{(2)}, \theta, t\theta, t^{(2)}\theta).
\]
Since the contact field generated by $f$ is (recall the expression \eqref{K_f})
\[
 K_f =(2-\theta\partial_\theta)(f)\partial_t+(-1)^{p(f)}\pderf{f}{\theta}\partial_\theta,
\]
it follows that
\[
 \begin{array}{l}
 \fg_0=\Span(K_f\mid f\in \Span(1,t,t^{(2)},\theta,t\theta)).\\
 \end{array}
\]
For completeness of the picture,
let us also consider $K_{t^{(2)}\theta} = t^{(2)}\theta\partial_t -t^{(2)}\partial_\theta$,
in order to show that this element does not belong not only to $(\fg_0)'$,
but also to $[\fg_{-1}, \fg_1]$ meaning that if we add this element to $\fg_0$,
the positive components of the Cartan prolongation will not change.

Let us embed $\fg$ into $\fk(3|3)$ realized by generating functions in $p, q, \xi, \eta, \zeta; z$.
The contact bracket in these indeterminates takes the form (cf. with eq.~\eqref{kb1})
\be\label{kb2}
\begin{array}{ll}
\{f,g\}=&
(2-E)(f)\frac{\der g}{\der z}-\frac{\der f}{\der z}(2-E)(g)-\frac{\der f}{\der p}\cdot\frac{\der g}{\der q}+\frac{\der f}{\der q}\cdot\frac{\der g}{\der p}\\
&+
(-1)^{p(f)}\left(\frac{\der f}{\der \xi}\cdot\frac{\der g}{\der \eta}+\frac{\der f}{\der \eta}\cdot\frac{\der g}{\der \xi}+\frac{\der f}{\der \zeta}\cdot\frac{\der g}{\der \zeta} \right),
\end{array}
\ee
where $E:=\sum y_i\del_{y_i}$ and the $y_i$ are all the indeterminates except $t$.
Accordingly, we get
\be\label{(1)}
\{p,q\} = -1, \ \ \{\xi,\eta\} = -1, \ \ \{\zeta,\zeta\} = -1, 
\ee
inducing the following isomorphism of $(-1)$st components:
\be\label{(2)} 
t\lra q, \ \ t^{(2)} \lra p, \ \ \theta\lra -\xi, \ \ t\theta\lra \zeta, \ \ t^{(2)}\theta \lra\eta.
\ee

To describe embedding of $\fg_0$ into $\fk(3|3)_0$ let us first collect in a~table how $\fk(3|3)_0$ acts on $\fk(3|3)_{-1}$:
\be\label{(T1)}
\begin{tabular}{|c|c|c|c|c|c|||c|c|c|c|c|c|}
\hline
 Operator& $p$ & $q$ & $\xi$ & $\eta$ & $\zeta$ & Operator &$p$ & $q$ & $\xi$ & $\eta$ & $\zeta$\\
 \hline \hline
 $\ad_{p^{(2)}}$ & 0 & $-p$ & 0 & 0 & 0 & $\ad_{p\xi}$ & 0 & $-\xi$ & 0 & $-p$ &0\\
 \hline
 $\ad_{pq}$ & $p$ & $-q$ & 0 &0 & 0 & $\ad_{p\eta}$ & 0& $-\eta$ & $-p$ & 0 & 0\\
 \hline
 $\ad_{q^{(2)}}$ & $q$& 0& 0& 0 & 0 & $\ad_{p\zeta}$ & 0 & $-\zeta$ & 0 &0 & $-p$\\
 \hline
 $\ad_{\xi\eta}$ & 0 &0& $-\xi$ & $\eta$& 0 & $\ad_{q\xi}$ & $\xi$& 0& 0& $-q$ & 0\\
 \hline
 $\ad_{\xi\zeta}$ & 0 &0& 0 & $\zeta$& $-\xi$ & $\ad_{q\eta}$ & $\eta$& 0& $-q$&0 & 0\\
 \hline
 $\ad_{\eta\zeta}$ & 0 &0& $\zeta$& 0 & $-\eta$ & $\ad_{q\zeta}$ & $\zeta$& 0&0 & 0& $-q$\\
 \hline
\end{tabular}
\ee

The next table describes the embedding of $\fk(1;1|1)$ into $\fk(3;\un|3)_0$ with unconstrained $\un$.
For this, we first describe how $\fk(1;1|1)$ acts on $\cO(1;1|1)/\Kee1$,
next this action we rewrite this action using isomorphism~\eqref{(2)} and the table~\eqref{(T1)}
allows us to recover the generating function~$\widehat f$ with the needed action of~$K_f$:
\be\label{(T2)}
\begin{tabular}{|l|l|l|l|}
\hline
$K_f\in$ & The action of $K_f$ in & The induced action $K_{\widehat f}$ in & $\widehat f\in $ \\
 $\fk(1;1|1)$ & $\Kee[t,\theta]$ & $\Kee[p,q,\xi,\eta,\zeta]$ &$\Kee[p,q,\xi,\eta,\zeta]$\\ 
\hline\hline
$K_1$ & $t\mapsto 0, \; t^{(2)}\mapsto -t,$ & $q\mapsto 0,\; p\mapsto -q$ & $-q^{(2)}-\xi\zeta$\\
 & $\theta \mapsto 0,\; t\theta \mapsto -\theta,\; t^{(2)}\theta \mapsto -t\theta $ & $\xi\mapsto 0,\; \zeta \mapsto \xi, \; \eta \mapsto -\zeta$ &\\
 \hline
$K_t$ & $t\mapsto -t, \; t^{(2)}\mapsto t^{(2)},$ & $q\mapsto -q,\; p\mapsto p$ & $pq-\xi\eta$\\
 & $\theta \mapsto \theta,\; t\theta \mapsto 0,\; t^{(2)}\theta \mapsto -t^{(2)}\theta $ & $\xi\mapsto \xi,\; \zeta \mapsto 0, \; \eta \mapsto -\eta$ &\\
 \hline
 $K_{t^{(2)}}$ & $t\mapsto -t^{(2)}, \; t^{(2)}\mapsto 0,$ & $q\mapsto -p,\; p\mapsto 0$ & $p^{(2)}-\eta\zeta$\\
 & $\theta \mapsto t\theta,\; t\theta \mapsto t^{(2)}\theta,\; t^{(2)}\theta \mapsto 0 $ & $\xi\mapsto -\zeta,\; \zeta \mapsto \eta, \; \eta \mapsto 0$ &\\
 \hline \hline
 $K_\theta$ & $t\mapsto \theta, \; t^{(2)}\mapsto t\theta,$ & $q\mapsto -\xi,\; p\mapsto \zeta$ & $p\xi+q\zeta$\\
 & $\theta \mapsto 0,\; t\theta \mapsto -t,\; t^{(2)}\theta \mapsto -t^{(2)}$ & $\xi\mapsto 0,\; \zeta \mapsto -q, \; \eta \mapsto -p$ &\\
 \hline
 $K_{t\theta}$ & $t\mapsto t\theta, \; t^{(2)}\mapsto -t^{(2)}\theta,$ & $q\mapsto \zeta,\; p\mapsto -\eta$ & $-p\zeta-q\eta$\\
 & $\theta \mapsto -t,\; t\theta \mapsto t^{(2)},\; t^{(2)}\theta \mapsto 0$ & $\xi\mapsto q,\; \zeta \mapsto p, \; \eta \mapsto 0$ &\\
 \hline
 $K_{t^{(2)}\theta}$ & $t\mapsto t^{(2)}\theta, \; t^{(2)}\mapsto 0,$ & $q\mapsto \eta,\; p\mapsto 0$ & $-p\eta$\\
 & $\theta \mapsto -t^{(2)},\; t\theta \mapsto 0,\; t^{(2)}\theta \mapsto 0$ & $\xi\mapsto p,\; \zeta \mapsto 0, \; \eta \mapsto 0$ &\\
 \hline
\end{tabular}
\ee

To find the elements of components of higher degrees we have to single out the component $\fg_0$ in $\fk(3;
\un|3)_0$ by means of equations expressed in terms of operators $\widetilde K_{\widehat f}$ commuting with all
operators $K_f \in\fg_-$ (the operators $\widetilde K_{\widehat f}$ are denoted by $Y$s in other examples):
\[
 K_p=p\der_z-\der_q, \;
 K_q=q\der_z+\der p, \;
 K_\xi=\xi\der_z-\der_\eta, \;
 K_\eta=\eta\der_z-\der_\xi, \;
 K_\zeta=\zeta\der_z-\der_\zeta, \;
 K_1= 2\der_z.
\]

Then,
\[
 \begin{array}{l}
 \widetilde K_p=p\der_z+\der_q, \ \ \widetilde K_q=q\der_z-\der p, \\ 
 \widetilde K_\xi=(\xi\der_z+\der_\eta)\Par, \ \ \widetilde K_\eta=(\eta\der_z+\der_\xi)\Par, \ \ \widetilde K_\zeta=(\zeta\der_z+\der_\zeta)\Par, \ \ K_1= 2\der_z,
 \end{array}
\]
where $\Par(f):=(-1)^{p(f)}f$.

The equations that single out the functions generating $\fg_0$ are
\[
 \begin{array}{l}
 (\widetilde K_p)^2(f)=\widetilde K_\zeta\widetilde K_\eta(f); \quad \widetilde K_q\widetilde K_p(f)=\widetilde K_\xi\widetilde K_\eta(f); \quad (\widetilde K_q)^2(f)=-\widetilde K_\zeta\widetilde K_\xi(f);\\
 \widetilde K_q\widetilde K_\xi(f)=-\widetilde K_p\widetilde K_\zeta(f); \quad \widetilde K_q\widetilde K_\zeta(f)=-\widetilde K_p\widetilde K_\xi(f); \quad \widetilde K_p\widetilde K_\eta(f)=0.
 \end{array}
\]

Observe that the operators $\widetilde K_f$ act on the functions independent of $z$ as partial derivatives
(up to a~sign), so the functions independent of $z$ are singled out by a~simpler system of equations:
\be\label{simpSys}
\nfrac{\der^2 f}{(\der q)^2}=\nfrac{\der^2 f}{\der \zeta\der \xi}; \quad \nfrac{\der^2 f}{\der p\der q}=-\nfrac{\der^2 f}{\der \eta\der \xi}; \quad \nfrac{\der^2 f}{(\der p)^2}=-\nfrac{\der^2 f}{\der \zeta\der \eta}; \quad \nfrac{\der^2 f}{\der p\der \xi}=\nfrac{\der^2 f}{\der q\der \zeta}; \quad \nfrac{\der^2 f}{\der p\der \zeta}=\nfrac{\der^2 f}{\der q\der \eta}; \quad \nfrac{\der^2 f}{\der q\der \xi}=0.
\ee

The component of least degree where we can observe the constrained coordinates of the shearing vector $\un$ is~$\fg_1$.
Solving the system~\eqref{simpSys} for 3rd degree functions we get a~unique, up to a~scalar multiple, solution:
\[
 G=q^{(2)}\eta+\xi\zeta\eta+pq\zeta+p^{(2)}\xi.
\]
In particular, this implies that the heights of~$p$ and~$q$ are automatically equal to 1.

As $\fg_0$-module, $\fg_1/\Kee\cdot G\simeq (\fg_{-1})^*$, and its basis is
\begin{align*}
 F_\eta :={} & \eta z-p^{(2)}\zeta+pq\eta,
 & F_q = {} & qz-pq^{(2)}+p\xi\zeta-q\xi\eta,\\
 F_\xi :={} & \xi z+q^{(2)}\zeta-pq\xi,
 & F_p :={} & pz+p^{(2)}q+p\xi\eta+q\eta\zeta, \\
 F_\zeta :={} & \zeta z-q^{(2)}\eta+p^{(2)}\xi.\\
\end{align*}

For $k\geq 2$ we have
\small
\begin{equation}
\renewcommand{\arraystretch}{1.4}
\begin{tabular}{|l|l|} 
 \hline
 $\fg_{2k}$
 & $2 z^{(k+1)} + (\eta\xi + pq) z^{(k)} + (p^{(2)}q^{(2)} + 2\eta\xi pq + \zeta\eta q^{(2)}) z^{(k-1)} + 2\eta\xi p^{(2)}q^{(2)} z^{(k-2)}$, \\
 & $(2\zeta\xi + q^{(2)}) z^{(k)} + (\zeta\xi pq + \eta\xi q^{(2)}) z^{(k-1)} + \zeta\xi p^{(2)} q^{(2)} z^{(k-2)}$, \\
 & $2z^{(k+1)} + (2pq + 2\eta\xi) z^{(k)} + (p^{(2)}q^{(2)}+ 2\zeta\xi p^{(2)}+2\eta\xi pq) z^{(k-1)} + \eta\xi p^{(2)}q^{(2)}z^{(k-2)}$, \\
 & $(\xi p + \zeta q) z^{(k)} + (\xi p^{(2)}q + 2\zeta pq^{(2)} + \zeta\eta \xi q) z^{(k-1)} + 2\zeta\eta\xi pq^{(2)} z^{(k-2)}$, \\
 & $(\zeta p + \eta q) z^{(k)} + (\zeta p^{(2)}q + 2\eta p q^{(2)} + 2\zeta\eta\xi p) z^{(k-1)} + 2\zeta\eta\xi p^{(2)}q z^{(k-2)}$,\\
 & $(\zeta\eta + p^{(2)}) z^{(k)} + (2\eta\xi p^{(2)} +\zeta \eta p q) z^{(k-1)} + 2\zeta\eta p^{(2)}q^{(2)} z^{(k-2)}$\\
 \hline
 $\fg_{2k+1}$
 & $\xi z^{(k+1)} + (2\xi pq + \zeta q^{(2)})z^{(k)} + (\zeta\eta\xi q^{(2)} + 2\xi p^{(2)}q^{(2)})z^{(k-1)}$, \\
 & $\zeta z^{(k+1)} + (\eta q^{(2)} + 2\zeta\eta\xi + 2\zeta pq)z^{(k)} + (2\zeta p^{(2)}q^{(2)} + \zeta\eta\xi pq)z^{(k-1)} + \zeta\eta\xi p^{(2)}q^{(2)}z^{(k-2)}$, \\
 & $q z^{(k+1)} + (2pq^{(2)} + 2\zeta\xi p + \eta\xi q) z^{(k)} + (2\zeta\xi p^{(2)}q + 2\eta\xi p q^{(2)}) z^{(k-1)}$, \\
 & $2\zeta z^{(k+1)} + (2\zeta\eta\xi + \xi p^{(2)} + 2\zeta pq)z^{(k)} + (2\zeta\eta\xi pq+ \zeta p^{(2)}q^{(2)})z^{(k-1)} + \zeta\eta\xi p^{(2)}q^{(2)}z^{(k-2)}$, \\
 & $2\eta z^{(k+1)} + (2\eta pq + \zeta p^{(2)}) z^{(k)} + (2\zeta\eta\xi p^{(2)} + \eta p^{(2)} q^{(2)}) z^{(k-1)}$, \\
 & $p z^{(k+1)} + (2\eta\xi p + 2\zeta\eta q + p^{(2)}q)z^{(k)} + (2\eta \xi p^{(2)}q + \zeta\eta pq^{(2)})z^{(k-1)}$ \\
 \hline
\end{tabular}
\end{equation}
\normalsize

\sssec{W-gradings of $\fMe(3;\un|3)$}\label{sss8.1.1}
Let $\fh$ be a~Lie superalgebra obtained from~$\fg$ by a~regrading defined by degrees assigned to
indeterminates $p, q, \xi, \eta, \zeta$.
By Lemma~\ref{L4.1}, for the grading of $\fh$ to be a~W-grading, it is necessary that $\fh_0 \cap\fg_-\neq 0$.
The element $1\in\fg_{-2}$ must remain in~$\fh_-$.
Then, the commutation relations~\eqref{(1)} imply that $\zeta\in \fh_-$, and only one element of each pair
$(p, q)$ and $(\xi, \eta)$ may appear in $\fh_0$. 

Since all the summands entering the generating functions in the right column of table~\eqref{(T2)} are
linearly independent, the functions themselves must remain homogeneous under regrading.
This means that the degrees of monomials entering each of these generating functions must coincide.
This yields 5 equations on the degrees of the regrading:
\be\label{eqs}
\begin{array}{l}
2\deg q=\deg \xi+\deg \zeta, \;\deg p+\deg q=\deg \xi+\deg \eta, \; 2\deg p=\deg \eta+\deg \zeta, \\
\deg p+\deg \xi=\deg q+\deg \zeta, \; \deg p+\deg \zeta=\deg q+\deg \eta.
\end{array}
\ee
Observe that $pq - \xi\eta\in\fh_0$. Hence, if 
\be\label{N}
N = \deg p + \deg q = \deg \xi + \deg \eta, 
\ee
then $\gr f=\deg f - N$.
Therefore, if an element of the pair $(p,q)$ (resp. $(\xi,\eta)$) belongs to $\fh_0$,
and hence its degree is equal to $N$, the degree of the other element in the pair should be equal to 0. 

Let us show that only one of the elements $p, q, \xi,\eta$ can occur in $\fh_0$,
and hence the degree of only one of these elements can be equal to 0.
Indeed, if, e.g., $\deg p = \deg\xi = 0$, then $\deg q = \deg\eta = N$.
The 4th of equations~\eqref{eqs} then implies that $\deg\zeta= -N$, and the 5th equation yields that $-N = 2N$, hence $N = 0$.
This is a~contradiction. The remaining cases are similarly considered. 

Now, observe that the rank of the system~\eqref{eqs} is equal to 3,
and hence the system has 2 linearly independent solutions,
one of which is $(1, 1, 1, 1, 1)$.
Thus, the condition that the selected indeterminate is of degree 0 uniquely fixes the grading.

\paragraph{The W-regrading $\fMe(3;\un|4)$}\label{Ex1} We set
\[
\deg p=0, \; \deg q=2, \; \deg \xi=3, \; \deg \eta=-1, \; \deg \zeta=1, \; \deg z=2 
\]
implying that $\gr:=\deg-2$ is the W-grading in the regraded algebra. We get:
\[
\begin{tabular}{|c|c|c|c|c|}
\hline
$\fh_{-3}$ & $\fh_{-2}$ & $\fh_{-1}$ & $\fh_0$ & $\fh_{1}$ \\
\hline\hline
$\eta$ & $p$ & $\zeta$ & $q$ & $\xi$ \\
& $p^{2}-\eta\zeta$ & $p\zeta+q\eta$ &$pq-\xi\eta$, $z$ & $p\xi+q\zeta$\\
 & 1 & $F_\eta$ & $F_p$ & $F_\zeta$, $G$\\
 \hline
\end{tabular}
\]
Thus, we got a~compatible (with parity) $\Zee$-grading in which
\[
\begin{array}{l}
\fh_0 = \fsl(2)\oplus \Kee z, \ \ \fh_{-1} = \Pi(\ad_{\fsl(2)})\boxtimes \Kee[-1], \ \ \fh_{-2} = [\fh_{-1}, \fh_{-1}] = \ad_{\fsl(2)}\boxtimes \Kee[-2], \\ 
\fh_{-3} = [\fh_{-1}, \fh_{-2}] = \Pi(\Kee[-3]) = \Pi(\Kee[0]).
\end{array}
\]
We see that this $\Zee$-grading is a~W-grading. We denote thus regraded superalgebra $\fMe(3;\un|4)$. 

By setting
\[
 \deg p=2, \ \ \deg q=0, \ \ \deg \xi=-1, \ \ \deg \eta=3, \ \ \deg \zeta=1, \ \ \deg z=2 
\]
(implying that $\gr=\deg-2$ is the W-grading in the regraded algebra) we get an isomorphic $\Zee$-graded
algebra.
We denote thus regraded superalgebra $\fMe(3;\un|4)$.

\paragraph{The W-regrading $\fMe(4;\un|3)$}\label{Ex2} We set
\[
 \deg p=3, \; \deg q=1, \; \deg \xi=0, \; \deg \eta=4, \; \deg \zeta=2, \; \deg z=4 
\]
implying that $\gr:=\deg-4$ is the W-grading in the regraded algebra:
\[
\begin{tabular}{|c|c|c|c|c|}
\hline
$\fh_{-4}$ & $\fh_{-3}$ & $\fh_{-2}$ & $\fh_{-1}$ & $\fh_0$ \\
\hline\hline
$\xi$ & $q$ & $q^{(2)}+\xi\zeta$ & $p$ & $pq-\xi\eta$, $z$ \\
1 & & $\zeta$ &$p\xi+q\zeta$ & $\eta$, $F_\xi$\\
 \hline
\end{tabular}
\]

\paragraph{Digression: Irreducible $\fgl(1|1)$-modules}\label{AN81}
A.N.~Sergeev described all irreducible finite-dimensional representations over finite-dimensional solvable Lie
superalgebras over $\Cee$, see~\cite{Sg}.
We need the following modest part of his result, explicitly.
Let $\fg=\fgl(1|1)$ with a~basis
\be\label{glBasis}
x=\begin{pmatrix}
0 & 1\\
0 & 0\\
\end{pmatrix}, \; y= \begin{pmatrix}
 0 & 0\\
 1 & 0\\
 \end{pmatrix}, \; h=\begin{pmatrix}
 1 & 0\\
 0 & 1\\
 \end{pmatrix}, \; d=\begin{pmatrix}
 1 & 0\\
 0 & 0\\
 \end{pmatrix} 
\ee

Let $T$ be an irreducible representation of $\fg$ in the superspace $V$,
let $X,Y,H,D$ be the images of the corresponding elements~\eqref{glBasis} in this representation.
Since $h$ is even and central, the operator $H$ acts as multiplication by a~number $\lambda$.
If $\lambda=0$, the representation $T$ is 1-dimensional: $T=\mu\cdot\str$, up to the reversal of parity of~$V$;
the corresponding action in $\Pi(V)$ will be denoted by $\Pi(T)$.

Let $\lambda\ne 0$.
Since $[x,x]=0$, there exists a~non-zero $v\in V$ such that $X(v)=0$.
Set $w=Y(v)$. Then, $Y(w)=0$, since $[y,y]=0$. Then, 
\[
 X(w)=XY(v)=-YX(v)+H(v)=\lambda\cdot v. 
\]
Therefore, $V=\Span(v,w)$ and for every $\lambda\ne 0$, there are two irreducible $\fsl(1|1)$-modues:
both are of superdimension $1|1$ and differ only by the parity of~$v$.
Let $T_\lambda$ be the one for which $p(v)=\ev$. Then, the other one is $\Pi(T_\lambda)$.

The matrices of $X,Y,H$ in the basis $v,w$ are of the form:
\footnotesize 
\[
X=\begin{pmatrix}
0 & \lambda\\
0 & 0\\
\end{pmatrix}, \; Y= \begin{pmatrix}
 0 & 0\\
 1 & 0\\
 \end{pmatrix}, \; H=\begin{pmatrix}
 \lambda & 0\\
 0 & \lambda\\
 \end{pmatrix}
\]
\normalsize 

Clearly, the operator
\footnotesize $D=\begin{pmatrix}
 1 & 0\\
 0 & 0\\
\end{pmatrix} $ \normalsize
extends $T_\lambda$ to a~representation of $\fgl(1|1)$, we denote this representation also by $T_\lambda$.
Thus, all irreducible representations of $\fg$ are of the form:
\be\label{subscr}
\begin{cases}
 \text{$\mu\cdot\str$ or $\Pi(\mu\cdot\str)$}&\text{if $\lambda=0$}\\
 T_{\lambda,\mu}:=T_\lambda\otimes (\mu\cdot\str) \text{ and } \Pi(T_{\lambda,\mu})&\text{if } \lambda\neq0.
\end{cases}
\ee

\textbf{Let us return to the W-regrading $\fMe(4;\un|3)$}.
In order to describe the action of $\fh_0$ on the components $\fh_i$ for $i$ negative,
let us first describe the correspondence between the elements of $\fh_0$ in terms of the above basis of $\fgl(1|1)$:
\[
 x\lra \eta, \quad y\lra F_\xi, \quad h=[x,y]\lra \{\eta,F_\xi\}=-z+pq-\xi\eta, \quad d\lra -z.
\]

Then, in $\fh_{-1}$:
\[
\begin{array}{l}
\{\eta,p\}=0\Lra p\lra v, \ \ \{F_\xi,p\}=p\xi+q\zeta\Lra p\xi+q\zeta \lra w, \\\
\ad_{-z+pq-\xi\eta}|_{\fh_{-1}}=-\id, \ \{-z,p\}=p, \ \{-z,p\xi+q\zeta\}=0.
\end{array}
\]
So, $\fh_0$ acts on~$\fh_{-1}$ via~$T_{-1,0}$; recall, see formula~\eqref{subscr},
that the first subscript in~$T_{\lambda,\mu}$ indicates the weight under the action of the central element
of~$\fg_0$, which is the grading operator on~$\fg_0$. 

In component $\fh_{-2}$:
\[
\begin{array}{l}
 \{\eta,\zeta\}=0\Lra \zeta\lra v, \ \ \{F_\xi,\zeta\}=-q^{(2)}-\xi\zeta\Lra -q^{(2)}-\xi\zeta \lra w,
 \\
 \ad_{-z+pq-\xi\eta}|_{\fh_{-1}}=-2\id, \ \ \{-z,\zeta\}=\zeta, \ \ \{-z,-q^{(2)}-\xi\zeta\}=0.
\end{array}
\]
So, $\fh_0$ acts on $\fh_{-2}$ via $\Pi(T_{-2,0})$, see formula~\eqref{subscr}.

The action of $\fh_0$ on $\fg_{-3}$ is trivial, by zero: $[\fh_0,\fh_{-3}]=0$.

Finally, on $\fg_{-4}$ we get:
\[
 \begin{array}{l}
 \{\eta,1\}=0\Lra 1\lra v, \ \ \{F_\xi,1\}=2\xi=-\xi \Lra -\xi \lra w,
 \\
 \ad_{-z+pq-\xi\eta}|_{\fh_{-4}}=-4\id, \ \ \{-z,1\}=2\cdot 1, \ \ \{-z,-\xi\}=(-\xi),
 \end{array}
\]
i.e., $\fh_0$ acts on $\fh_{-4}$ via $T_{-4,1}$, see formula~\eqref{subscr}.

Let us summarize (the components and the $\fh_0$-action in them):
\[
\begin{tabular}{|c|c|c|c|c|}
\hline
$\fh_{-4}$ & $\fh_{-3}$ & $\fh_{-2}$ & $\fh_{-1}$ & $\fh_0$ \\
\hline\hline
$T_{-4,1}$ & $0\cdot \str$ & $\Pi(T_{-2,0})$ & $T_{-1,0}$ & $\fgl(1|1)$ \\
 \hline
\end{tabular}
\]
We see that this $\Zee$-grading is a~W-grading. We denote thus regraded superalgebra $\fMe(4;\un|3)$. 

If we interchange $p\lra q$ and $\xi\lra \eta$, we get an isomorphic $\Zee$-graded Lie superalgebra. 

\sssbegin{Theorem}\label{TmBy1}
The above three W-regrading of the Melikyan superalgebra $\fMe(3;\un|3)$ exhaust all its W-gradings.
The growth vectors $v_\fg$ of non-integrable distributions preserved by the W-regraded Melikyan superalgebra
$\fMe(3;\un|3)$ are as follows
\begin{equation}\label{GrVme33}
\renewcommand{\arraystretch}{1.4}
\begin{tabular}{|c|c|}
\hline $\fg$&$v_\fg$\\
\hline
$\fMe(3;\un|3)$&$(2|3, 3|3)C$\\
\hline
$\fMe(3;\un|4)$&$(0|3, 3|3, 3|4)$\\
\hline
$\fMe(4;\un|3)$&$(1|1, 2|2, 3|2, 4|3)$\\
\hline
\end{tabular}
\end{equation}
\end{Theorem}

\begin{proof}[Proof\nopoint] follows from the arguments in Subsection~\ref{sss8.1.1}.
\end{proof}

\ssec{$\fg:=\fBj(3;(n, 1,1)|3)$}\label{Bj33}
Let $\fBj(3;(n, 1,1)|3)$ be a~super analog of $3\fme(\un)$.
Consider the Lie superalgebra $\fosp(2|4)$, see~\cite{BGL}, 
\[
\text{with Cartan matrix \footnotesize $\begin{pmatrix}
 0&1&0 \\
 -1&2&-2 \\
 0&-1&2
 \end{pmatrix}$ \normalsize and grading $(010)$}
\]
and consider the Lie superalgebra $\fk(3;\un|2)$ realized by generating functions in $p,q,\xi,\eta, t$
in the standard grading (we do not indicate $\un$ in what follows):
\[
 \deg t=2, \; \deg p=\deg q=\deg \xi=\deg\eta=1, \quad \gr(f)=\deg(f) -2.
\]

Then, for a~basis of the non-positive components of $\fk(3|2)$ we take
\[
\begin{tabular}{|c|c|c|}
\hline
$\fk(3|2)_{-2}$ & $\fk(3|2)_{-1}$ & $\fk(3|2)_0\simeq \fosp(2|2)\oplus \Kee$ \\
\hline
$1$ & $p,\ q,\ \xi, \ \eta$ & $t, \ p^{(2)},\ pq, \ q^{(2)}, \ \xi\eta, \ p\xi,\ q\xi,\ p\eta, \ q\eta $\\
\hline
\end{tabular}
\]

If $p>3$, the component $\fk(3|2)_1$ is the direct sum of two irreducible $\fk(3|2)_0$-modules:
\[
V=\Span(tp,\ tq,\ t\xi,\ t\eta)\simeq (\fg_{-1})^* \text{
and~~} W=S^3(<p,\ q,\ \xi,\ \eta>).
\]

The Cartan prolongation of the subspace $\fk(3|2)_{\le 0}\oplus V$ adds to it a~1-dimensional
space $\Kee t^{(2)}$, and the prolong is isomorphic to $\fosp(2|4)$.
The Cartan prolongation of $\fk(3|2)_{\le 0}\oplus W$ yields $\fpo(2|2)$. 

If $p=3$, then $\fk(3|2)_0$-module $\fk(3|2)_1$ is more complicated. First, $W$ contains a~submodule 
\[
 \text{$W_0$ generated by $f$ such that $\deg_p f<3$ and $\deg_q f<3$.}
\]
The Cartan prolongation of the subspace $\fk(3|2)_{\le 0}\oplus W_0$ yields $\fpo(2; \One|2)$
whose shearing vector has the minimal values of heights of~$p$ and~$q$.
Adding to the first component of the ``initial part" to be prolonged any element of $W\setminus W_0$ increases
the height of the corresponding indeterminates.

But the most interesting is the fact that if $p=3$, then $W_0$ splits into the direct sum of irreducible
$\fk(3|2)_0$-modules, call them $W_\xi$ and $W_\eta$:
\[
\begin{array}{l}
 W_\xi=\Span(
 p^{(2)}q - p\xi\eta,\ pq^{(2)}+ q\xi\eta,\ p^{(2)}\xi ,\ pq\xi,\ q^{(2)}\xi
 ), \\
 W_\eta=\Span(
 pq^{(2)}- q\xi\eta,\ p^{(2)}q + p\xi\eta,\ p^{(2)}\eta,\ pq\eta,\ q^{(2)}\eta
 ).
\end{array}
\]

If we add any element of $W\setminus W_0$ to any of the submodule $W_\xi$ or $W_\eta$ we generate the whole
$\fk(3|2)_0$-module $W_0$. 

The Cartan prolongs of the subspaces $\fk(3|2)_{\le 0}\oplus V\oplus W_\xi$, and $\fk(3|2)_{\le 0}\oplus
V\oplus W_\eta$, are infinite-dimensional (if $\un$ is unconstrained) simple Lie superalgebras.
These prolongs are isomorphic.
We denote this Lie superalgebra $\fBj:=\fBj(3;\un|2)$, where the heights of all indeterminates, except $t$,
are constrained as we will see. 

Consider $\fBj$ as the Cartan prolong of the subspace $\fk(3|2)_{\le 0}\oplus V\oplus W_\eta$ and list all its
$W$-gradings.
As earlier, we consider the regradings obtained by assigning degrees to the indeterminates $p,\ q,\ \xi,\
\eta,\ t$. Let $\fg$ be a~$W$-regrading of $\fBj$ thus obtained. 

First of all, observe that $t, pq, \xi\eta$ must remain in $\fg_0$. Hence (recall formula~\eqref{N}),
\[
 \deg t=\deg p+\deg q=\deg\xi+\deg\eta=N \text{~~and~~} \gr(f)=\deg(f)-N.
\]

Since there are no restrictions on the height of $t$, the function $1$ must remain in the negative part
$\fg_{<0}$ of the $\Zee$-grading.
Since $\{p,q\}=\{\xi,\eta\}=-1$, at least one element of each pair $(p,q)$ and $(\xi,\eta)$ also must remain
in $\fg_{<0}$. 

By Lemmas~\ref{L4.1} and~\ref{L3.3} the condition $\fBj_-\cap \fg_0\ne 0$ is necessary for $W$-grading.
We see that if $p\in \fg_0$, then $\deg p=N$, see formula~\eqref{N}, so $\deg q=0$, and vice versa.
The same is true for the pair $(\xi,\eta)$.
However, the indeterminates~$p$ and~$q$ enter $\fBj$ on equal footing,
and hence it suffices to consider one of the two possible cases,
the indeterminates $\xi$ and $\eta$ are not on equal footing, so we have to consider both cases.

\sssec{Remark}\label{remBj}
The isomorphic partial Cartan prolongs of either of the subspaces $\fk(3|2)_{\le 0}\oplus W_\xi$ and
$\fk(3|2)_{\le 0}\oplus W_\eta$ are not~simple Lie superalgebras: they have $9|9$-dimensional ideals
$\fg_{-1}\oplus\fg_{0}\oplus\fg_{1}$; observe that the 2nd components of the prolongs vanish. 

\sssec{The W-regrading $\fBj(3;\un|3)$}\label{Ex1Bj}
We set 
\[
 \deg \xi=0, \; \deg t=\deg \eta=2, \deg p=\deg q=1, \quad \gr(f)=\deg(f)-2.
\]
Then,
\[
\begin{array}{l}
\fg_{-2}=\Span(1, \ \xi)=\Lambda(\xi), \\
\fg_{-1}=\Span(p,\ q,\ p\xi,\ q\xi)=\Span(p,q)\otimes \Lambda(\xi), \\
\fg_0=\Span(p^{(2)}, \ pq, \ q^{(2)},\ t,\ t\xi,\ \eta,\ \xi\eta), \\
\fg_1=\Span(tp,\ tq,\ pq^{(2)}+q\xi\eta,\ p^{(2)}q-p\xi\eta, \ p\eta,q\eta, \ tp\xi-p^{(2)}q\xi, \ tq\xi+pq^{(2)}\xi).
\end{array}
\]

Observe that $\Span(p^{(2)}, \ pq, \ q^{(2)})\simeq \fsl(2)$, and the following elements act on $\fg_{-1}$ as 
\[
\begin{array}{l}
\eta \mapsto-\der_\xi, \ \ \xi\eta\mapsto -\xi\der_\xi\\
t+\xi\eta\mapsto \text{multiplication by $-1$},\\
t\xi\mapsto \text{multiplication by $-\xi$}.
\end{array}
\]
Hence, $\Span(t,\ t\xi,\ \eta,\ \xi\eta)=\Lambda(\xi)\inplus \fvect(\xi)\simeq \fgl(1|1)$, where $\fm\inplus
\fa$ denotes the semidirect sum of an algebra $\fa$ and an $\fa$-module $\fm$ which is an ideal in the sum.

In notation of representations of $\fgl(1|1)$, see Subsection~\ref{AN81}, we get the following description of
the negative components as $\fg_0$-modules, recall formula~\eqref{subscr}: 
\[
\fg_0=\fsl(2)\oplus \fgl(1|1), \quad \fg_{-1}=\id_{\fsl(2)}\boxtimes T_{-1,-1}, \quad \fg_{-2}=T_{-2,-1}.
\]

We see that the regrading obtained is a~$W$-grading.
This is the form the Lie superalgebra $\fBj$ appeared in~\cite{BGL}, where it was denoted $\fBj(3;\un|3)$;
we will retain this notation. 

\sssec{The W-regrading $\widetilde \fBj(3;\un|3)$}\label{Ex2Bj}
We set 
\[
\deg \eta=0, \ \deg t=\deg \xi=2, \ \deg p=\deg q=1, \quad \gr(f)=\deg(f)-2.
\]

Then,
\[
\begin{array}{l}
\fg_{-2}=\Span(1,\ \eta)=\Lambda(\eta), \\ 
\fg_{-1}=\Span(p,\ q,\ p\eta,\ q\eta)=\Span(p,\ q)\otimes \Lambda(\eta), \\
\fg_0=\Span(p^{(2)},\ pq,\ q^{(2)},\ p^{(2)}\eta, \ pq\eta, \ q^{(2)}\eta,\ t, \ t\eta,\ \xi,\ \xi\eta).
\end{array}
\]

Thus, the part $\fg_{\le 0}$ coincides with the non-positive part of $\fk(3|2)$ in this grading. Their first
components are, however, different:
\[
\fg_1=\Span(tp, \ tq,\ pq^{(2)}+q\xi\eta, \ p^{(2)}q-p\xi\eta, \ p\xi, \ q\xi,\ tp\eta+p^{(2)}q\eta, \ tq\eta-pq^{(2)}\eta).
\]

In this case, $\fg$ is obtained as the Cartan prolong of $\fg_{\le 0}\oplus\fg_1$.
We denote it by $\widetilde\fBj(3;\un|3)$. 

\sssbegin{Theorem}\label{TmBy2}
The above examples exhaust all W-regradings of the Bouarroudj superalgebra $\fBj(3;\un|3)$.
The growth vectors $v_\fg$ of non-integrable distributions preserved by the W-regradings of $\fBj(3;\un|3)$
are as follows, hence both preserve a~contact distribution preserved by $\fk(3;\un|2;1)$
\begin{equation}\label{GrVbj33}
\renewcommand{\arraystretch}{1.4}
\begin{tabular}{|c|c|}
\hline $\fg$&$v_\fg$\\
\hline
$\fBj(3;\un|3)$&$(2|2, 3|3)C$\\
\hline
$\widetilde \fBj(3;\un|3)$&$(2|2, 3|3)C$\\
\hline
\end{tabular}
\end{equation}
\end{Theorem}

\begin{proof}
Let us consider the three other possible regrading
 
\textbf{Case 1.}
We set
\[
\deg p=0, \ \deg q=\deg t=2,\ \deg\xi=\deg\eta=1, \ \ \gr(f)=\deg(f)-2..
\]
Then, we get a~compatible $\Zee$-grading:
\[
\begin{array}{l}
\fg_{-2}=\Span(1,\ p,\ p^{(2)}), \\
 \fg_{-1}=\Span(\xi,\ p\xi, \ \eta,\ p\eta,\ p^{(2)}\eta),\\
 \fg_0=\Span(t,\ tp,\ tp^{(2)}-p^{(2)}\xi\eta, \ q,\ pq,\ p^{(2)}q-p\xi\eta, \ \xi\eta.)
\end{array}
\]

Observe that, on one hand, the element $\xi\eta$ is central in~$\fg_0$; on the other hand, it has two distinct
eigenvalues on~$\fg_{-1}$.
Hence, the $\fg_0$-action on $\fg_{-1}$ is definitely reducible and does not correspond to any $W$-grading.

\textbf{Case 2.}
We set
\[
\deg p=\deg\xi=0, \; \deg q=\deg\eta=\deg t=1, \quad \gr(f)=\deg(f)-1.
\]

Hence, $\fg_{-1}=\Span(1,\ p,\ p^{(2)}, \ \xi,\ p\xi)$.

To understand what~$\fg_0$ is we write its elements in the bigrading:
\[
\begin{tabular}{|c|c|c|c|}
\hline
$\fg_{0,-1}$ & $\fg_{0,0}$ & $\fg_{0,1}$ & $\fg_{0,2}$ \\
\hline\hline
$q$ & $t,\; \xi\eta$ & $tp, \; t\xi$ & $tp^{(2)}-p^{(2)}\xi\eta$\\
$\eta$ & $pq,\; p\eta,\; q\xi$ & $p^{(2)}q-p\xi\eta,\; p^{(2)}\eta$ & $tp\xi-p^{(2)}q\xi$\\
\hline
\end{tabular}
\]
This representation suggests, and the direct verification confirms, that $\fg_0=\fvect(1|1)\oplus \Kee\ z$,
where $z=t-pq+\xi\eta$ is the central element.
Since $q\in\fg_0$ acts on $\fg_{-1}$ as $\der_p$, while $\eta\in\fg_0$ acts as $(-\der_\xi)$, we can identify
$\fg_0$ with $\fvect(p,\xi)\oplus \Kee\ z$.
We denote this correspondence by~$D$.
Then, the $\fg_0$-action on $\fg_{-1}$ is given by the formula:
\[
D_f: \vf\mapsto D_f(\vf)+(\Div D_f)\vf, \quad z(\vf)=-2\vf=\vf \text{~~for any~~} f\in\fg_0, \; \vf\in\fg_{-1},
\]
i.e.,. $\fg_{-1}\simeq \Vol_0$ as $\fg_0$-module, see definition~\eqref{vol}.
We denote this Cartan prolong by $\fBj(3;\un|2)$.
Since its depth is equal to 1, it does not preserve any distribution.

\textbf{Case 3}.
We set
\[
\deg p=\deg\eta=0, \; \deg q=\deg\xi=\deg t=1, \quad \gr(f)=\deg(f)-1.
\]
Then, $\fg_{-1}=\Span(1,\ p,\ p^{(2)}, \ \eta,\ p\eta, \ p^{(2)}\eta)$.

Again, to understand what $\fg_0$ is we write its basis in the bigrading:
\[
\begin{tabular}{|c|c|c|c|c|}
\hline
$\fg_{0,-1}$ & $\fg_{0,0}$ & $\fg_{0,1}$ & $\fg_{0,2}$ &\\
\hline\hline
$q$ & $t,\; \xi\eta$ & $tp, \; t\eta$ & $tp^{(2)}-p^{(2)}\xi\eta$ & $tp^{(2)}\eta$\\
$\xi$ & $pq,\; q\eta,\; p\xi$ & $p^{(2)}q-p\xi\eta,\; pq\eta$ & $tp\eta+p^{(2)}q\eta$ &\\
\hline
\end{tabular}
\]
Direct calculations show that this Lie superalgebra is a~subsuperalgebra of the semi-direct sum 
$\fh=\cO(p;1|\eta)+\fvect(p;1|\eta)$ with the natural action on $\fg_{-1}=\cO(p;1|\eta)$.
The subsuperspace
\[\cO=\Span(t-pq-\xi\eta, \; tp+p^{(2)}q-p\xi\eta, \; t\eta-pq\eta, \; tp^{(2)}-p^{(2)}\xi\eta, \; tp\eta+p^{(2)}q\eta,\ tp^{(2)}\eta)
\]
is a~commutative ideal in $\fg_0$, naturally identified with $\cO(p;1|\eta)$ via the formula $f\mapsto \frac{\der f}{\der t}$.

The subsuperspace $S$ complimentary to this ideal is a~subalgebra $\fsl(2|1)\subset \fvect(p;1|\eta)$.
Let us write a~basis of $S$ in the bigrading:
\[
\begin{tabular}{|c|c|c|}
\hline
$S_{0,-1}$ & $S_{0,0}$ & $S_{0,1}$ \\
\hline\hline
$q$ & $t,\; pq-\xi\eta$ & $tp - p^{(2)}q+p\xi\eta$ \\
$\xi$ & $q\eta,\; p\xi$ & $ t\eta+pq\eta$ \\
\hline
\end{tabular}
\]
We denote this Cartan prolong by $1\fBj(3;\un|2)$.
Since its depth is equal to 1, this Cartan prolong does not preserve any distribution.
Lemmas~\ref{L4.1} and~\ref{L4.1.1} guarantee that we have described all $W$-gradings.

\end{proof}

\ssec{$\fg:=\fBj(1;\un|7)$}\label{Bj17}
The Lie superalgebra $\fBj(1;\un|7)$ is one more super analog of $3\fme(\un)$, see~\cite{BL}.
Indeed, while $3\fme(\un)$ is the Cartan prolong of the non-positive part of $p=5$ version of the Lie algebra
$\fg(2)$ over $\Cee$, the Lie superalgebra $\fBj(1;\un|7)$ is the Cartan prolong of the non-positive part of
$p=3$ version of the Lie superalgebra $\fh:=\fag(2)$ over $\Cee$ 
\[
\text{with Cartan matrix \footnotesize $\begin{pmatrix}
 0&1&0 \\
 -1&2&-3 \\
 0&-1&2
 \end{pmatrix}$ \normalsize and grading $(100)$}.
 \]
 Let
$\fh_{\leq 0}=\fag(2)_{\leq 0}=\mathop{\oplus}_{-2\leq i\leq 0}\fh_i$, where
$\dim\fh_{-2}=1$ and $\sdim \fh_{-1}=0|7$. Therefore, the Cartan prolong
$\fh_*$ of $\fh_{\leq 0}$ is a~subalgebra in $\fk(1; \un|7)$, where $\un$ is considered only if $p>0$, corresponding to the contact form \begin{equation}\label{alpha}
\alpha=dt+2udu-\sum\limits_{i=1,\ 2, \ 3}(v_idw_i+w_idv_i),
\end{equation}
where the $v_i$, $w_i$ and $u$
are odd and notations match~\cite[p.~354]{FH} while our $X_i^+$
and $X_i^-$, see eq.~\eqref{our}, correspond to $X_i$ and $Y_i$ of~\cite[p.~340]{FH}, respectively.
Recall that we denote the simple Lie algebra usually denote by $\fg_2$ by $\fg(2)$ in order not to confuse
with other of our notation, and observe that $\fh_0=\fg(2)\oplus\fc$, where 
$\fc=\Span(t)$ is the center of $\fh_0$. Hereafter the elements of
the Lie superalgebra of contact vector fields are given in terms
of their generating functions in indeterminates introduced in eq.~\eqref{alpha}.

We set 
\be\label{our}
\text{$X_3^\pm :=[X_1^\pm , X_2^\pm ]$, \ \ $X_4^\pm
:=[X_1^\pm , X_3^\pm ]$, \ \ $X_5^\pm :=[X_1^\pm , X_4^\pm ]$, \ \ $X_6^\pm :=[X_2^\pm , X_5^\pm ]$.}
\ee
The $\fh_0$-module $\fh_{-1}$ is irreducible with the highest weight vector $v_3$.
\[
\footnotesize
\renewcommand{\arraystretch}{1.4}
\begin{tabular}{|l|l|}
\hline
$\fh_{i}$&the generating functions of generators, as $\fh_{0}$-modules \\
\hline \hline
$\fh_{-2}$&$1$\\
\hline
$\fh_{-1}$& \text{$u$ and $v_i,
w_i$ for $i=1, 2, 3$}\\ 
\hline
$\fh_0$& $t$ \text{~~(spans $\fc$)}\\
&$X_1^+= v_3w_2+uv_1$\\ 
&$X_2^+=v_2w_1$\\ 
&$X_1^-=v_2w_3+uw_1$\\ 
&$X_2^-=v_1w_2$\\ 
\hline
\end{tabular}
\normalsize
\]
As expected, for $p=0$ and $p>3$, the Cartan prolong of $\fh_{\leq 0}$ is isomorphic to $\fag(2)$.

For $p=3$, the Lie subalgebra $\fg(2)\subset \fh_0$ is not simple, but has a~simple
Lie subalgebra isomorphic to $\fpsl(3) =\Span(X^\pm_1, X^\pm_3, X^\pm_4, H_1)$ generated by
$X_1^\pm$, and $X_3^\pm =[X_{1}^\pm, X_2^\pm]$. Let
$\widetilde\fg_0:=\fpsl(3)\oplus\fc$, where $\fc=\Span(t)$ 
is the center of
$\widetilde\fg_0$.

The $\widetilde\fg_0$-module $\widetilde\fg_{1}$ splits into two irreducible
components of dimensions $0|1$ and $0|7$.
For bases of $\widetilde\fg_0$ and $\widetilde\fg_{1}$ we take the following elements
\small
\be\label{notBj17}
\renewcommand{\arraystretch}{1.4}
\begin{array}{ll}
\begin{tabular}{|l|l|}
\hline
$\widetilde\fg_0$ & $t$ \\
{}& $ u v_1+v_3 w_2$ \\
{}& $ u w_1+v_2 w_3$ \\
{}& $ 2 u v_2+v_3 w_1$ \\
{}& $ v_1 w_1+v_2 w_2+2 v_3 w_3$ \\
{}& $ u w_2+2 v_1 w_3$ \\
{}& $ u v_3+v_1 v_2$ \\
{}& $ u w_3+2 w_1 w_2$ \\
\hline
\end{tabular} &\begin{tabular}{|l|l|}
\hline
$\widetilde\fg_1$ & $ 2 u v_1 w_1+2 u v_2 w_2+u v_3 w_3+v_3 w_1w_2+v_1 v_2 w_3$ \\
{} & $ t v_1+2 u v_3 w_2+v_1 v_2 w_2+2 v_1 v_3 w_3$ \\
{} & $2 t v_2+2 u v_3 w_1+v_1 v_2 w_1+v_2 v_3 w_3$ \\
{} & $ t v_3+2 u v_1 v_2+v_2 v_3 w_2+v_1 v_3 w_1$ \\
{} & $ t u+2 u v_1 w_1+2 u v_2 w_2+u v_3 w_3+2 v_1 v_2 w_3$ \\
{} & $ t w_1+2 u v_2 w_3+v_2 w_1w_2 +2 v_3 w_1w_3 $ \\
{} & $ t w_2+u v_1 w_3+2 v_1 w_1 w_2+2 v_3 w_2w_3 $ \\
{} & $ t w_3+u w_1 w_2+v_2 w_2 w_3+v_1 w_1 w_3$ \\
\hline
\end{tabular} 
\end{array} 
\ee

\normalsize

\sssec{On W-regradings of $\fg:=\fBj(1;\un|7)$}\label{WgradBj17}
We work in notation of table~\eqref{notBj17}. 

First of all, we observe that in any W-grading $1\in \fg_-$ since the height of $t$ is unbounded.
Then, the commutation relations $\{u,u\}=\{v_i,w_i\}=-1$ imply that $u\in\fg_-$, and at least one of the
elements of each pair $(v_i, w_i)$ must also lie in $\fg_-$.

Let
\[
\deg u=k, \quad \deg v_i=m_i, \quad \deg w_i=n_i.
\]
Then, the arguments analogous to those given in Subsection~\ref{sss8.1.1} yield a~system of homogeneous linear
equations for the degrees $k,\ m_i,\ n_i$:
\[
k= m_1+ n_1= m_2+ n_2= m_3+n_3 (=N, \text{~~see eq.~\eqref{N}}),
\]
which mean that $t$ and $v_1w_1+v_2w_2+2v_3w_3$ must remain in degree-0 component.
The next 6 equations mean that the remaining elements in $\widetilde\fg_0$, see Table~\eqref{notBj17},
must remain homogeneous (since all monomials-summands in them are linearly independent whereas we want to get
a~bigrading).
For example, just under $t$ we see $uv_1+v_3w_2$.
The degree of the first summand is equal to $k+m_1$, that of the second summand is $m_3+n_2$.
The equality of these degrees is the first equation in the left column below.
The other equations are similarly obtained. 
\[
\begin{tabular}{l}
$k+m_1=m_3+n_2, \quad \quad k+n_1=m_2+n_3,$\\
$k+m_2=m_3+n_1, \quad \quad k+n_2=m_1+n_3,$\\
$k+m_3=m_1+m_2, \quad \quad k+n_3=n_1+n_2$.\\
\end{tabular}
\]
The corank of this system is equal to 2, and therefore it has two linearly independent solutions one of which is 
\[
k=2, \ \ m_i=n_i=1 \text{~~for $i=1,2,3$}. 
\]
This means that the condition on the degree of a~given indeterminate to be equal to 0 whereas $\deg u<N$,
see~\eqref{N}, together with the requirement ``the grading is incompressible" uniquely determine the degrees
of the indeterminates. Recall that $\gr(f)=\deg(f)-N$. 

\sssbegin{Theorem}\label{TmBj17}
The following examples exhaust all W-gradings of the Bouarroudj superalgebra $\fBj(1; \un|7)$.
The growth vectors $v_\fg$ of non-integrable distributions preserved by the W-regradings of the Bouarroudj
superalgebra $\fBj(1; \un|7)$ are as follows.
\begin{equation}\label{GrVbj17}
\renewcommand{\arraystretch}{1.4}
\begin{tabular}{|c|c|}
\hline $\fg$&$v_\fg$\\
\hline
$\fBj(1;\un|7)$&$(0|7, 1|7)C$\\
\hline
$\fBj(4;\un|5)$&$(2|2, 4|4, 4|5)$\\
\hline
\end{tabular}
\end{equation}
\end{Theorem}
\begin{proof}[Proof\nopoint] follows from arguments in Section~\ref{WgradBj17} and the description of
 Examples.
\end{proof}

\paragraph{Example 1}\label{Ex1Bj17} Let
\[
\begin{tabular}{|c|c|c|c|c|c|c|c|c|}
\hline
$f$ & $t$ & $u$ & $v_1$ & $v_2$& $v_3$& $w_1$& $w_2$& $w_3$\\
\hline
$\deg f$ & 4 & 2& 3 & 3 & 4 & 1 & 1 & 0 \\
\hline
\end{tabular}
\]

Then, $\gr(f)=\deg(f)-4$ and
\[
\begin{tabular}{|c|c|c|c|c|}
\hline
$\fg_{-4}$ & $\fg_{-3}$ & $\fg_{-2}$ & $\fg_{-1}$ & $\fg_0$ \\
\hline
$ 1$ & $w_1$ & $u$ & $V_1=\Span(v_1, uw_2-v_1w_3) $ & $t$, $v_1w_1+v_2w_2-v_3w_3-t$\\
$w_3$ & $w_2$ & $uw_3-w_1w_2$ & $V_2=\Span(v_2, uw_1+v_2w_3)$ & $v_3, tw_3+uw_1w_2+v_1w_1w_3+v_2w_2w_3$\\
\hline
\end{tabular}
\]

Therefore, $\fg_0\simeq \fgl(1|1)$, and the space $\fg_{-1}$ splits into a~sum of the two $\fg_0$-invariant
subsuperspaces: $V_1$ and~$V_2$.
Hence, this grading is not a~Weisfeiler grading.

We have considered a~regrading for which the indeterminate $v_3$ moves to the degree-0 component.
Since the indeterminates $v_3$ and $w_3$ enter $\fBj(1;\un|7)$ absolutely symmetrically, the regrading for
which $w_3\in\fg_0$ is isomorphic to the above one.

\paragraph{Example 2}\label{Ex2Bj17}
Now, let
\[
\begin{tabular}{|c|c|c|c|c|c|c|c|c|}
\hline
$f$ & $t$ & $u$ & $v_1$ & $v_2$& $v_3$& $w_1$& $w_2$& $w_3$\\
\hline
$\deg f$ & 2 & 1& 2 & 2 & 3 & 0 & 0 & -1 \\
\hline
\end{tabular}
\]
Then, $\gr(f)=\deg(f)-2$ and for a~basis we take
\[
\begin{tabular}{|c|c|c|c|}
\hline
$\fg_{-3}$ & $\fg_{-2}$ & $\fg_{-1}$ & $\fg_0$ \\
\hline
&$ 1$, & $u,$ & $v_1, v_2$, $t$, $v_1w_1+v_2w_2-v_3w_3-t$\\
$w_3$ & $w_1$,$w_2$& $uw_1+v_2w_3, uw_2-v_1w_3$ & $tw_1-uv_2w_3+v_2w_1w_2-v_3w_1w_3$\\
 & $uw_3-w_1w_2$ &$tw_3+uw_1w_2+v_1w_1w_3+v_2w_2w_3$ & $tw_2+uv_1w_3-v_1w_1w_2-v_3w_2w_3$\\
\hline
\end{tabular}
\]

In this case, the component $\fg_0$ is a~semi-direct sum of the Heisenberg superalgebra and its 1-dimensional
outer derivation
\[
\fg_0=\mathfrak{hei}(0|4)\inplus \Kee\cdot t, 
\]
the representation $\fg_0$ in $\fg_{-1}$ is irreducible and $\fg_{-1}$ generates the whole $\fg_-$.
Therefore, this is a~$W$-regrading.

Since the pairs $(v_1,v_2)$ and $(w_1,w_2)$ symmetrically enter the Lie superalgebra $\fBj(1;\un|7)$, then the
regrading for which $w_1,w_2\in\fg_0$ is isomorphic to the just considered.

Thus, $\fBj(1;\un|7)$ has only one $W$-regrading; denote it $\fBj(4;\un|5)$.

\section{When $p=2$}\label{Sp=2}

In this section, $\Kee$ is an algebraically closed field of characteristic 2.

There are more types of simple Lie algebras for $p=2$ than for $p\neq 2$.
In particular, in~\cite{BGLLS2}, there are described characteristic~2 versions of all simple vectorial Lie
superalgebras over $\Cee$, and some deformations of these versions; then, the desuperization
functor~$\textbf{F}$ forgets the squaring and turns a~given Lie superalgebra over a~field of characteristic 2
into a~Lie algebra.
This desuperization, under which a~constrained coordinate might become unconstrained, gave
us several new infinite-dimensional Lie algebras, serial as well as exceptional (not belonging to a~series).
By considering only finite coordinates of their shearing vector we got new finite-dimensional simple Lie
algebras.
For details and examples, see~\cite{BLLS2, BGLLS2}. 

Although the classification of simple Lie (super)algebras is most difficult for $p=2$, the simple Lie
superalgebras are, miraculously, classified \textbf{modulo non-existing at the moment and very tough
 classification problem of simple Lie algebras}.
Since the list of simple Lie superalgebras is implicit, we consider in this section only several examples
which at the moment look outstanding.
We recall the two methods which turn any given simple Lie algebra into a~simple Lie superalgebra.

\ssec{Growth vectors $v_\fg$ of distributions preserved by desuperized exceptional Lie
 superalgebras}\label{GrVDesupExc}
The following are the growth vectors of the desuperized exceptional simple vectorial Lie superalgebras,
compare with Table \eqref{grVex} and table in \cite[\S25.5]{BGLLS2}.
\begin{equation}\label{grVexF}\footnotesize
\renewcommand{\arraystretch}{1.3}
\begin{array}{lll}
\begin{tabular}{|c|c|}
\hline 
Lie algebra $\fg$&$v_\fg$\cr

\hline
\hline

$\textbf{F}\fv\fle(4|3; 1)$&$(8, 9)C$\cr
\hline

$\textbf{F}\fv\fle(4|3; K)$&$(6, 9)$\cr
\hline \hline

$\textbf{F}\fk\fle(9|6)$&$(14, 15)C$\cr
\hline

$\textbf{F}\fk\fle(9|6; 2)$&$(16, 20)$\cr
\hline

$\textbf{F}\fk\fle(9|6; K)$&$(10, 15)$\cr
\hline

$\textbf{F}\fk\fle(9|6; CK)$&$(12, 18, 20)$\cr
\hline
\end{tabular} &
\begin{tabular}{|c|c|}
\hline 
Lie algebra $\fg$&$v_\fg$\cr
\hline
\hline

$\textbf{F}\fk\fas$&$(6, 7)C$\cr
\hline

$\textbf{F}\fk\fas(; 1\xi)$&$(8, 10)$\cr
\hline \hline

$\textbf{F}\fm\fb(4|5)$&$(8, 9)C$\cr 
\hline

$\textbf{F}\fm\fb(4|5; 1)$&$(8, 11)$\cr
\hline

$\textbf{F}\fm\fb(4|5; K)$&$(6, 9, 11)$\cr
\hline 
\end{tabular}
&
\end{array}
\end{equation}

\ssec{Digression: on restricted Lie (super)algebras}\label{ssRestr}
To advance further, recall the notion of classical a restricted Lie algebra.
Let the ground field $\Kee$ be of characteristic $p>0$, and $\fg$ a~Lie algebra.
For every $x\in\fg$, the operator $(\ad_x)^{p}$ is a~derivation of~$\fg$.
The Lie algebra~$\fg$ is said to be \textit{restricted}
or \textit{having a~$p$-structure} ${}[p]:\fg\tto\fg$, $x\mapsto x^{[p]}$ if $(\ad_x)^{p}$ is an inner
derivation $\ad_{x^{[p]}}$ for every $x\in\fg$ so that
\begin{equation*} \label{restricted-3}
\begin{array}{l}
[x^{[p]}, y]=(\ad_x)^{p}(y)\quad \text{~for any~}x,y\in\fg,
\\
(ax)^{[p]}=a^px^{[p]}\quad \text{~for any~}a\in\Kee,~x\in\fg,
\\
(x+y)^{[p]}=x^{[p]}+y^{[p]}+\mathop{\sum}\limits_{1\leq i\leq
p-1}s_i(x, y) \quad\text{~for any~}x,y\in\fg,
\end{array}
\end{equation*}
where $is_i(x, y)$ is the coefficient of $\lambda^{i-1}$ in
$(\ad_{\lambda x+y})^{p-1}(x)$.

Given a~Lie superalgebra $\fg$ of characteristic $p>0$, let the Lie
algebra $\fg_\ev$ be restricted and
\begin{equation} \label{restr3}
[x^{[p]}, y]=(\ad_x)^{p}(y)\quad \text{~for any~}x\in\fg_\ev,~y\in\fg.
\end{equation}
This gives rise to the map 
\begin{equation*}\label{2p}
{}[2p]:\fg_\od\to\fg_\ev,~~~x\mapsto(x^2)^{[p]},
\end{equation*}
satisfying the condition
\begin{equation*}\label{2p2}
{}[x^{[2p]},y]=(\ad_x)^{2p}(y)\quad\text{~for any~}x\in\fg_\od,~y\in\fg.
\end{equation*}
The pair of maps $[p]$ and $[2p]$ is called a~$p$-\textit{structure}
(or, sometimes, a~$p|2p$-\textit{structure}) on~$\fg$, and $\fg$ is said to be \textit{restricted}.

Observe that if $p=2$, there are several more indigenous versions of restrictedness, see~\cite{BLLS2}, but
only the classical restrictedness is needed in this paper.

\ssec{The two methods to get simple Lie superalgebras from simple Lie algebras}\label{2methods}
In~\cite{BLLS2}, there are described two methods that turn a~given simple Lie algebra into a~simple Lie
superalgebra; any simple Lie superalgebra is obtained by one of these methods.

\sssec{Method 1: the queerification}\label{Method1}
We define the Lie superalgebra $\fq(\fg)$ --- queerified Lie algebra $\fg$ --- by setting
$\fq(\fg)_\ev=\fg$ and $\fq(\fg)_\od=\Pi(\fg)$; define the
multiplication involving the odd elements as follows:
\begin{equation}\label{q(g)}
{}[x,\Pi(y)]=\Pi([x,y]);\quad (\Pi(x))^2=x^{[2]}\quad\text{~for
any~~}x,y\in\fg.
\end{equation}
Clearly, if $\fg$ is restricted and $\fii\subset\fq(\fg)$ is an ideal,
then $\fii_\ev$ and $\Pi(\fii_\od)$ are ideals in $\fg$. So, if $\fg$ is
restricted and simple, then $\fq(\fg)$ is a~simple Lie superalgebra.
(Note that $\fg$ has to be simple as a~Lie algebra, not just as a
\textit{restricted} Lie algebra, i.e., $\fg$ is not allowed to have
\textbf{any} ideals, not only restricted ones.) A generalized
queerification is the following procedure producing as many
simple Lie superalgebras as there are simple Lie algebras.

Let the \textit{$1$-step
restricted closure} $\fg^{<1>}$ of the 
Lie algebra $\fg$ be
the minimal subalgebra of the restricted closure
$\overline{\fg}$ containing $\fg$ and all the elements $x^{[2]}$,
where $x\in\fg$. To any 
Lie algebra $\fg$ the \textit{generalized queerification} assigns the Lie superalgebra
\begin{equation*}\label{tildeQ}
\tilde \fq(\fg):=\fg^{<1>}\oplus \Pi(\fg)
\end{equation*}
with squaring given by $(\Pi(x))^2=x^{[2]}$ for any $x\in\fg$.
Obviously, for $\fg$ restricted, the generalized queerification coincides with the
queerification: $\tilde \fq(\fg)=\fq(\fg)$.

\sssec{Method $2$}\label{Method2}
This method for any simple Lie algebra $\fg$ with a~$\Zee/2$-grading $\gr$, constructs a~simple Lie
superalgebra $\fs(\fg,\gr)$ whose depth might be bigger than the depth of $\fg$.
Let $\fg=\fg_+\oplus\fg_-$ be a~Lie algebra with a~$\Zee/2$-grading $\gr$.
Let $(\fg,\gr)$ be the minimal Lie subalgebra of the restricted closure $\overline{\fg}$ containing $\fg$ and
all the elements $x^{[2]}$, where $x\in\fg_-$.
Clearly, there is a~single way to extend the grading $\gr$ from $\fg$ to $(\fg,\gr)$. 

Let $\fs(\fg,\gr)$ be
the Lie superalgebra structure (squaring) on the space $(\fg,\gr)$ given by
\[
x^2:=x^{[2]}\text{~~for any $x\in\fg_-$.}
\]

\sssec{Fact}(\cite[Prop.7.4]{BLLS2})\label{Fact}
If $\bar\fg $ is the restricted closure of $\fg$, then
$[\bar\fg , \bar\fg]\subset[\fg,\fg]$.
Hence, $[\bar\fg, \fg]\subset\fg$.
If $\fg$ is of depth 3, then $\fs(\fg,\gr3)_{-4}=\fs(\fg,\gr3)_{-5}=0$.
Also $[\fg_{-1}, (\fg_{-1})^{[2]}]=0$.
To obtain $\bar\fg$, we only have to add to $\fg$ the squares of basis elements of negative odd degrees
corresponding to the indeterminates whose heights are $> 1$.

\sssbegin{Theorem}\label{Statement} \textup{(\cite{BLLS2})}
For any simple Lie algebra~$\fg$, the Lie superalgebras $\tilde \fq(\fg)$ and $\fs(\fg,\gr)$ for
a~$\Zee/2$-grading $\gr$ of~$\fg$ are simple.
\end{Theorem}

Since there is no classification of simple Lie algebras in characteristic~2, we will consider superizations
only of several examples, selected \textit{ad hoc} from examples in~\cite{BLLS2}. In particular, we do not
consider finite-dimensional simple vectorial superalgebras found in~\cite{BLLS2, BGLLS2, CSS, KL} and
references cited therein.

\ssec{Growth vectors $v_\fg$ of distributions preserved by queerifications of desuperized exceptional Lie
 superalgebras}\label{GrVSupDesupExc}
Let 
\[
\begin{array}{l}
D_1=\dim \Span(X^{[2]}\in\fg_{-1}\mid X \text{~corresponds to $\un_i$ unconstrained}),\\
D_3=\dim \Span(Y^{[2]}\in\fg_{-3}\mid Y \text{~corresponds to $\un_j$ unconstrained})
\end{array}
\]
If $\dim \fg_i=a$, then $\sdim\fq(\fg))_i=a|a$ in the $\Zee$-grading induced from $\fg$, whereas the vectors
$v_{\tilde \fq(\fg)}$ are as follows thanks to table~\eqref{grVexF} and Fact~\ref{Fact}:
\[
\begin{array}{l}
v_{\fg}(a, b) \Lra v_{\tilde \fq(\fg)}=(a\mid a,\ \ b+D_1\mid b+D_1); \\ 
v_{\fg}(a, b, c) \Lra v_{\tilde \fq(\fg)}=(a\mid a, \ \ b+D_1\mid b+D_1, \ \ c+D_1+D_3\mid c+D_1+D_3).
\end{array}
\]

\ssec{Growth vectors $v_\fg$ of distributions preserved by superizations of desuperized exceptional Lie
 superalgebras}\label{GrVSupDesupExcSup}
Let the depth of a~W-grading $\gr$ of $\fg$ be equal to 1 or 2.
Then, the depth of the corresponding W-grading $\fs(\fg, \gr)$ is equal to 2 with a~larger $(-2)$nd component:
\[
\fs(\fg, \gr)_{-2}=(\fg_{-1})^{[2]}\oplus \fg_{-2}, \text{~~where $(\fg_{-1})^{[2]}:=\Span(x^{[2]}\mid
x\in\fg_{-1})$}. 
\]

If $(\fg_{-1})^{[2]}\neq 0$, then the corresponding distribution preserved by $\fs(\fg, \gr)$ is new: 
\be\label{GrSg}
\begin{tabular}{|c|c|}
\hline 
Lie superalgebra $\fs(\fg, \gr)$&$v_\fg$\cr
\hline
\hline
$\fs(\textbf{F}\fv\fle(4|3))$&$(0|7, 7|7)$\cr
\hline
$\fs(\textbf{F}\fvas(4|4))$&$(0|8,8|8)$\cr
\hline
$\fs(\textbf{F}\fk\fas(; 3\xi))$&$(0|8,8|8)$\cr
\hline
$\fs(\textbf{F}\fk\fas(; 3\eta))$&$(0|7, 7|7)$\cr
\hline

\end{tabular}
\ee

\sssec{The two superizations $\fs(\fg, \gr)$ of simple Lie algebras $\fg$ of depth 3}\label{2depth3}

We know two such exciting examples: $\fg=\textbf{F}(\fk\fs\fle(9;\un|11))$ and $\fg=\textbf{F}(\fm\fb(3;\un|8))$ whose components are explicitly
described as modules over $\fg_0$ in~\cite[Tables 25.3, 25.5]{BGLLS2} reproduced in Table~\eqref{table32}.
These Lie algebras $\fg$ have $(-3)$rd components in their $\Zee$-gradings $\gr3$ considered here modulo~2, so
their superizations corresponding to such gradings $\gr3$ are of depth 6, the squares of the basis elements of
degree $-3$ corresponding to the indeterminates whose heights are $>1$ span the trivial 
$\fg_0$-module which we denote by $\fs(\fg,\gr3)_{-6}:=(\fg_{-3})^{[2]}$, see Fact~\ref{Fact}. 

\sssbegin{Theorem}\label{TmS(g)}
The growth vectors $v_{\fs(\fg)}$ of non-integrable distributions preserved by the W-regradings of the Lie
superalgebras $\fs(\fg,\gr3)$ of depth $6$ are as follows:
\begin{equation}\label{GrVd=6}
\renewcommand{\arraystretch}{1.4}
\begin{tabular}{|l|l|}
\hline $\fs(\fg,\gr3)$&$v_{\fs(\fg,\gr3)}$\\
\hline
$\fs(\textbf{F}(\fk\fs\fle(9;\un|11)))$&$(0|12, 6|12, 8|12, 10|12)$\\
\hline
 $\fs(\textbf{F}(\fm\fb(3;\un|8)))$&$(0|6, 3|6, 5|6, 7|6)$\\
\hline

\end{tabular}
\end{equation}
\end{Theorem}


\section{Miscellanea}\label{Sequiv} 

\ssec{Frobenius criterion of integrability of distributions}\label{Frob}
The classical formulation of the integrability criterion is as follows:

\textit{A given distribution $\cD$ on a~smooth manifold is integrable if and only if the sections of $\cD$ form a~Lie algebra with respect to the commutator.}

Since $\cD$ is a~solution of the system of Pfaff equations $\alpha_j(Y)=0$ for $j= 0,\dots , k$, where the
$\alpha_j$ are 1-forms, an equivalent form of the integrability criterion is 

\textit{A given distribution $\cD$ is integrable if and only if the ideal $I \subset \Omega^{\bcdot}$, generated by the forms $\alpha_j$, is closed with respect to the exterior differential, i.e., $dI \subset I$}.

J.~Bernstein generalized this criterion from smooth manifolds to smooth supermanifolds,
see~\cite[D.5.3, p.159]{Lpet};
M.~Kuznetzov extended this criterion to the case of characteristic $p>0$, see~\cite{Ku};
his proof can be directly superized. 

\ssec{A definition of equivalence of distributions and an example}\label{Equiv}
Two distributions on $M$ are considered \textit{equivalent} if there is an invertible diffeomorphism of $M$
sending one distribution into the other one. 

Let
$\fg$ be the Lie (super)algebra preserving the distribution $\cD$ given by a~system of Pfaff equations 
\be\label{pff}
\text{$\alpha_j(Y)=0$ for $Y\in\fvect(m|n)$ and a~set of indexes $j\in J$.}
\ee
If the filtration of $\fg$ is Weisfeiler, then $\sdim \fg_{-1}$ is equal to the (super)rank of the
distribution, every basis of $\fg_{-1}$ consists of the fields $X_i$ commuting with all fields $Y$
satisfying~\eqref{pff}, and the algebra $\fg_-$ is isomorphic to $\fh_-$. 
Therefore,
\begin{equation}\label{eqviv} 
\begin{minipage}[l]{14cm}
the two given distributions are equivalent if and only if the negative parts of the 
W-graded Lie (super)algebras preserving these distributions are isomorphic.
\end{minipage}
\end{equation} 

The growth vector is a~rather rude invariant of the distribution. 
For example, consider the regradings~\eqref{newDeg} of Lie superalgebras of series $\fk$ and $\fm$ defined in
eq.~\eqref{simpEx}.
The negative parts of the regraded components are not isomorphic, although the growth vectors of the
corresponding distributions coincide.
Indeed, consider the simplest example illustrating the equality \eqref{obs}: 
\be\label{simpEx}
\begin{array}{ll}
\text{$\fk:=\fk(3|2)=\fk(t,p,q,\xi,\eta)$}& \text{and $\alpha =dt + pdq-qdp +\xi d\eta+ \eta d\xi$},\\
\text{$\fm:=\fm(2)=\fm(\tau, u,\sigma, v,\eps)$}& \text{and $\beta =d\tau +u d\sigma+ \sigma du+v d\eps+ \eps dv$}
\end{array}
\ee
with the grading given by the formulas
\be\label{newDeg}
\begin{array}{lll}
\text{in $\fk$:}& \text{$\det t=\deg\eta=2$,}& \text{$\deg \xi=0$ and $\deg p=\deg q=1$,}\\
\text{in $\fm$:}& \text{$\deg\tau=\deg u=2$,}& \text{$\deg \sigma=0$ and $\deg \eps=\deg v=1$.}
\end{array}
\ee
Then, in the grading~\eqref{newDeg}, we have
\[
\text{$\sdim\fk_{-1}=\sdim \fm_{-1}=(2|2)$, while $\sdim\fk_{-2}=\sdim \fm_{-2}=(1|1)$.}
\]
Then, 
\be\label{however}
\begin{array}{ll}
\text{$\fk_{-1}=\Span(p,q,p\xi,q\xi)$}& \text{and $[(\fk_{-1})_{\od}, (\fk_{-1})_{\od}]=0$,}\\
\text{$\fm_{-1}=\Span(\eps, v\sigma, v,\eps\sigma)$}& \text{and $[(\fm_{-1})_{\od},(\fm_{-1})_{\od}]\ne 0$.}
\end{array}
\ee


\ssec{The curvature tensor for non-holonomic distributions}\label{CurvTens}

Under the classical definition, the contact distribution has a~non-zero curvature. Nevertheless, any contact
distribution can be locally reduced to a~normal shape of the contact form defining the distribution
($dt-\sum p_idq_i$ as in analytical mechanics and over fields of characteristic $p=2$, see~\cite{BGLLS2}, or
$dt-\sum (p_idq_i-q_idp_i) $ as in other mathematical papers, especially in characteristic
$p> 2$), so in some sense the contact structure is ``flat" --- there exists only one equivalence class of
normal shapes of the contact form or, equivalently, distribution --- and the same is true 
for the two super versions of the contact structure, with the same proof as in~\cite{GL1}. 
Skryabin proved the flatness of the contact structure for $p>0$
(existence of only one equivalence class of normal shapes of contact forms), see~\cite{Sk1}. 

The contradiction between the classical definition of the curvature and the manifest flatness of the contact distributions disappears under the definition of
the curvature tensor that naturally takes non-holonomicity of the distribution into account (see~\cite{L,
 GL1}), the contact and pericontact distributions are flat, see \cite{BGLS}. This ``non-holomomic'' analog of
the curvature tensor is defined (despite Vershik's doubts in its existence, see Appendix 3 in~\cite{SV}) and
computed in~\cite{L} for the two cases of Theorem formulated on minipage~\eqref{VG}: in particular, the Engel
distribution is not flat.

The requirement~\eqref{assum} that all nilpotent Lie (super)algebras $\fg(*)$ should be isomorphic in
a~neighborhood of a~given point $*$ is an analog of flatness and vanishing of what physicists call Wess-Zumino
constraints, see~\cite{GL0, BGLS}.

In~\cite{GL1}, the analogs of the curvature tensor are computed for several distributions preserved by the Lie
algebra $\fsl(n)$ considered with various $\Zee$-gradings, and for the exceptional Lie algebra~$\fg(2)$.

In~\cite{GLS}, there are computed the analogs of the curvature tensor for the \textbf{super}manifolds
corresponding to those 11 of 15 exceptional simple W-graded vectorial Lie superalgebras over $\Cee$ which
preserve non-holonomic distributions.

\ssec{Open problems}\label{OP}
1) For each of the distributions described above for $p>0$, compute the corresponding ``non-holonomic'' analog
of the Riemann curvature tensor --- continuation of the results in~\cite{GLS, L, GL1, BGLS}.
The example~\eqref{simpEx} shows that the computations of~\cite{GLS} performed for the \textit{exceptional}
simple W-graded vectorial Lie superalgebras should be extended also to embrace the \textit{serial} W-graded Lie superalgebras for any $p$.
The first to consider are regradings of Lie superalgebras of series $\fk$ and $\fm$ from formula~\eqref{obs}. 

The $\fg_0$-invariant part of $H^2(\fg_-; \fg)$ has an interpretation formulated, perhaps, for the first time, in~\cite{CK1, CK1a}:
the elements of $(H^2(\fg_-; \fg))^{\fg_0}$ are obstructions to the filtered deformations of $\fg$ --- contributing to the classification of simple Lie (super)algebras.
Ladilova computed the (obstructions to) filtered deformations for certain W-gradings of Skryabin algebras~\cite{LaY, LaD, LaZ}.
Clearly, the spaces $H^2(\fg_-; \fg)$ and $H^2(\fh_-; \fh)$ for two distinct W-gradings of the same Lie (super)algebra have little in common, see~\cite{GLS}; it is an \textbf{Open problem} to compute both $H^2(\fg_-; \fg)$ and the invariants $(H^2(\fg_-; \fg))^{\fg_0}$ for all new W-gradings found in this paper.

2) For modular Lie algebras considered for $p>2$, the W-gradings with the smallest codimension of the maximal subalgebra --- the analogs of the standard gradings of the serial Lie algebras over~$\Cee$ --- are provided by
the normalizers of \textit{sandwich subalgebras}, see~\cite{Kir, KirS}. What are analogs of such constructions for $p=2$ and for Lie superalgebras for any~$p$?

\ssec{Another local invariant of non-integrable distributions}\label{ssVar}
Observe that the local invariant Varchenko introduced to distinguish non-equivalent distributions,
see~\cite{Va}, concerns a~situation totally different from the one we consider here. Given a~W-graded Lie
(super)algebra with negative part $\fg_-:=\oplus_{-d\leq i<0} \ \fg_i$ we consider only distributions of the
form $\cD:=\fg_{-1}$, whereas Varchenko allowed distributions $\cD:=\oplus_{-k\leq i<0}\ \fg_i$ for $k>1$, and
introduced a~local invariant of such distributions. For example, he considered a~3-dimensional distribution
$\cD$ such that $\cD+[\cD,\cD]=\Ree^5$ with \textbf{filtered} Lie algebras generated by $\cD$, not
$\Zee$-graded, as we consider, and the associated graded Lie algebra in Varchenko's example is generated not
just by $\fg_{-1}$.
In more details: let the basis sections of the tangent bundle be $X_1, \dots , X_5$,
with $\cD$ spanned by $X_1, X_2 , X_3$.
Varchenko considered the two cases: with commutation relations 
\be\label{cases}
[X_1, X_3]=X_4, \ \ [X_2, X_3]=X_5\ \ \text{(in both cases)}
\ee
and, in case 1, there are no more non-zero relations, whereas in case 2 there is one more:
\[
 [X_1, X_2]=X_3.
\]

In case 1, the distribution spanned by $X_1, X_2$ is integrable. In case 2, there are no such
subdistributions. (Actually, this property distinguishes these two cases without any additional invariants.) 
 
The symmetry algebras of these distributions constructed via Vershik-Gershkovich in these two cases are
isomorphic, described by relations~\eqref{cases} only.

\ssec{Remark}\label{Rems}
Apart from the description of the Lie (super)algebra in terms of Cartan prolongation 
one can, sometimes, describe the vectorial Lie (super)algebra by other means more convenient for description
of multiplication than in terms of Cartan prolongation we use 
in this work, cf.~\cite{Sk, Ku} and~\cite[arXiv version, p.28]{GL}.


\subsection*{Acknowledgements} We are thankful to Pavel Grozman for his wonderful code \textit{SuperLie}, see
\cite{Gr} that helped us a~lot.
A.K. was supported by the GA\v{C}R projects 24-10887S and by HORIZON-MSCA-2022-SE-01-01 CaLIGOLA.
This article is based upon work from COST Action CaLISTA CA21109 supported by COST (European
Cooperation in Science and Technology, \url{www.cost.eu}).

\def\eightit{\it}
\def\bib{\bf}
\bibliographystyle{amsalpha}

\end{document}